\documentclass[11pt,a4]{article}

\usepackage[numbers]{natbib}

\usepackage{subfigure}% Support for small, `sub' figures and tables
\RequirePackage{fix-cm}
 
\usepackage[a4paper]{geometry}

\usepackage{latexsym}
\usepackage{amsmath}
\usepackage{amssymb}
\usepackage{theorem}
\usepackage{xspace}
\usepackage{array}
\usepackage{subfigure}
\usepackage{epsfig}
\usepackage{hhline}
\usepackage{bbm}
\usepackage{color}
\usepackage{graphicx}
\usepackage{multirow}
\usepackage{dcolumn}
\usepackage{caption}

\usepackage{xcolor,colortbl}

\definecolor{Gray}{gray}{0.85}

\def\pth#1{\left(#1\right)}

\def\cro#1{\left[#1\right]}

\def\eeX{\mathbb{X}}

\newcommand{\eY}{{\bf Y}}
\newcommand{\eX}{{\bf X}}

\newcommand{\R}{\mathbb{R}}
\newcommand{\PP}{\mathbb{P}}

\def\argmin{\mathop{\mathrm{arg\,min}}}

\def\ee1{\textrm{\mathversion{bold}$\mathbf{\varepsilon}$\mathversion{normal}}}  
  
\def\eeps{\textrm{\mathversion{bold}$\mathbf{\varepsilon}$\mathversion{normal}}}

\def\eb{\textrm{\mathversion{bold}$\mathbf{\beta}$\mathversion{normal}}}

\def\emu{\textrm{\mathversion{bold}$\mathbf{\mu}$\mathversion{normal}}}  

\def\eE{I\!\!E}

\def\e1{1\!\!1}

\theoremstyle{plain}

\newtheorem{theo}{Theorem}[section]
\newtheorem{lem}{Lemma}[section]
\newtheorem{rem}{Remark}[section]

\definecolor{gray2}{rgb}{0.4,0.4,0.4}

\usepackage[title]{appendix}

\renewcommand{\eqref}[1]{(\ref{#1})}

\begin{document}

 \title{Change-point Detection by the  Quantile LASSO Method}

 \author{Gabriela CIUPERCA$^1$ and Mat\'u\v{s} MACIAK$^{2}$}
 
 \maketitle
 \setcounter{footnote}{1}
\footnotetext{\noindent \small{Universit\'e de Lyon, Université Lyon 1, CNRS, UMR 5208, Institut Camille Jordan, Bat.  Braconnier, 43, blvd du 11 novembre 1918, F - 69622 Villeurbanne Cedex, France\\
	\textit{Email address}: Gabriela.Ciuperca@univ-lyon1.fr\\
	\indent $^2$Charles University, Faculty of Mathematics and Physics, Department of Probability and Mathematical Statistics, Sokolovsk\'a 83, Prague, 186 75, Czech Republic\\
	\textit{Email address}: Matus.Maciak@mff.cuni.cz}}

% Abstract, keywords, and classification codes
%\startabstract{%

% MSC 2010 subject classification codes
 
\begin{abstract}

A simultaneous change-point detection and  estimation  in a piece-wise constant model is a common task in modern statistics. If, in addition, the whole estimation can be performed automatically, in just one single step without going  through any hypothesis tests for non-identifiable models, or  unwieldy classical a-posterior methods, it becomes an interesting, but also challenging idea.  In this paper we introduce the estimation method based on the quantile LASSO approach. Unlike standard LASSO approaches,  our  method does not rely on typical assumptions  usually required for the model errors, such as sub-Gaussian or Normal distribution. The proposed quantile LASSO method can effectively handle heavy-tailed random error distributions, and, in general, it offers a more complex view of the data as one can obtain any conditional quantile of the target distribution, not just  the conditional mean. It is proved that under some reasonable assumptions the number of change-points is not underestimated with probability tenting to one, and, in addition, when the number of change-points is estimated correctly, the change-point estimates provided by the quantile LASSO  are consistent. Numerical simulations are used to demonstrate these results and to illustrate the empirical performance robust favor of the proposed quantile LASSO method.
\end{abstract}
%\makechaptertitle
\textit{Keywords:} quantile LASSO; change-points; sparsity; piece-wise constant model; automatic detection; consistency.
% Email address for corresponding author
%\correspondingauthor[*]{\\\email{Insert your email address here only after your paper has been accepted}}

 %\end{frontmatter}

\section{Introduction}
Change-points in statistical models attract a lot of attention in recent years. The reason is 
that a continuous or even smooth favor of standard modeling approaches is not what we usually observe in real life situations. In many applications it is quite common that the mechanism producing  data can suddenly change. This usually happens due to some known or unknown event which caused this change. In such situations  we  refer to change-points and we are interested in their detection, estimation, and  statistical inference.

The change-point detection and estimation is  typically performed using a standard $L_{2}$-norm minimization, therefore the estimated structure can be interpreted as a conditional mean value of the target variable. There are many approaches proposed in the statistical literature to handle structural breaks (change-points respectively) from various perspectives (e.g., \cite{Antoch06, Csorgo1997, desmet2011, gao08, Horvath2002, pestova.pesta.16} to name a few). Such methods are either based on a segmentation principle (e.g, \cite{kim2009}) or a two stage approach (e.g., \cite{Csorgo1988, qiu_yandell}), where in both one firstly needs to detect potential change-point locations, and later, in the second phase---if there are some change-points detected---the overall dependence structure is estimated using the $L_{2}$-norm objective function and  the knowledge about the existing change-points gained in the first phase. An alternative idea was recently proposed in \cite{Frick.14} where the authors utilized a two stage non-convex minimization based on likelihood approach to recover piece-wise constant trend in exponential family models.

In order to avoid the two (and more) stage estimation techniques mentioned above, an effective algorithm can be obtained when taking an advantage of some recent  developments in the area of machine learning approaches and atomic pursuit techniques, the LASSO regularization in particular.

Although the pioneering idea of the LASSO penalization originates in sparse signal recovering problems (see  \cite{chen2001} and \cite{Tibshirani:96}) it can be also effectively used for the change-point detection and estimation. There is enormous literature available on LASSO in general (see \cite{Tibshirani:11} for a nice summary) with various LASSO modifications (e.g., fused LASSO proposed in \cite{tibs2005}; adaptive LASSO introduces in \cite{zou2006}; or elastic LASSO presented in \cite{zou2005}), which can be also used for the trend filtering (e.g., \cite{tibs2014}). On the other hand, there is only very little work available on the automatic change-point detection using the LASSO type methods. A simple change-point in location problem in a piece-wise constant model within the LASSO estimation framework was firstly considered in \cite{Harchaoui.Levy.10}, but an alternative insight on the same model can be also found in \cite{boysen2009} and \cite{mammen1997}.   A generalization of the piece-wise constant change-point model into a piece-wise linear and continuous case was considered in \cite{Maciak.Mizera.16} and a~more general linear scenarios are presented in \cite{Ciuperca.14}, \cite{Qian.Su.16}, and \cite{tibs2014}. Some post-selection inference tools in such models are discussed in \cite{Frick.14} and \cite{Tibshirani.16}. In addition, a high-dimensional regression scenario for detecting change-points by employing the LASSO penalty is  investigated in  \cite{Lee.Seo.Shin.16}. However,  in all the aforementioned situations the authors consider the standard $L_{2}$-norm based approach for estimating the conditional mean and, moreover, the results are derived under the assumptions on the Gaussian (or sub-Gaussian respectively) distribution of random errors. 

On the other hand, modeling the conditional mean may not be sufficient from the practical point of view. The reason is that there is only a  limited information provided about the target distribution when referring to its mean value. Ideally, one should be interested in estimating the whole conditional distribution which, unfortunately, turns out to be a quite complex problem. Instead, the quantile LASSO approach  allows us to estimate any conditional quantile and therefore, we can still obtain a complex and overall insight into the distribution of the data. The main idea presented in this paper follows as a generalization of the approach presented in  \cite{Harchaoui.Levy.10}  and further elaborated in \cite{lin2016}. We consider the same model, however, with one key difference: the authors in both aforementioned papers work either with the normally distributed random error terms or the zero mean errors with a sub-Gaussian distribution. Unlike their work, the results derived in this paper are free of such distributional assumptions imposed on the random error terms. We utilize the LASSO regularized estimation approach together with the standard check function $\rho_{\tau}(v) = v(\tau - \e1(v < 0))$, for $v \in \mathbb{R}$, and $\tau \in (0,1)$ (see \cite{Koenker.05}), which allows us to work with various error distributions accounting also for  random error terms with outliers  or heavy-tailed distributions with no direct specification on their moments. 

A posteriori detection of the change-points (their number and locations) by the quantile LASSO model was already considered in \cite{Ciuperca-16}, but it is done by a rather unwieldy technique to put into practice: in order to find the number of change-points one firstly needs to minimize a Schwarz-type criterion, locate the change-points, and estimate the model between two consecutive change-points. Moreover, the approach presented in \cite{Ciuperca-16} does not cover the piece-wise constant model due to the non-singularity of the design matrix. Therefore, the method presented in this paper has the advantage of overcoming this issue and, in addition, it simultaneously estimates the number of change-points, detects their locations, and recovers the overall quantile structure in a robust manner.

Considering the quantile LASSO estimator we mainly focus on providing some precision for the performance of the change-point location detection, similarly as in \cite{Harchaoui.Levy.10}, rather than proving the consistency in terms of the sign consistency or the oracle properties as, for instance, considered in \cite{Meinshausen2006} or \cite{zhao2006}. This allows us to use less strict assumptions for the design matrix, which has a very specific form in our case, and, otherwise, does not satisfy stricter irrepresentable conditions, or the eigenvalue restriction required for the sign consistency or the mean consistency in the standard $L_{2}$-norm sense (see \cite{Harchaoui.Levy.10, Meinshausen2006} or \cite{zhao2006} for more details). 

The main contribution of this paper lies in the new robust quantile LASSO proposal for a simultaneous change-point detection and estimation: this method is free of any restrictive distributional assumptions common for the standard LASSO approach which is possible due to the different loss function employed in the minimization problem, analogously to \cite{he1996} or \cite{Maciak.Huskova.17}. Moreover, the estimation method presented in this paper is proved to be consistent with respect to the change-point detection and estimation and the consistency results do not depend on such strict assumptions as one needs to require for the sign consistency or the oracle properties. Therefore, the modeling framework presented in this paper is much widely applicable in practical situations and the final model can be easily obtained by using common estimating techniques and standard optimization toolboxes. 
 
This paper is organized as follows: in the next section we introduce the quantile LASSO model and we propose the estimation approach for fitting the model. The main theoretical results are presented in Section \ref{results}  and  the empirical performance is investigated via an extensive simulation study in  Section \ref{simulations}. Some remarks and comments are given in Section \ref{Conclusion}.   All proofs of the theorems are given in the appendix section.

\section{Model and  Notations}
\label{section:motivation}
Let us consider a sample $Y_{1}, \dots, Y_{n}$, for  $n \in \mathbb{N}$, with a specific location structure  with $K^* \in \mathbb{N}$ change-points, located in $t_{1}^{*}  \dots t_{K^*}^* \in \{1, \dots, n\}$,  such that $1 < t_{1}^* < t_{2}^* < \dots < t_{K^*}^* < n$, and
\begin{equation}\label{eq1}
Y_t=\mu^*_k +\varepsilon_t, \quad \textrm{for} \quad t =1, \cdots, n, \quad k=1, \cdots , K^*+1, \quad t^*_{k-1} \leq t \leq t^*_{k} -1,
\end{equation}
where $t_{0}^* = 1$ and $t_{K^* + 1}^* = n + 1$. The model can be equivalently  expressed as 
\begin{equation}
\label{eq1_alt}
Y_t= \sum_{k = 1}^{K^* + 1} \mu^*_k \e1_{\{ t_{k - 1}^* \leq t \leq t_{k}^* - 1\}} + \varepsilon_t, \quad \textrm{for} \quad t = 1, \dots, n,
\end{equation}
with $K^* + 1$ unknown parameters (phases) to be estimated and the corresponding change-point locations $t_1^*, \dots, t_{K^*}^*$, which are also left unknown. Alternatively, we can also use the formulation 
\begin{equation}
\label{eq2}
Y_t=u^*_t +\varepsilon_t, \qquad \textrm{for} \quad t =1, \cdots, n,
\end{equation}
where $u^*_t=\mu^*_k$, for $t=t^*_{k-1}, \cdots , t^*_k-1$, and $k = 1, \dots, K^* + 1$ (see Figure 1 for an illustration). The random error terms $\{\varepsilon_{t}\}_{t= 1}^{n}$ are assumed to be 
independent and identically distributed random variables with some (unknown) continuous distribution function $F$.
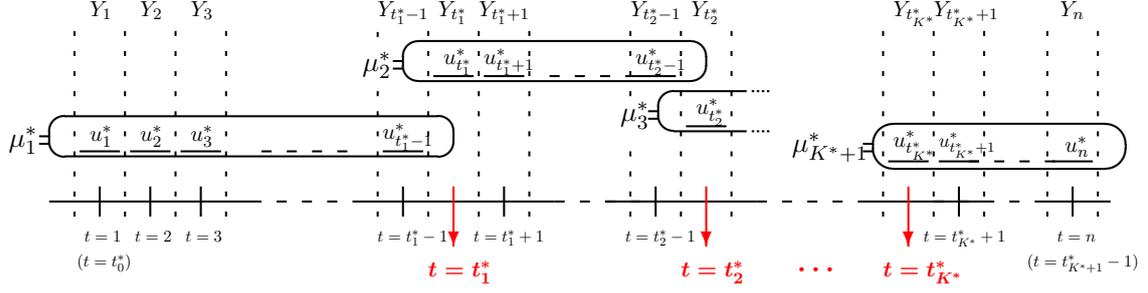
\begin{figure}
	\label{fig:illustration}
	{\centering \scalebox{0.95}{
		\setlength{\unitlength}{.7cm}
		\begin{picture}(25,6)(-5.5,-1.5)
		\thicklines
		\put(-5,0){\line(1,0){4}}
		\multiput(-4.5,-0.3)(0,0.4){10}{\line(0,1){0.1}}
		\put(-4,-0.3){\line(0,1){0.6}}
		\put(-4.4,1){\line(1,0){0.8}}
		\put(-4.3,-0.8){\scalebox{0.6}{$t = 1$}}
		\put(-4.45,-1.3){\scalebox{0.6}{$(t = t_0^*)$}}
		\put(-4.2,1.2){\scalebox{0.8}{$u_{1}^*$}}
		\multiput(-3.5,-0.3)(0,0.4){10}{\line(0,1){0.1}}	
		\put(-3,-0.3){\line(0,1){0.6}}
		\put(-3.4,1){\line(1,0){0.8}}
		\put(-3.3,-0.8){\scalebox{0.6}{$t = 2$}}
		\put(-3.2,1.2){\scalebox{0.8}{$u_{2}^*$}}
		\multiput(-2.5,-0.3)(0,0.4){10}{\line(0,1){0.1}}	
		\put(-2,-0.3){\line(0,1){0.6}}
		\put(-2.4,1){\line(1,0){0.8}}
		\put(-2.3,-0.8){\scalebox{0.6}{$t = 3$}}
		\put(-2.2,1.2){\scalebox{0.8}{$u_{3}^*$}}
		\multiput(-1.5,-0.3)(0,0.4){10}{\line(0,1){0.1}}
		
		\multiput(-1,0)(0.5,0){4}{\line(1,0){0.2}}
		\multiput(-0.8,1)(0.5,0){4}{\line(1,0){0.2}}
		\put(-1,1.3){\oval(8.0,0.8)}
		\put(-5.8,1.1){$\mu_{1}^*$}
		\put(-5.2,1.15){\line(1,0){0.2}}
		\put(-5.2,1.30){\line(1,0){0.2}}
		
		\put(1,0){\line(1,0){4}}	
		\multiput(1.5,-0.3)(0,0.4){7}{\line(0,1){0.1}}
		\multiput(1.5,3.28)(0,0.4){1}{\line(0,1){0.1}}
		\put(2,-0.3){\line(0,1){0.6}}
		\put(1.6,1){\line(1,0){0.8}}
		\put(1.5,-0.8){\scalebox{0.6}{$t = t_1^* - 1$}}
		\put(1.6,1.25){\scalebox{0.8}{$u_{t_{1}^* -1}^*$}}
		\multiput(2.5,-0.3)(0,0.4){10}{\line(0,1){0.1}}
		\put(3,0.3){\textcolor{red}{\vector(0,-1){1.2}}}
		\put(2.6,2.5){\line(1,0){0.8}}
		\put(2.5,-1.5){\textcolor{red}{\scalebox{0.8}{$\boldsymbol{t = t_1^*}$}}}
		\put(2.8,2.75){\scalebox{0.8}{$u_{t_{1}^*}^*$}}
		\multiput(3.5,-0.3)(0,0.4){10}{\line(0,1){0.1}}
		\put(4,0.3){\line(0,-1){0.6}}
		\put(3.6,2.5){\line(1,0){0.8}}
		\put(3.4,-0.8){\scalebox{0.6}{$t = t_1^* + 1$}}
		\put(3.6,2.75){\scalebox{0.8}{$u_{t_{1}^* + 1}^*$}}
		\multiput(4.5,-0.3)(0,0.4){10}{\line(0,1){0.1}}	
		
		\multiput(5,0)(0.5,0){4}{\line(1,0){0.2}}
		\multiput(4.9,2.5)(0.5,0){4}{\line(1,0){0.2}}	
		\put(5,2.8){\oval(6.0,0.8)}
		\put(1.2,2.6){$\mu_{2}^*$}
		\put(1.8,2.65){\line(1,0){0.2}}
		\put(1.8,2.80){\line(1,0){0.2}}
		
		\put(6,0){\line(1,0){3}}
		\multiput(6.5,-0.3)(0,0.4){4}{\line(0,1){0.1}}
		\multiput(6.5,2.1)(0,0.4){4}{\line(0,1){0.1}}
		\put(7,-0.3){\line(0,1){0.6}}		
		\put(6.6,2.5){\line(1,0){0.8}}
		\put(6.4,-0.8){\scalebox{0.6}{$t = t_2^* - 1$}}
		\put(6.6,2.75){\scalebox{0.8}{$u_{t_{2}^* -1}^*$}}
		\multiput(7.5,-0.3)(0,0.4){10}{\line(0,1){0.1}}	
		\put(8,0.3){\textcolor{red}{\vector(0,-1){1.2}}}
		\put(7.6,1.5){\line(1,0){0.8}} 
		\put(7.5,-1.5){\textcolor{red}{\scalebox{0.8}{$\boldsymbol{t = t_2^*}$}}}
		\put(7.8,1.75){\scalebox{0.8}{$u_{t_{2}^*}^*$}}
		\multiput(8.5,-0.3)(0,0.4){10}{\line(0,1){0.1}}		
		
		\multiput(9,0)(0.5,0){4}{\line(1,0){0.2}}
		\put(8.8,1.8){\oval(3.5,0.8)[l]}
		\multiput(8.9,2.2)(0.1,0){4}{\line(1,0){0.05}}
		\multiput(8.9,1.4)(0.1,0){4}{\line(1,0){0.05}}
		\put(6.3,1.6){$\mu_{3}^*$}
		\put(6.85,1.65){\line(1,0){0.2}}
		\put(6.85,1.80){\line(1,0){0.2}}
		
		\put(11,0){\line(1,0){2.5}}
		\multiput(11.5,-0.3)(0,0.4){10}{\line(0,1){0.1}}
		\put(12,0.3){\textcolor{red}{\vector(0,-1){1.2}}}
		\put(11.6,0.8){\line(1,0){0.8}}
		\put(11.5,-1.5){\textcolor{red}{\scalebox{0.8}{$\boldsymbol{t = t_{K^*}^*}$}}}
		\put(11.7,1.1){\scalebox{0.8}{$u_{t_{K^*}^*}^*$}}
		\multiput(12.5,-0.3)(0,0.4){10}{\line(0,1){0.1}} 
		\put(13,-0.3){\line(0,1){0.6}}
		\put(12.6,0.8){\line(1,0){0.8}}
		\put(12.3,-0.8){\scalebox{0.6}{$t = t_{K^*}^* + 1$}}
		\put(12.6,1.1){\scalebox{0.7}{$u_{t_{K^*}^* + 1}^*$}}
		\multiput(13.5,-0.3)(0,0.4){10}{\line(0,1){0.1}} 
		
		\multiput(12.5,0)(0.5,0){4}{\line(1,0){0.2}}  
		\multiput(13.2,0.8)(0.5,0){4}{\line(1,0){0.2}} 
		\put(13.8,1.1){\oval(5,0.9)}
		\put(9.7,1){$\mu_{K^* + 1}^*$}
		\put(11.1,1){\line(1,0){0.2}}
		\put(11.1,1.15){\line(1,0){0.2}}
		
		\put(14.5,0){\line(1,0){1.5}}
		\multiput(14.75,-0.3)(0,0.4){10}{\line(0,1){0.1}} 
		\put(15.25,-0.3){\line(0,1){0.6}}
		\put(14.85,0.8){\line(1,0){0.8}}
		\put(14.95,-0.8){\scalebox{0.6}{$t = n$}}
		\put(14.35,-1.3){\scalebox{0.6}{$(t = t_{K^* + 1}^* - 1)$}}
		\put(15.1,1.0){\scalebox{0.8}{$u_n^*$}}
		\multiput(15.75,-0.3)(0,0.4){10}{\line(0,1){0.1}}
		
		\put(-4.2,3.7){\scalebox{0.8}{$Y_{1}$}}
		\put(-3.2,3.7){\scalebox{0.8}{$Y_{2}$}}
		\put(-2.2,3.7){\scalebox{0.8}{$Y_{3}$}}
		%\put(-1.2,3.7){\scalebox{0.8}{$Y_{4}$}}
		
		\put(1.5,3.7){\scalebox{0.8}{$Y_{t_1^* - 1}$}}
		\put(2.7,3.7){\scalebox{0.8}{$Y_{t_1^*}$}}
		\put(3.5,3.7){\scalebox{0.8}{$Y_{t_1^* + 1}$}}
		%\put(4.5,3.7){\scalebox{0.8}{$Y_{t_1^* + 2}$}}
		
		\put(6.5,3.7){\scalebox{0.8}{$Y_{t_2^* - 1}$}}
		\put(7.7,3.7){\scalebox{0.8}{$Y_{t_2^*}$}}
		%\put(8.5,3.7){\scalebox{0.8}{$Y_{t_2^* + 1}$}}
		
		\put(11.7,3.7){\scalebox{0.8}{$Y_{t_{K^*}^*}$}}
		\put(12.5,3.7){\scalebox{0.8}{$Y_{t_{K^*}^* + 1}$}}
		\put(15,3.7){\scalebox{0.8}{$Y_n$}}
		
		\put(9.8, -1.4){\textcolor{red}{$\boldsymbol{\dots}$}}
		
		\end{picture}}
	}
	\caption{\footnotesize An illustration of the model  defined by  \eqref{eq1}, and \eqref{eq1_alt} respectively. The change-point locations are denoted as $t_{1}^{*}, \dots, t_{K}^{*}$ and there are $K^* + 1$ model phases represented by the  parameters $\mu_1^*, \dots, \mu_{K^* + 1}^*$.} 
\end{figure}
\begin{rem} 
 	The model above can be also seen a sampling scheme within some fixed domain, for instance, interval $(0,1)$. In such case the change-point locations can be understood as some specific points $\tau_k^*$, for $k = 1, \dots, K^*$, such that $ {t_{k}^*}/{n} \to \tau_k^* \in (0,1)$ for $n \to \infty$ and any $k \in \{1, \dots, K^*\}$. The unknown model segments $\mu_{1}^*, \dots, \mu_{K^{*} + 1}^*$ are determined by a fixed sequence of the true change-point locations $0 < \tau_{1}^* < \dots < \tau_{K^*}^* < 1$ and $K^* \in \mathbb{N}$, which is also fixed.
\label{remark1}
\end{rem}
The formulation in \eqref{eq2} introduces a kind of sparsity principle in parameters $u_{t}^*$, for $t = 1, \dots, n$, as we assume that $u_{t}^* = u_{t - 1}^*$, for all $t = 2, \dots, n$, but only $K^*$ specific exceptions for $t \in \{t_{1}^{*}, \dots, t_{K^*}^*\}$.
In order to estimate the vector of unknown parameters $\boldsymbol{u}^{*} = (u_{1}^*, \dots, u_{n}^*)^\top \in \mathbb{R}^{n}$, and the locations where $u_{t}^* \neq u_{t - 1}^*$, we solve the  minimization problem
\begin{equation}
\label{eq3}
\widehat{\boldsymbol{u}} = \argmin_{(u_{1}, \dots, u_{n}) \in \R^n} \bigg( \sum^n_{i=1}\rho_\tau (Y_i-u_i) + n \lambda_n\sum^{n-1}_{i=1} |u_{i+1}- u_i | \bigg),
\end{equation}
with $\widehat{\boldsymbol{u}} = (\widehat{u}_{1}, \dots, \widehat{u}_{n})^\top$, and $\rho_{\tau}(v) = v(\tau - \e1(v < 0))$, for some $\tau \in (0,1)$, and any  $v \in \mathbb{R}$. The regularization parameter $\lambda_{n} > 0$ controls for the overall number of change-points in the final model: for $\lambda_{n} = 0$ the minimization in \eqref{eq3} results in  $\widehat{\boldsymbol{u}}$ where $\widehat{u}_{t} \neq \widehat{u}_{t - 1}$, for each $t = 2, \dots, n$, while for $\lambda_{n} \to \infty$ we have $\widehat{u}_{t} = \widehat{u}_{t - 1}$, for all $t = 2, \dots, n$, and thus, the final model corresponds to a standard quantile linear regression model for the given $\tau \in (0,1)$. 

Using a parameter substitution and some  algebra calculations (analogously to \cite{tibs2014}, where it was applied to the linear (and higher order) trend filtering) we can rewrite the model in (\ref{eq1}) in terms of an ordinary 
linear regression model as
\begin{equation}
\label{eq4}
\eY^n=\eeX_n \eb^n +\eeps^n,
\end{equation}
where $\eY^n \equiv (Y_1, \cdots, Y_n)^\top$,  $\eb^n \equiv (d_{t^*_0}, 0, \cdots , 0,d_{t^*_1}, 0, \cdots , 0,d_{t^*_{K^*}}, 0, \cdots , 0 )^\top$, and $\eeps^n \equiv (\varepsilon_1, \cdots , \varepsilon_n)^\top$ 
with  
$ d_{t^*_k}$ on the position $t^*_k$, for  $k = 0, \dots, K^*$, $d_{t^*_0}=\mu^*_1$, and $d_{t^*_k}=\mu^*_k-\mu^*_{k-1}$, for $k=2, \cdots , K^*+1$. The model matrix, of the type $n \times n$, takes the from
\[
\eeX_n \equiv \left[
\begin{array}{ccccccccc}
1& & 0&& 0&& \cdots && 0 \\
1 && 1& & 0& &  \cdots & & 0 \\
1 & &1 & &1&&  \cdots && 0 \\
\vdots && \vdots && \vdots && \cdots && 0 \\
1 & &1 && 1& & \cdots & &1 \\
\end{array}
\right].
\]

Let $\eX_i$ denotes the $i$-th row of $\eeX_n$ and let $\widehat{\eb^n} \equiv \big(\widehat \beta_1, \cdots , \widehat \beta_n  \big)^\top$ be the  solution of
the quantile LASSO minimization problem 
\begin{equation}
\label{eq5}
\widehat{\eb^n} = \argmin_{\eb \in \R^n} \bigg[ \sum^n_{i=1}\rho_\tau (Y_i-(\eeX_n \eb)_i) + n \lambda_n \sum_{i = 2}^n | \beta_{i} | \bigg],
\end{equation}
where $(\eeX_n \eb)_i=\eX_i \eb$. Let $\widehat{\cal A}_n  \equiv \big\{i \in \{ 2, \cdots , n \}; \;\; \widehat \beta_i  \neq 0 \big\} = \big\{\widehat t_1, \cdots , \widehat t_{|\widehat{\cal A}_n |} \big\}$ be the set of  estimated change-point locations  and the corresponding  estimates of $u_i$ are defined as
\begin{equation}
\label{eq7}
\widehat u_i \equiv \big( \eeX_n \widehat{\eb^n} \big)_i=\eX_i  \widehat{\eb^n}, \qquad \textrm{ for } i=1, \cdots , n.
\end{equation}

\begin{rem}
	For brevity, we use  the notation were we suppress the dependence of the estimates $\widehat{\eb^n}$, $\widehat{\cal A}_n$, and $\widehat u_i$ on the value of the regularization parameter  $\lambda_n > 0$.
\end{rem}
The minimization problem defined in \eqref{eq5} is convex and it can be effectively solved using 
some standard optimization toolboxes. However, the parameter estimates for the vector of parameters $\eb^n$ are not given explicitly and iterative algorithms need to be employed to obtain the final solution. In the next section we  consider the model defined in \eqref{eq1} and we  derive and prove some theoretical properties for the estimation procedure defined by the minimization problem in \eqref{eq5}.

\section{Theoretical Results}
\label{results}
Let us start with introducing some necessary notation which will be used throughout this paper. 
Let $I^*_{min} \equiv \min_{1 \leqslant k \leqslant K^*}(t^*_{k+1}-t^*_k)$ and $I^*_{max} \equiv \max_{1 \leqslant k \leqslant K^*}(t^*_{k+1}-t^*_k)$. Analogously, for the change-point magnitudes, we define 
$J^*_{min} \equiv  \min_{1 \leqslant  k \leqslant K^*} | \mu^*_{k+1}-\mu^*_k|,$ and $J^*_{max} \equiv  \max_{1 \leqslant  k \leqslant K^*} | \mu^*_{k+1}-\mu^*_k|$. Obviously, we have $\mu^*_k \ne \mu^*_{k+1}$, for any  $k=1, \cdots, K^*$. Moreover, $C$ is used to denote a universal positive constant which does not depend on the sample size and which may take different values in different formulas. Let the model in \eqref{eq1} hold.  Then, in order to prove the results in this section, the following assumptions need to be satisfied:
\noindent
\begin{itemize}
\item[\textbf{(A1)}] The true parameters $\mu^*_k \in \mathbb{R}$, for any $k=1, \cdots, K^*+1$ do  not depend on $n \in \mathbb{N}$.\\[-0.6cm]
\item[\textbf{(A2)}] Random error  terms $\{\varepsilon_i\}$ are i.i.d., with some absolutely continuous distribution function $F(x)$, such that $\PP[\varepsilon <0]=\tau$, for the given quantile level $\tau \in (0,1)$, with the corresponding density function $f(x) > 0$,  for all $x \in \R$, which is continuously differentiable, such that $|f'(x)| < \infty$;  \\[-0.6cm]
\item[\textbf{(A3)}] Let $I^*_{min} \geq n \delta_n$, for some decreasing sequence $\{\delta_n\}$, such that $\delta_n \rightarrow 0$,   for $n \rightarrow \infty$; \\[-0.6cm]
\item[\textbf{(A4)}] Let, in addition, the following holds: $ \lambda_n / \delta_{n} \rightarrow 0$,  for $n \rightarrow \infty$;\\[-0.6cm]
\item[\textbf{(A5)}] We assume, that the number of change-points $K^* \in \mathbb{N}$ is fixed  and does not depend on the sample size $n \in \mathbb{N}$; \\[-0.6cm]
\item[\textbf{(A6)}] Let $\lambda_{n} = C(n^{-1} \log n)^{1/2}$, for some $C > 0$.

\end{itemize}

The assumption in (A1) specifies the model defined in \eqref{eq1} while Assumption (A2) is standard for the high-dimensional quantile regression models (see \cite{Koenker.05}).  
Assumption (A3) is considered, for instance, by \cite{Harchaoui.Levy.10} and \cite{Qian.Su.16} to ensure a proper change-point detection by the classical LASSO estimation approach: the authors in both these papers assume, among other assumptions, that  $(n \delta_n J^*_{min})^{-1} n \lambda_n \rightarrow 0$, for $n \rightarrow \infty$.  Thus, for $0 < J^*_{min} < \infty$ fixed, Assumption (A4)  in our paper corresponds to Assumption (A4) of \cite{Harchaoui.Levy.10} and also Assumption (A3)(iii) of \cite{Qian.Su.16}. 
Assumption (A5) on the true number of jumps $K^* \in \mathbb{N}$ is, for instance, considered in \cite{Harchaoui.Levy.10} for a least squares  model with $L_1$-penalty it is also quite reasonable in all practical applications. Assumption (A6) is needed in order to apply the results of \cite{Fan.Fan.Barut.14} on the convergence rate of the quantile LASSO estimator. Assumptions (A4) and (A6) imply that for the sequence $(\delta_n)$ from  Assumption (A3) that
\begin{equation*}
\label{A3bis}
\sqrt{\frac{n}{\log n}} \delta_n {\underset{n \rightarrow \infty}{\longrightarrow} }\infty.
\end{equation*}
This last relation implies that $(n \delta_n) \rightarrow \infty$ as $n \rightarrow \infty$.

\begin{rem}
Concerning the jump magnitudes, the assumptions imposed on $\{\delta_{n}\}$, $\{\lambda_{n}\}$, and  $J_{min}^{*}$ in \cite{Harchaoui.Levy.10} are the following: $n \delta_{n} (J_{min}^{*})^2 / \log n \to \infty$, and $(n \delta_{n} J_{min}^{*})^{-1} n \lambda_{n} \to 0$, for $n \to \infty$. Then, it is easy to see that for $n \delta_{n} = n \lambda_{n} = \log n$, it is necessary that $J_{min}^{*} \to \infty$. 
Thus, the smallest jump magnitude can not be bounded from above which obviously facilitates the detection of changes.  Therefore, the method presented in \cite{Harchaoui.Levy.10} requires the jump sizes to converge to infinity when the sample size increases.  In our present paper the jump magnitudes are all fixed.
\end{rem}

The main results of this paper are presented in the next three theorems. Theorem  \ref{Proposition 5} gives the convergence rate of the change-point location estimates if the number of the estimated change-points coincides with the true number of change-points. Theorem  \ref{Proposition 6}  covers the situation when the estimated number of change-points is greater than $K^*$, and finally, Theorem \ref{theorem 3.3} deals with a scenario where $\widehat{K}$ is smaller than $K^*$. All proofs are postponed to the appendix part in Section \ref{Append}. Let us firstly consider a situation when the estimated number of  change-points coincides with the reality---the true number of change-points $K^{*}$. In this case, with a probability converging to 1 as $n \rightarrow \infty$, the distance between the true location $t^*_k$ and the estimated location $\widehat t_k$ is smaller than $I^*_{min}$, which is the smallest distance between two consecutive true change-points.

\begin{theo}	
\label{Proposition 5}
	Let  $|\widehat{\cal A}_n|=K^*$. Then, under Assumptions (A1) -- (A6), it holds that 
	\[
	\PP \bigg[\max_{1 \leqslant k \leqslant K^*} |\widehat t_k -t^*_k | \geq n \delta_n  \bigg]\rightarrow 0, \qquad \textrm{ for } n \rightarrow\infty.
	\]
\end{theo}
 For the purpose of the second theorem, let us  introduce (similarly as in \cite{Harchaoui.Levy.10}), a distance between two sets, $A$ and $B$,  defined as
\[
{\cal E} \big(A || B \big) \equiv \sup_{b \in B} \inf_{a \in A} |a - b|.
\]
Let us also define two sets ${\cal T}^* \equiv \{t^*_1, \cdots , t^*_{K^*}   \}$ and $ \widehat{{\cal T}}_{|\widehat{\cal A}_n|} \equiv \{\widehat t_1, \cdots , \widehat t_{|\widehat{\cal A}_n|}  \}$.  In fact, the set  $ \widehat{{\cal T}}_{|\widehat{\cal A}_n|}$ is identical with  $\widehat{\cal A}_n$.  Thus, in the following theorem we show that  if the estimated number of change-points is greater than $K^*$ then the distance between $\widehat{{\cal T}}_{|\widehat{\cal A}_n|}$ and ${\cal T}^*$ is, with probability converging to one,  less than $n \delta_n$. Then, we can say that $ \widehat{{\cal T}}_{|\widehat{\cal A}_n|}$ is a weakly consistent estimator of ${\cal T}^*$. Le us start by studying the cardinality of the set $\widehat{\cal A}_n$.\\
We suppose that $\mu^*_1 \neq 0$. Otherwise, the reasoning is the same.
Let ${\cal A} \equiv \{1, t^*_1, \cdots , t^*_{K^*} \}$ which contains the indices where the vector $\eb^n$ has non-zero components, the elements $t^*_1, \cdots , t^*_{K^*}$ being also the observations where the model (\ref{eq1}) changes. 

Consider the following $n \times (K^*+1)$-matrix: $\textbf{S} \equiv \eeX_{n, {\cal A} }$, where $\eeX_{n, {\cal A} }$ is the submatrix formed by columns of $\eeX_{n }$ with indices in ${\cal A} $. Then, the $(K^*+1)$-square matrix
\[
\frac{1}{n} \textbf{S}^t \textbf{S}= 
\left[
\begin{array}{ccccccccc}
1& & 1- \frac{t^*_1}{n}&&1- \frac{t^*_2}{n}&& \cdots & &1- \frac{t^*_{K^*}}{n} \\
1- \frac{t^*_1}{n} && 1- \frac{t^*_1}{n} && 1- \frac{t^*_2}{n}& & \cdots && 1- \frac{t^*_{K^*}}{n} \\
\vdots & &\vdots & &\vdots && \cdots & &\vdots \\
1- \frac{t^*_{K^*}}{n} & &1- \frac{t^*_{K^*}}{n} && 1- \frac{t^*_{K^*}}{n}& & \cdots && 1- \frac{t^*_{K^*}}{n} \\
\end{array}
\right]
\] 
has all the leading principal minors equal to $1$, $\left(1- \frac{t^*_1}{n} \right)\frac{t^*_1}{n}$, $\left(1- \frac{t^*_2}{n} \right)\left(\frac{t^*_2}{n}- \frac{t^*_1}{n} \right)\frac{t^*_1}{n}$, $ \cdots$, 
$\left(1- \frac{t^*_{K^*}}{n} \right) \left(\frac{t^*_{K^*}}{n}- \frac{t^*_{K^*-1}}{n} \right)\left(\frac{t^*_{K^*-1}}{n}- \frac{t^*_{K^*-2}}{n} \right) \cdots \frac{t^*_1}{n}$.

By the characterization of Sylvester symmetric matrices, we have that $n^{-1} \textbf{S}^\top \textbf{S} $ is positive-definite. Moreover, since $ {t^*_k}/{n} {\underset{n \rightarrow \infty}{\longrightarrow} } \tau^*_k$, for all $k=1, \cdots , K^*$ (see Remark \ref{remark1}), we have that there exists a constant $C >0$ such that
\[
0 < C \leq \lambda_{min}\Big(\frac{1}{n} \textbf{S}^\top\textbf{S}\Big) \leq \lambda_{max}\Big(\frac{1}{n} \textbf{S}^\top \textbf{S}\Big) \leq \frac{1}{C},
\]
where $ \lambda_{min}(n^{-1} \textbf{S}^\top \textbf{S})$ and  $\lambda_{max}(n^{-1} \textbf{S}^\top \textbf{S}) $ are the smallest and the largest eigenvalues of the matrix $n^{-1} \textbf{S}^\top \textbf{S}$.
Let us consider also the $n \times (n-K^*)$-matrix $\textbf{Q} \equiv \eeX_{n, \overline{ \cal A} }$, where $ \overline{ \cal A}$ is the complementary of ${ \cal A}$. Then, there again exists a positive constant   $C$ such that: 
\[
\Big\| \frac{1}{n} \textbf{Q}^t \textbf{S} \Big\|_{2,\infty} < C,
\]
with $\|\textbf{A} \|_{2, \infty}=\sup_{\textbf{x} \neq \textbf{0}}\|\textbf{A}\textbf{x}\|_{\infty}/\|\textbf{x} \|_2$. 
Taking also into account Assumptions (A2), (A5), and (A6), and  applying Theorem 2 of \cite{Fan.Fan.Barut.14}, we have
\begin{equation}
 \label{Bel3}
 \widehat{\mu}_k -\mu^*_k=O_{\PP}\pth{\sqrt{\frac{\log n}{n}}},  \qquad \textrm{for any } k=1, \cdots , K^*+1, 
 \end{equation}
 and
 \begin{equation}
 \label{Bel2bis}
 \PP \big[ |\widehat{\cal A}_n| \leq C K^*  \big] {\underset{n \rightarrow \infty}{\longrightarrow} }1,
  \end{equation}
  for some constant $0 < C < \infty$. Therefore,  using Assumption (A5), we can conclude that the number of estimated change-points, $|\widehat{\cal A}_n|$, is bounded with probability converging to one.
  
\begin{theo} \label{Proposition 6}
	Let $|\widehat{\cal A}_n|\geq K^*$. Then, under Assumptions (A1) -- (A6), it holds that
	\begin{equation}
	\label{eq15}
	\PP \bigg[ {\cal E} \big( \widehat{{\cal T}}_{|\widehat{\cal A}_n|} || {\cal T}^* \big) \leq n \delta_n  \bigg]\rightarrow 1, \qquad \textrm{ as } n \rightarrow\infty.
	\end{equation}
\end{theo}

Taking into account Assumptions (A3), (A4), and (A6), we obtain that the minimum distance between two consecutive change-points, $I_{min}^{*}$, has to satisfy $(n \log n)^{-1/2}I_{min}^{*} \rightarrow \infty$. Thus, since $(n \log n)^{-1/2}I_{min}^{*} =\frac{I_{min}^{*} }{\log n} \sqrt{\frac{\log n}{n}}$, we have
\begin{equation}
\label{A7}
(\log n)^{-1} I_{min}^{*} % \geq n \delta_n (\log n)^{-1}
 \to \infty, \qquad \textrm{for $n \to \infty$.}
\end{equation}

Note, that (\ref{A7}) indicates that in order to avoid underestimation of the number of change-points the minimum distance between two consecutive change-points must be larger than $\log n$.

Finally, the last theorem proves that the estimated number of change-points is not underestimated, but, it is rather overestimated with probability tending to one as $n$ increases. In such cases, however, for each true change-point location $t_{k}^* \in  \cal T^*$ there is at least one change-point location estimate in $\widehat{{\cal T}}_{|\widehat{\cal A}_n|}$, such that the distance between the true location and the corresponding estimate is less than $n \delta_n$, again with probability tending to 1, for $n \to \infty$ (assertion of Theorem \ref{Proposition 6}). 

\begin{theo}\label{theorem 3.3}
		Under Assumptions (A1) -- (A6), we have that 
	$$\PP [ |\widehat{\cal A}_n | <K^*] \rightarrow 0, \qquad \textrm{as $n \rightarrow \infty$.}$$
	\end{theo}
	
As already mentioned before, the theorem above suggests that the number of estimated change-points is more likely to be overestimated. This is, indeed, a common nature of standard LASSO regularization approaches. On the other hand, the overestimation can be amended, for instance, by adopting the adaptive LASSO approach which is well-known for being able to solve the overestimation problem and, moreover, can achieve the oracle properties.

In the next section we compare the proposed quantile LASSO method with some other common estimation techniques and the presented theoretical results will be illustrated in terms of the empirical performance.

\section{Simulation Study}
\label{simulations}
In this section we investigate the finite sample properties of the proposed quantile LASSO estimator.  For  the simulation purposes we consider a location model defined as 
\begin{equation}
\label{sim1}
Y_{t} = \left\{\begin{array}{cc}
\mu_{1}^* + \varepsilon_{t} & \textrm{for $0 \leq t \leq t_{1}^*$}\\
\mu_{2}^* + \varepsilon_{t} & \textrm{for $t_{1}^* \leq t \leq t_{2}^*$}\\
\mu_{3}^* + \varepsilon_{t} & \textrm{for $t_{2}^* \leq t \leq 1$}
\end{array}\right., \quad \textrm{for $t = 1, \dots, n$, and $n \in \mathbb{N}$,}
\end{equation}
with two distant change-points $t_{1}^{*}, t_{2}^{*} \in \{1, \dots, n\}$, with the corresponding parameters $\mu_1^* = 0$, $\mu_{2}^* = 2$, and $\mu_{3}^* = 1$. The random error terms 
$\{\varepsilon_t\}_{t = 1}^n$ are considered to be independent but, in order to investigate a signal-to-noise performance and the robust favor of the proposed quantile LASSO method, we consider various error distributions (the standard normal distribution, Student's $t$ distribution with three degrees of freedom, and finally, the Cauchy distribution). 

Three different sample sizes  $n \in \{20, 100, 500\}$ are used but in order to be able to easily compare different models with various sample sizes the model is rescaled in terms of Remark \ref{remark1}, such that $Y_{\tilde{t}} \equiv Y_{\frac{t}{n}}$, for $\tilde{t} \in (0,1)$, with the corresponding change-points being located at $t_1^* / n \to \tau_{1}^* = 0.2$ and $t_{2}^*/n \to \tau_2^* = 0.7$. 

A set of quantile levels for $\tau \in \{0.05, 0.10, 0.25, 0.50, 0.75, 0.90, 0.95\}$ is considered and the final model is  estimated using Equation \eqref{eq5}, while three different approaches are  applied to determine the value of the regularization parameter $\lambda_{n} > 0$: firstly, we considered the asymptotically appropriate value fulfilling Assumption (A6), denoted as $\lambda_{AS}$, where $\lambda_{AS} = C \cdot (n^{-1} \log n)^{1/2}$. For the second model, we use the prior knowledge that there are two true change-points in the model (thus, the final model always contains two change-points and three segments and the corresponding regularization parameter is denoted as $\lambda_{(2)}$). Finally, for the third model, we consider the parameter denoted as $\lambda_{MS}$ which is determined by the minimum Mean Squared Error (MSE) quantity $n^{-1}\sum_{i = 1}^{n} (\widehat{u}_{i}^{*} - u_{i}^{*})^2$.
In addition, we compare the quantile LASSO method with the standard LASSO approach and the SMUCE estimator proposed by \cite{Frick.14}. The final models are compared with respect to the averaged estimation bias given by $n^{-1}\sum_{i = 1}^{n} ( \widehat{u}_{i}^{*} - u_{i}^{*} )$, the  MSE quantity, and the change-point detection error expressed as $\frac{1}{2} \sum_{k = 1}^2 |\widehat{t}_{k} - t_{k}^*|$. The change-point detection error is, however,  only obtained for models where at least two change-points are detected (otherwise, the quantity is not reported). The simulations are based on 1000 Monte Carlo repetitions for every possible model scenario and the results are reported in Tables \ref{tab1}, \ref{tab2}, and \ref{tab3}. 
 
\begin{figure}
	\centering
	\subfigure[$n = 20$ and $\lambda_{n} = \lambda_{AS}$]{\label{fig1:b}\includegraphics[width=0.48\textwidth, height = 0.25\textwidth]{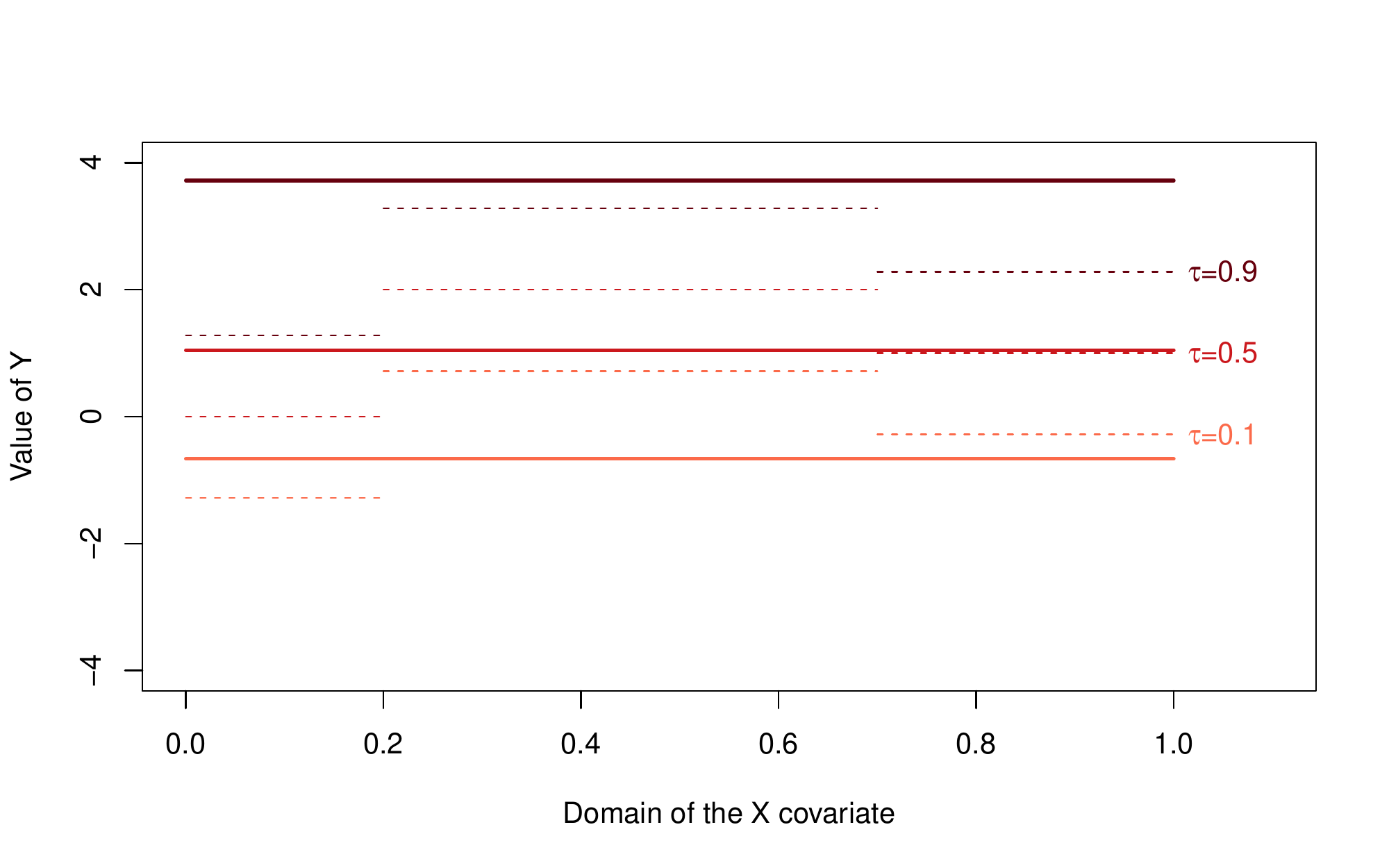}}
	\subfigure[$n = 20$ and $\lambda_{n} = \lambda_{(2)}$]{\label{fig1:a}\includegraphics[width=0.48\textwidth, height = 0.25\textwidth]{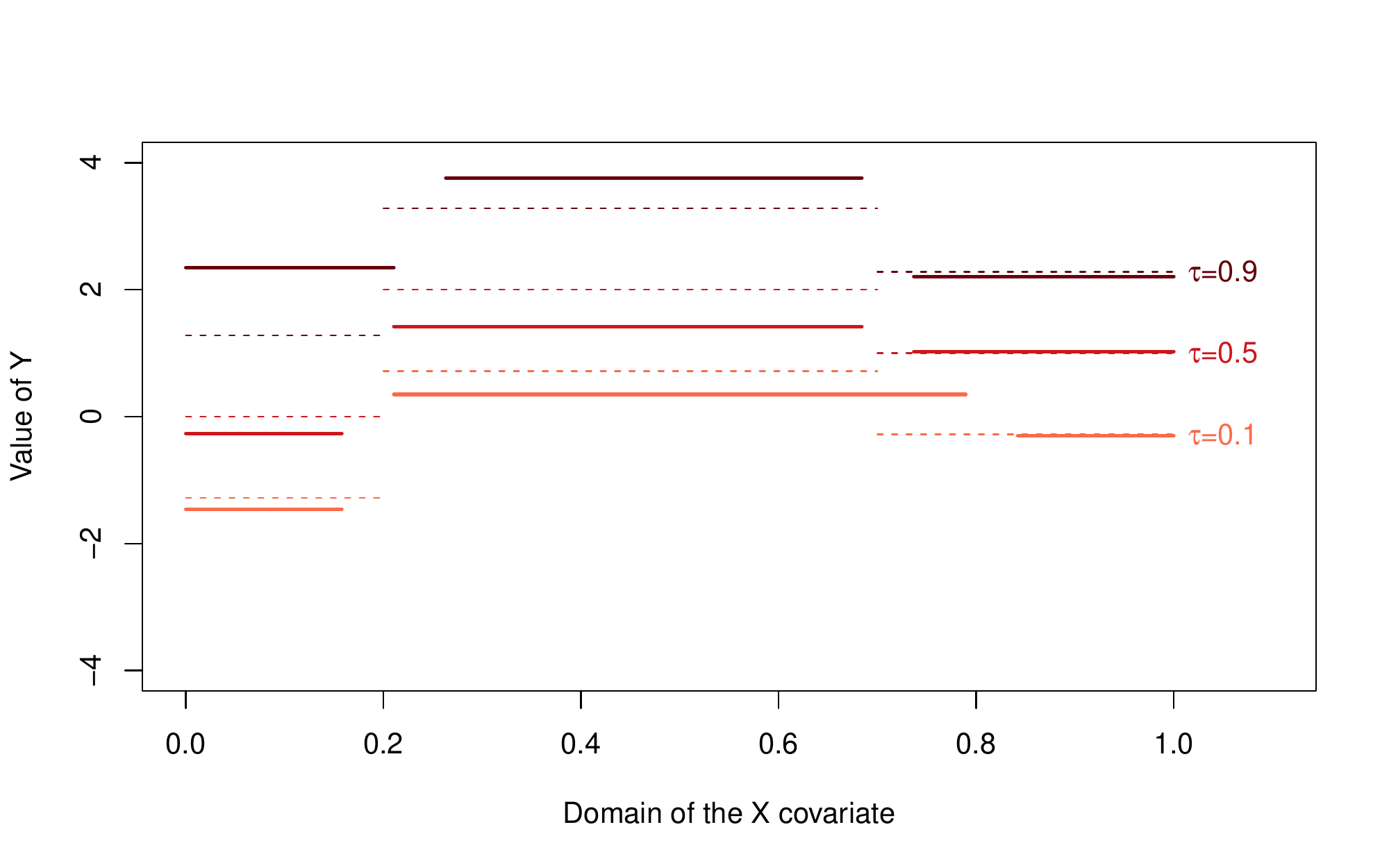}}
	
	\subfigure[$n = 100$ and $\lambda_{n} = \lambda_{AS}$]{\label{fig1:d}\includegraphics[width=0.48\textwidth, height = 0.25\textwidth]{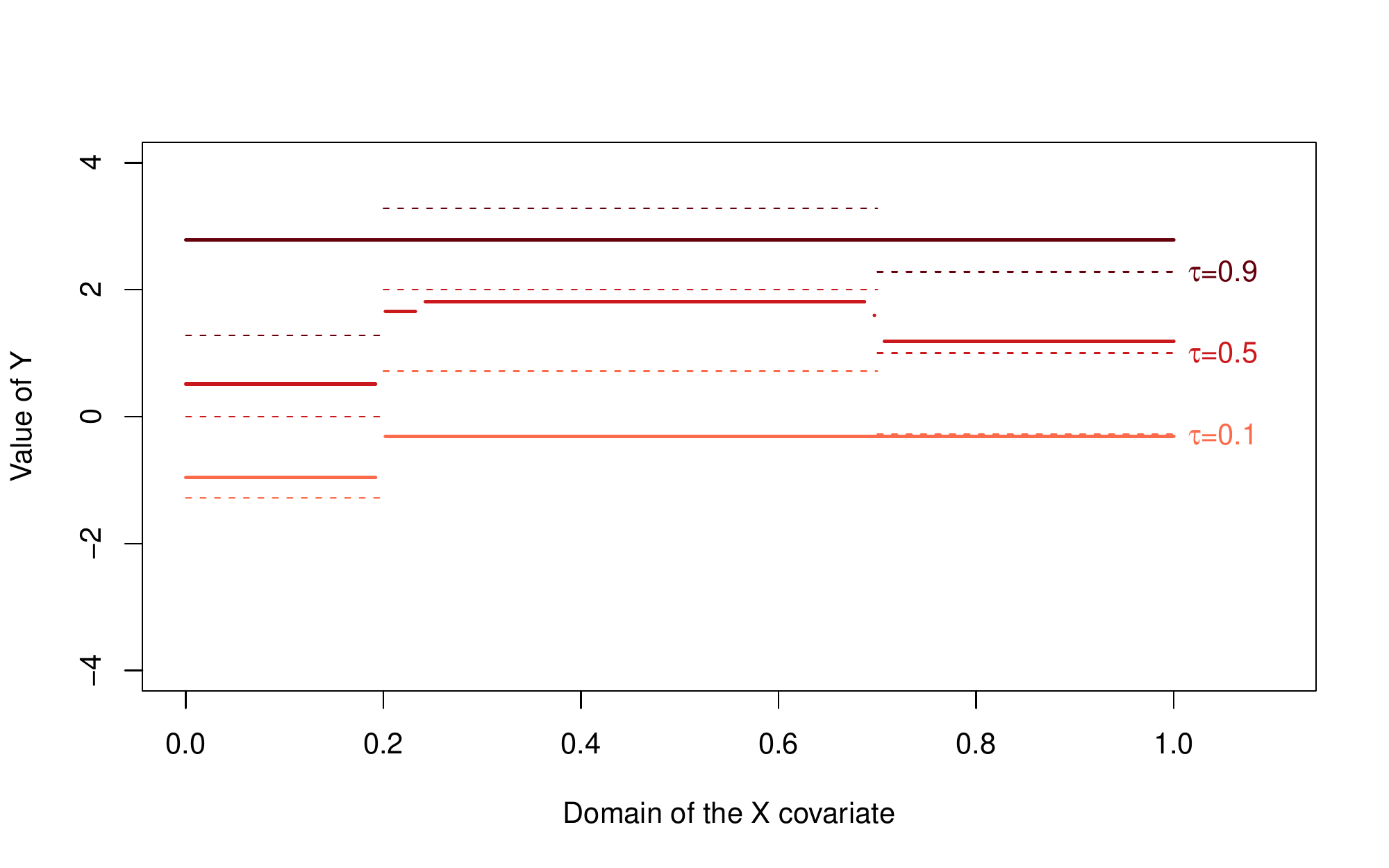}}
	\subfigure[$n = 100$ and $\lambda_{n} = \lambda_{(2)}$]{\label{fig1:c}\includegraphics[width=0.48\textwidth, height = 0.25\textwidth]{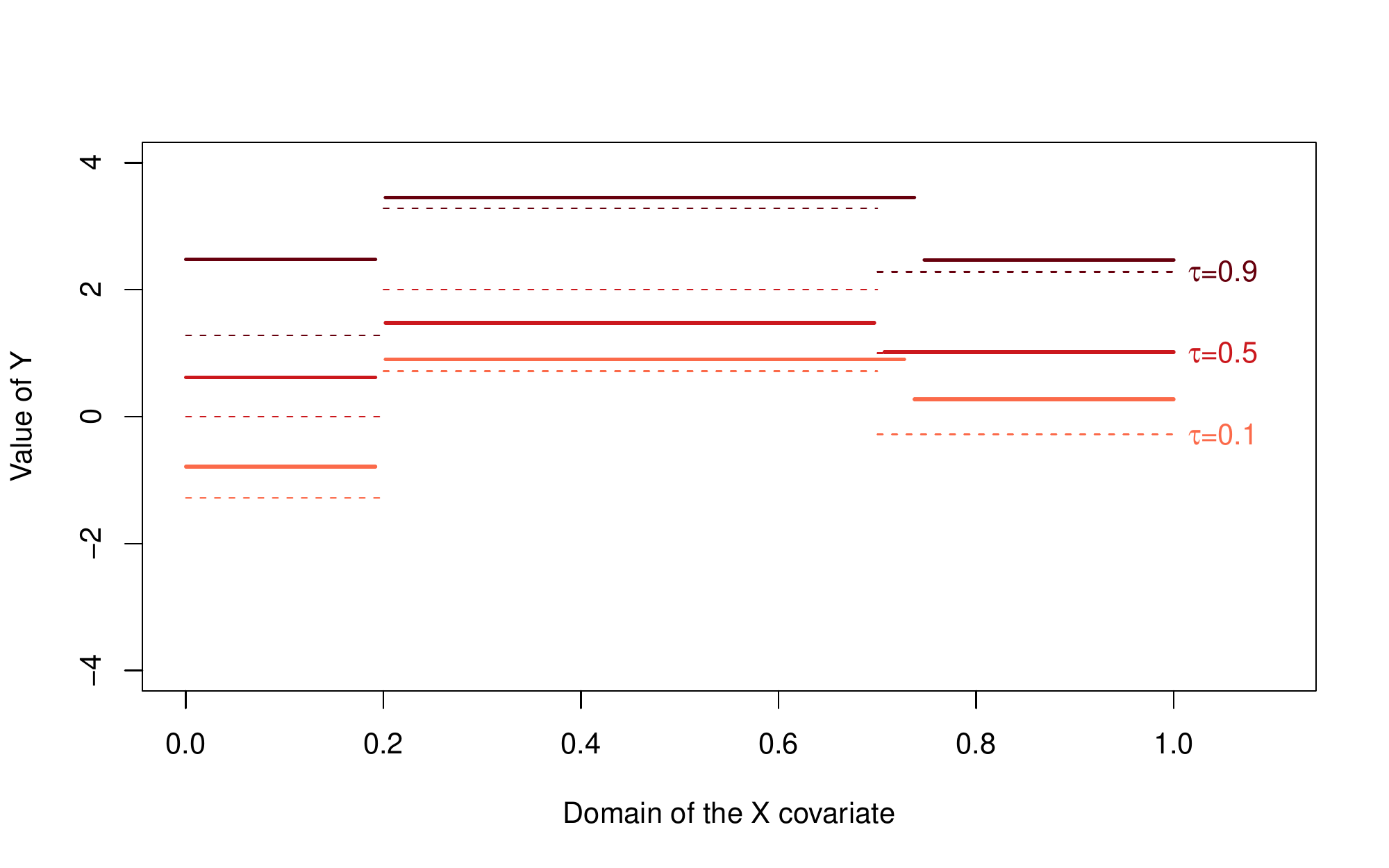}}
	
	\subfigure[$n = 500$ and $\lambda_{n} = \lambda_{AS}$]{\label{fig1:f}\includegraphics[width=0.48\textwidth, height = 0.25\textwidth]{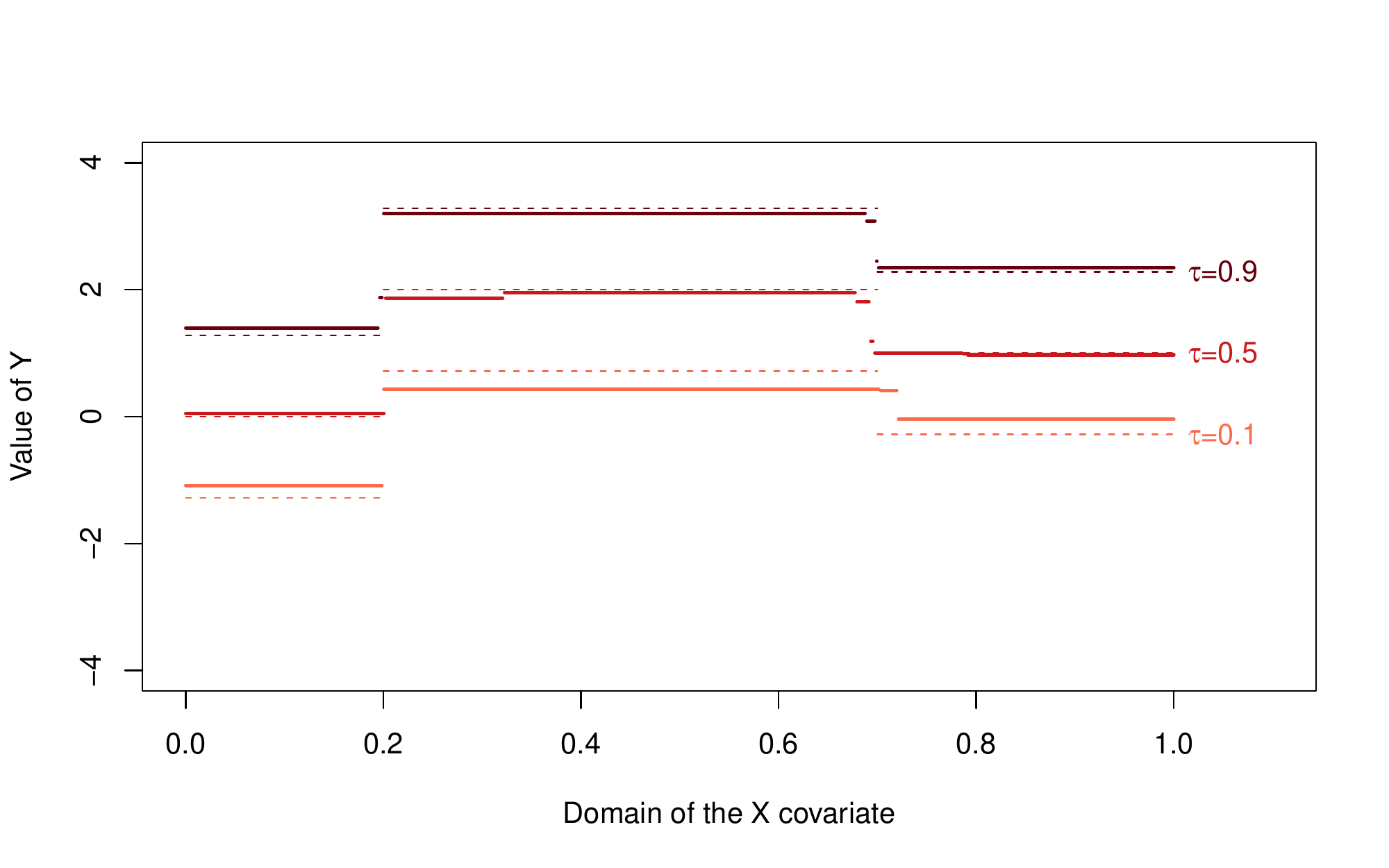}}
	\subfigure[$n = 500$ and $\lambda_{n} = \lambda_{(2)}$]{\label{fig1:e}\includegraphics[width=0.48\textwidth, height = 0.25\textwidth]{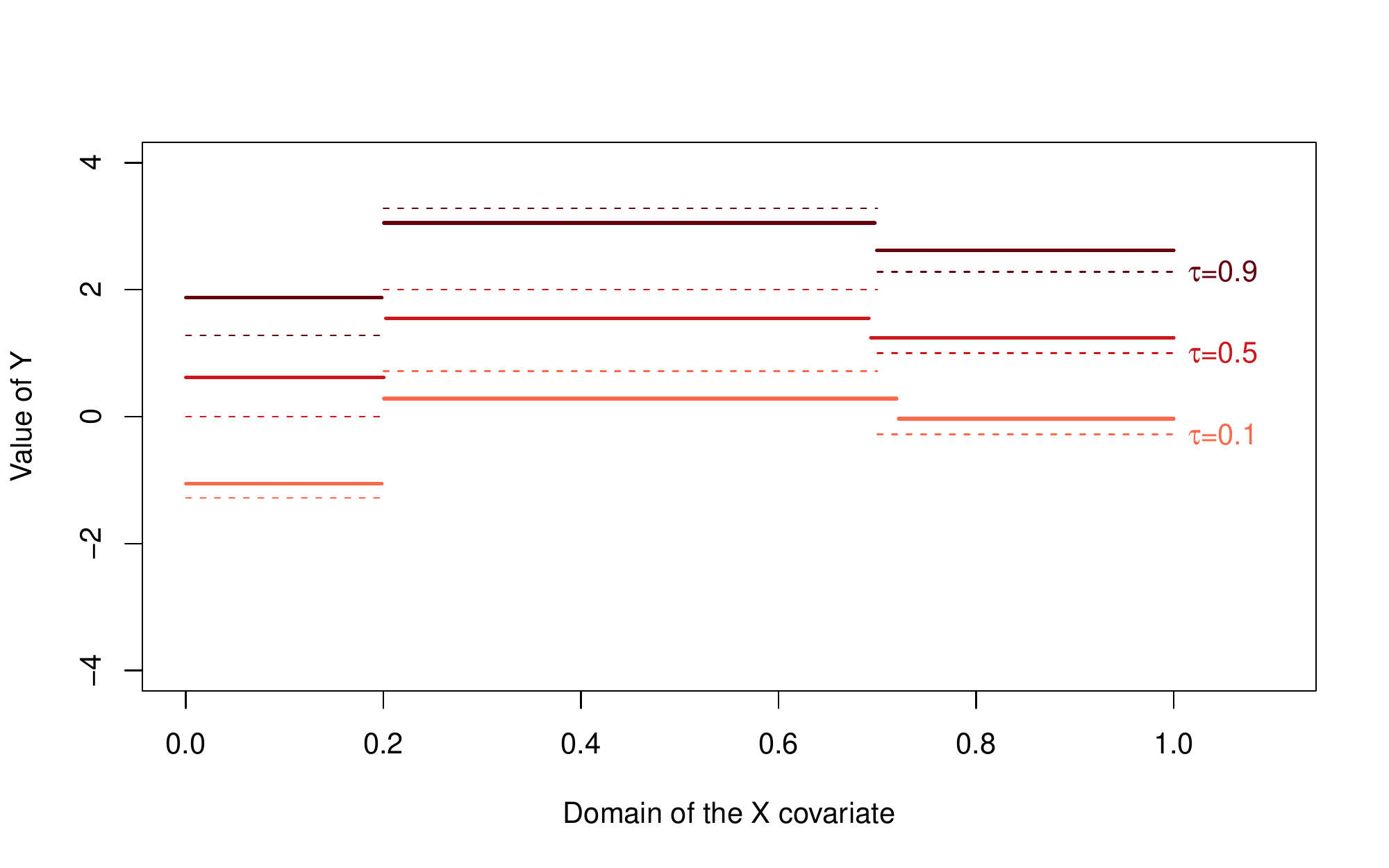}}

	\caption{\footnotesize An illustration of the performance of the quantile LASSO estimator for the asymptotically appropriate value of the regularization parameter $\lambda_{AS} = C \cdot (n^{-1} \log n)^{1/2}$ on the left-hand side, and the prior knowledge that two change-points should be recovered in the final model used on the right-hand side (with the corresponding regularization parameter denoted as $\lambda_{(2)}$). The true regression quantiles for $\tau \in \{0.1, 0.5, 0.9\}$ are given by dashed lines and the estimated quantiles are plotted by solid lines. The model with $\lambda_{AS}$ (left-side figures) outperforms the model with $\lambda_{(2)}$ (right-side figures) with respect to a smaller bias but, on the other hand, it overfits the data with respect to the number of detected change-points.}
	\label{fig1}
\end{figure}

\begin{figure}
	\centering
	\subfigure[$\lambda_{AS} \quad | \quad \tau = 0.10$]{\label{fig2:a}\includegraphics[width=0.48\textwidth, height = 0.25\textwidth]{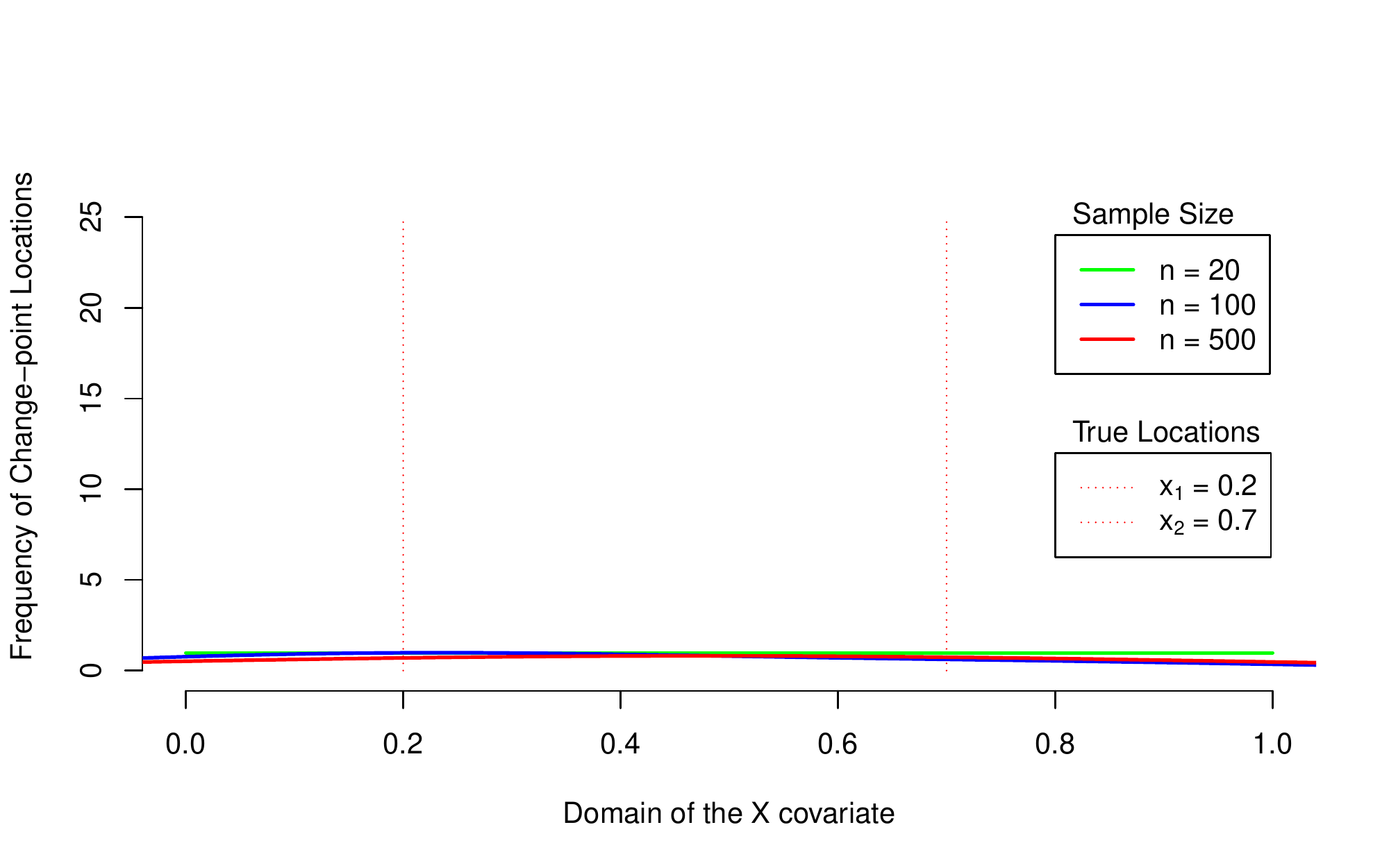}}
	\subfigure[$\lambda_{(2)} \quad | \quad \tau = 0.10$]{\label{fig2:b}\includegraphics[width=0.48\textwidth, height = 0.25\textwidth]{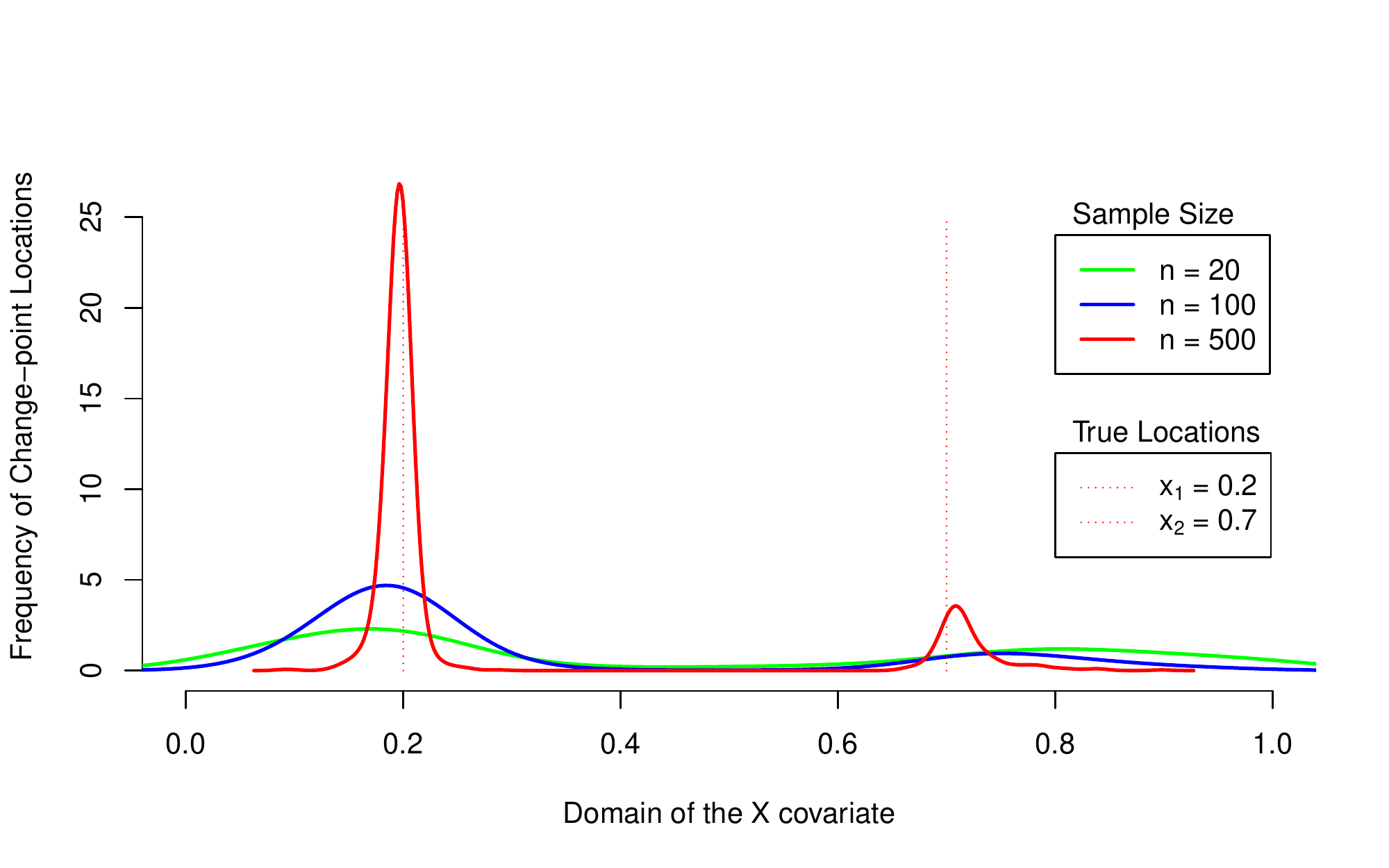}}
	
	\subfigure[$\lambda_{AS} \quad | \quad \tau = 0.25$]{\label{fig2:c}\includegraphics[width=0.48\textwidth, height = 0.25\textwidth]{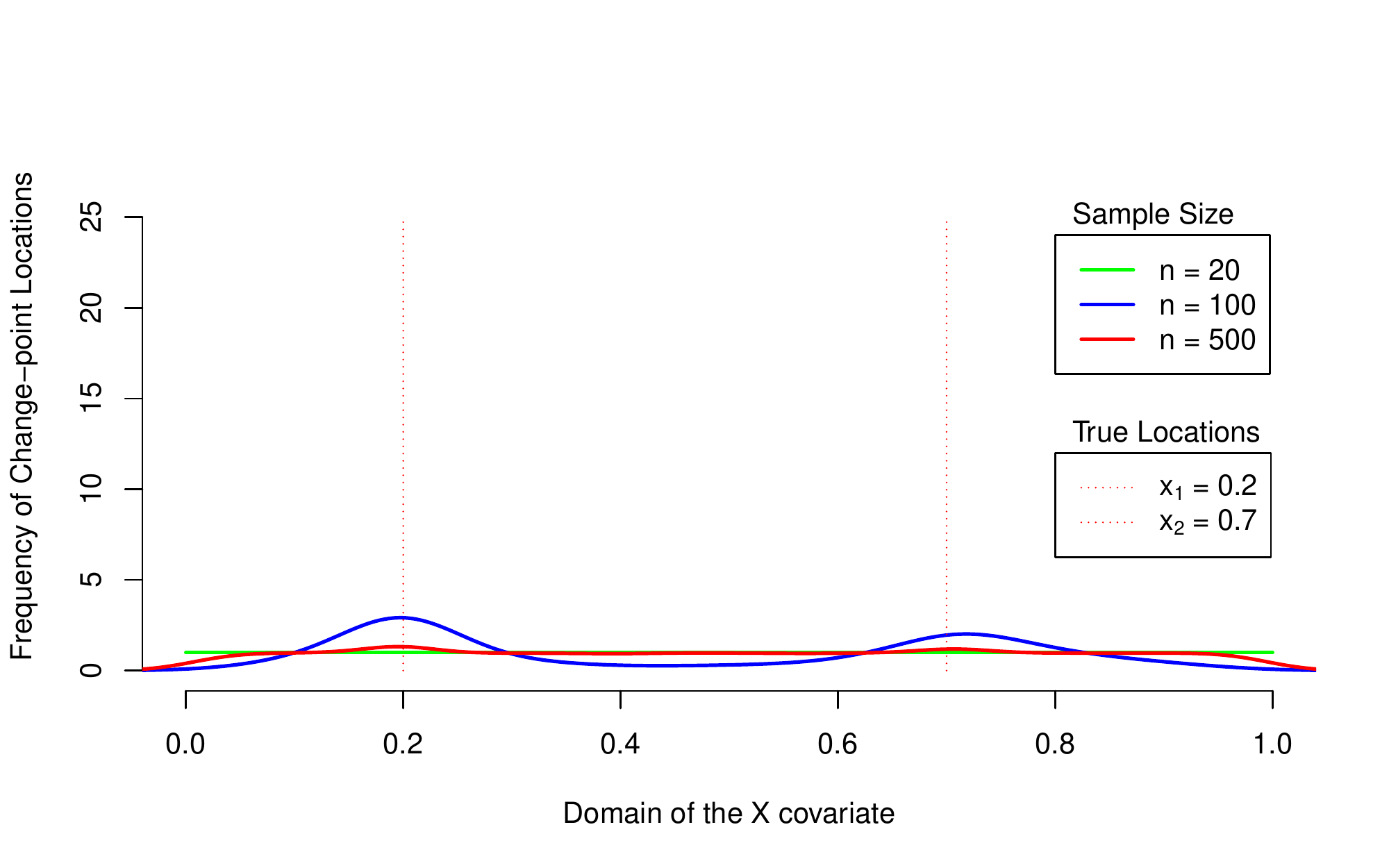}}
	\subfigure[$\lambda_{(2)} \quad | \quad \tau = 0.25$]{\label{fig2:d}\includegraphics[width=0.48\textwidth, height = 0.25\textwidth]{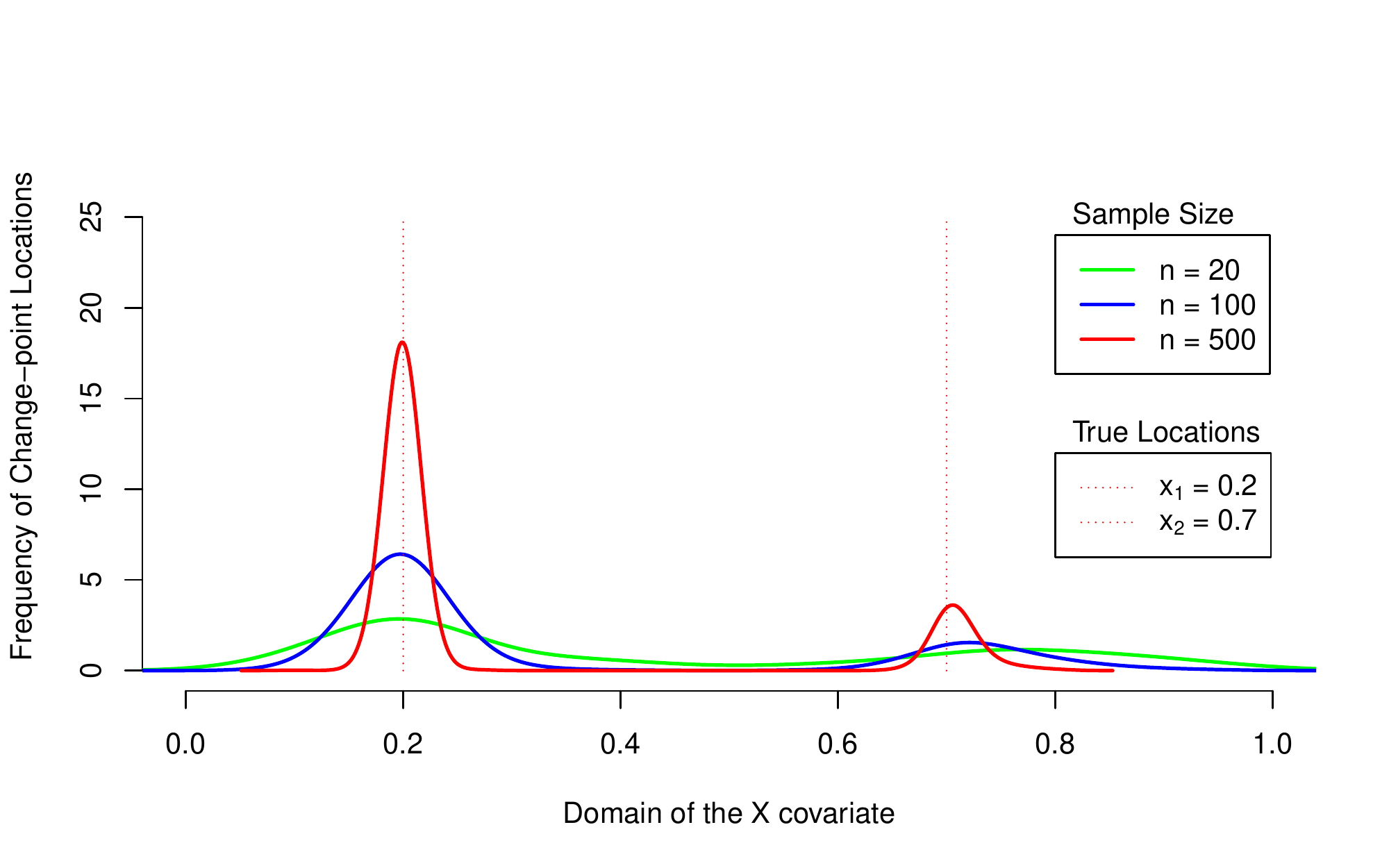}}
	
	\subfigure[$\lambda_{AS} \quad | \quad \tau = 0.50$]{\label{fig2:e}\includegraphics[width=0.48\textwidth, height = 0.25\textwidth]{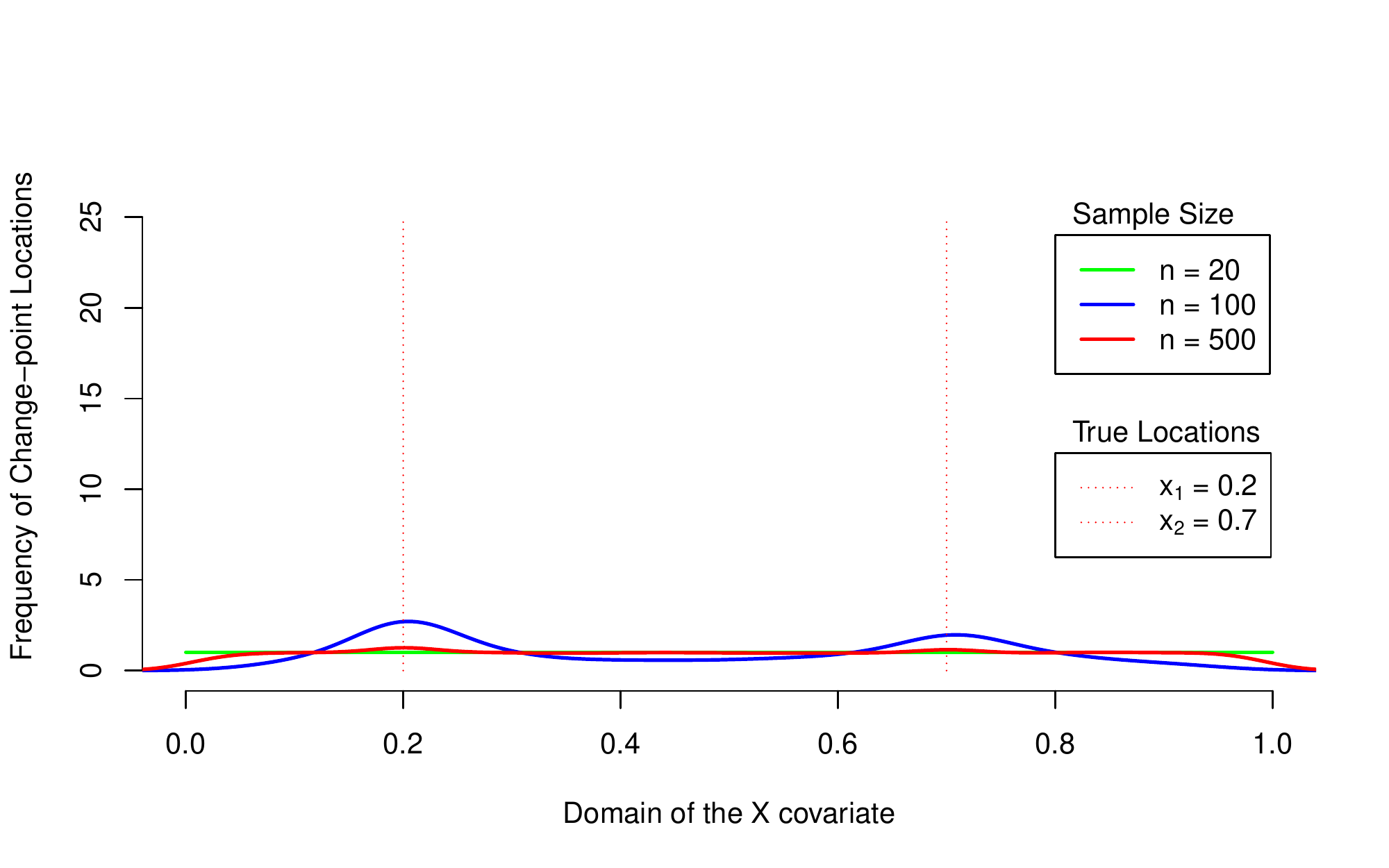}}
	\subfigure[$\lambda_{(2)} \quad | \quad \tau = 0.50$]{\label{fig2:f}\includegraphics[width=0.48\textwidth, height = 0.25\textwidth]{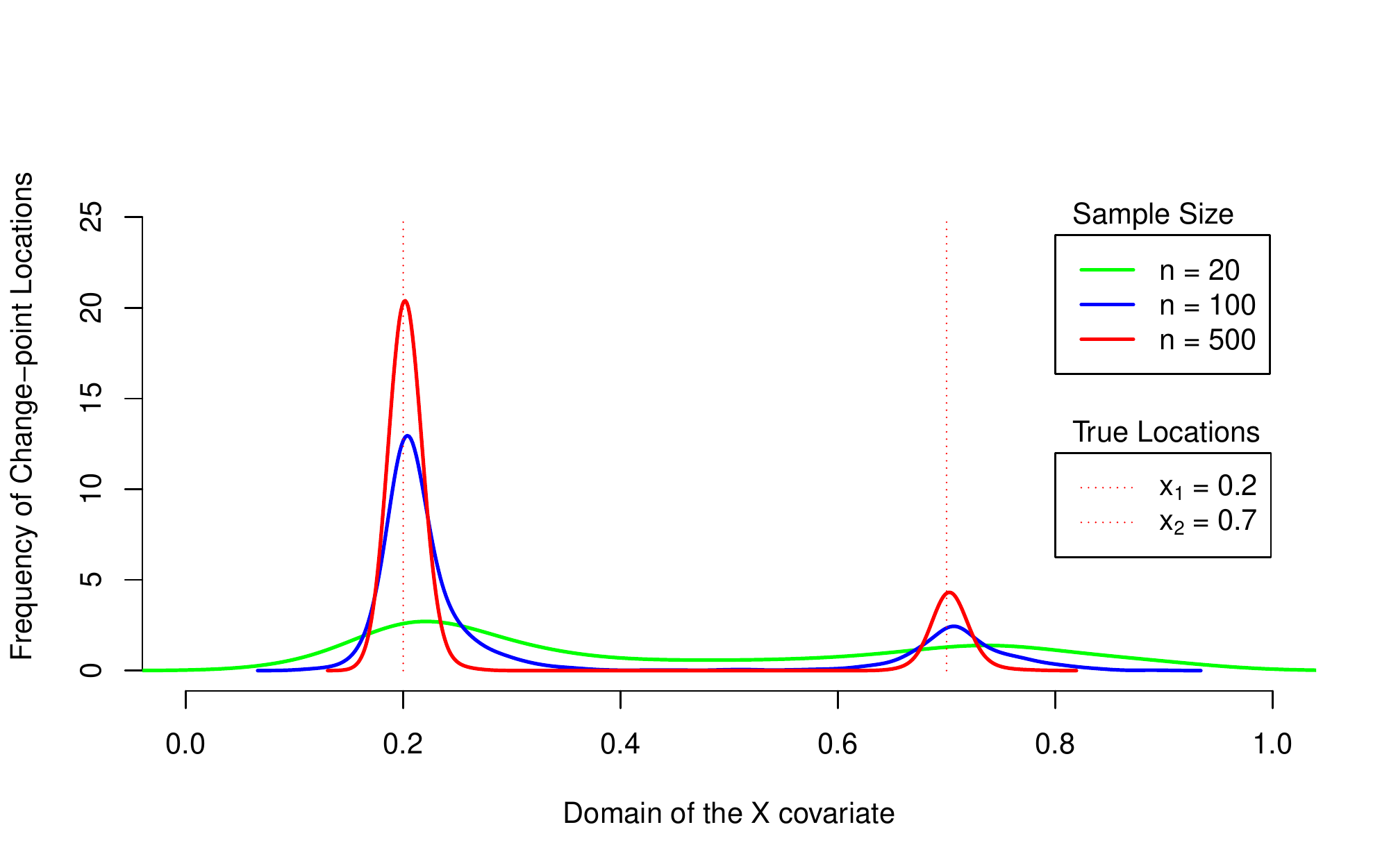}}
	\caption{\footnotesize The jump detection performance of the quantile LASSO for various sample sizes and quantile levels $\tau \in \{0.10, 0.25, 0.50\}$ (the results for $\tau = 0.75$ and $\tau = 0.90$ are analogous to the results for $\tau = 0.25$ and $\tau = 0.10$). The asymptotically appropriate value $\lambda_{AS} = C \dot (n^{-1} \log n)^{1/2}$  is considered on the left-hand side and the prior knowledge about two change-points in the model is used in the right-hand side figure (with the corresponding regularization parameter denoted as $\lambda_{(2)}$). The overestimation of the model with $\lambda_{AS}$ is obvious by the more spread densities around the true change-point locations (left-side panels). Also, the model with $\lambda_{AS}$ performs poorly for small sample sizes: no change-points are recovered for $n = 20$. In contrast to that, the models with $\lambda_{(2)}$ seem to be able to recover true jump locations consistently (right-side panels).}
	\label{fig2}
\end{figure}

First of all, we are primarily interested in the quantile LASSO performance when estimating different quantile levels  (the results are summarized in Table \ref{tab1} and illustrated in Figures \ref{fig1} and \ref{fig2}). From the asymptotical point of view, the model with $\lambda_{AS}$ outperforms the model with two change-points (the model with the regularization parameter $\lambda_{(2)}$): the estimation bias and the MSE quantity are both  much smaller  for larger sample sizes. The model with $\lambda_{AS}$ selects more than just two change-points and therefore, it allows for more augmentation of the sparse vector of parameters and thus, a smaller bias. On the other hand, the model with two change-points is more reliable  when detecting the true change-point locations: the model with $\lambda_{AS}$ tends to select more non-zero parameters---change-points---than actually needed. This is, however, a common property of the LASSO methods in general and it could be slightly reduced by adopting, for instance, an~adaptive LASSO approach.
It is also worth to mention, that the quantile LASSO performs much better when estimating quantile levels close to the median value rather than the levels on the tails (see Table \ref{tab1}). This is however, a common fact  and such behavior is quite expected.

\begin{table}[!ht]\footnotesize
	\begin{center}
		\scalebox{0.95}{
		\begin{tabular}{cc|cr@{ }lr@{ }lr@{ }lr@{ }lr@{ }lc}
			\multirow{2}{*}{$\boldsymbol{n}$} & \multirow{2}{*}{$\boldsymbol{\tau}$} & $\boldsymbol{\lambda_{AS}}$ & \multicolumn{2}{c}{$\boldsymbol{\lambda_{(2)}}$} & \multicolumn{4}{c}{\textbf{Model with $\boldsymbol{\lambda_{(2)}}$}} & \multicolumn{4}{c}{\textbf{Model with $\boldsymbol{\lambda_{AS}}$}} & \multicolumn{1}{c}{$\boldsymbol{|\widehat{\mathcal{A}_{n}}|}$} \\
			~ & ~ & \scalebox{0.8}{Value} & \scalebox{0.8}{Mean} & \scalebox{0.8}{Std.Err.} & \multicolumn{2}{c}{\scalebox{1}{Est. Bias}} & \multicolumn{2}{c}{MSE} & \multicolumn{2}{c}{\scalebox{1}{Est. Bias}} & \multicolumn{2}{c}{MSE} & \multicolumn{1}{c}{\scalebox{0.8}{[M $|$ M $|$ M]}}  \\\hline\hline
			\textbf{20} & 0.05 & 3.87 & 0.30 & \textit{(0.09)} & -0.32 & \textit{(0.41)} & 0.66 & \textit{(0.35)} & 0.49 & \textit{(0.46)} & 1.06 & \textit{(0.55)} & [0$|$0$|$0]\\
			~ & 0.10 & 3.87 & 0.56 & \textit{(0.18)} & -0.02 & \textit{(0.40)} & 0.55 & \textit{(0.28)} & 0.39 & \textit{(0.39)} & 0.92 & \textit{(0.38)} & [0$|$0$|$0]\\
			~ & 0.25 & 3.87 & 1.08 & \textit{(0.29)} & 0.03 & \textit{(0.33)} & 0.45 & \textit{(0.24)} & 0.20 & \textit{(0.33)} & 0.76 & \textit{(0.21)} & [0$|$0$|$0]\\
			~ & 0.50 & 3.87 & 1.47 & \textit{(0.46)} & -0.01 & \textit{(0.31)} & 0.52 & \textit{(0.22)} & -0.04 & \textit{(0.30)} & 0.70 & \textit{(0.13)} & [0$|$0$|$0]\\
			~ & 0.75 & 3.87 & 1.04 & \textit{(0.18)} & -0.03 & \textit{(0.32)} & 0.56 & \textit{(0.27)} & -0.21 & \textit{(0.32)} & 0.76 & \textit{(0.23)} & [0$|$0$|$0]\\
			~ & 0.90 & 3.87 & 0.52 & \textit{(0.17)} & 0.07 & \textit{(0.41)} & 0.80 & \textit{(0.35)} & -0.33 & \textit{(0.40)} & 0.88 & \textit{(0.40)} & [0$|$0$|$0]\\
			~ & 0.95 & 3.87 & 0.28 & \textit{(0.09)} & 0.34 & \textit{(0.44)} & 1.04 & \textit{(0.45)} & -0.40 & \textit{(0.49)} & 1.01 & \textit{(0.60)} & [0$|$0$|$0]\\
			\multicolumn{2}{c}{~} & \multicolumn{12}{c}{~}\\[-0.2cm]
			\textbf{100} & 0.05 & 2.15 & 1.67 & \textit{(0.57)} & 0.19 & \textit{(0.25)} & 0.45 & \textit{(0.22)} & 0.31 & \textit{(0.26)} & 0.61 & \textit{(0.27)} & [0$|$1$|$6]\\
			~ & 0.10 & 2.15 & 2.63 & \textit{(0.85)} & 0.17 & \textit{(0.20)} & 0.40 & \textit{(0.17)} & 0.09 & \textit{(0.19)} & 0.31 & \textit{(0.14)} & [1$|$3$|$10]\\
			~ & 0.25 & 2.15 & 4.73 & \textit{(1.31)} & 0.10 & \textit{(0.15)} & 0.35 & \textit{(0.13)} & 0.00 & \textit{(0.14)} & 0.14 & \textit{(0.07)} & [1$|$6$|$15]\\
			~ & 0.50 & 2.15 & 6.22 & \textit{(1.51)} & -0.01 & \textit{(0.14)} & 0.41 & \textit{(0.15)} & -0.01 & \textit{(0.13)} & 0.12 & \textit{(0.06)} & [2$|$9$|$21]\\
			~ & 0.75 & 2.15 & 4.40 & \textit{(0.91)} & -0.12 & \textit{(0.16)} & 0.44 & \textit{(0.18)} & 0.00 & \textit{(0.14)} & 0.16 & \textit{(0.09)} & [1$|$6$|$15]\\
			~ & 0.90 & 2.15 & 2.29 & \textit{(0.45)} & -0.19 & \textit{(0.19)} & 0.56 & \textit{(0.19)} & -0.17 & \textit{(0.19)} & 0.53 & \textit{(0.18)} & [0$|$3$|$9]\\
			~ & 0.95 & 2.15 & 1.36 & \textit{(0.27)} & -0.20 & \textit{(0.22)} & 0.64 & \textit{(0.20)} & -0.37 & \textit{(0.22)} & 0.80 & \textit{(0.18)} & [0$|$0$|$3]\\
			\multicolumn{2}{c}{~} & \multicolumn{12}{c}{~}\\[-0.2cm]
			\textbf{500} & 0.05 & 1.11 & 7.73 & \textit{(2.51)} & 0.24 & \textit{(0.13)} & 0.41 & \textit{(0.15)} & -0.06 & \textit{(0.09)} & 0.10 & \textit{(0.03)} & [7$|$15$|$25]\\
			~ & 0.10 & 1.11 & 12.82 & \textit{(4.44)} & 0.19 & \textit{(0.10)} & 0.36 & \textit{(0.13)} & -0.06 & \textit{(0.07)} & 0.09 & \textit{(0.03)} & [15$|$24$|$41]\\
			~ & 0.25 & 1.11 & 22.60 & \textit{(6.33)} & 0.10 & \textit{(0.07)} & 0.33 & \textit{(0.11)} & -0.04 & \textit{(0.06)} & 0.10 & \textit{(0.02)} & [37$|$53$|$69]\\
			~ & 0.50 & 1.11 & 30.23 & \textit{(6.61)} & -0.01 & \textit{(0.06)} & 0.38 & \textit{(0.13)} & 0.00 & \textit{(0.05)} & 0.09 & \textit{(0.02)} & [55$|$79$|$103]\\
			~ & 0.75 & 1.11 & 20.26 & \textit{(3.35)} & -0.13 & \textit{(0.08)} & 0.41 & \textit{(0.15)} & 0.04 & \textit{(0.06)} & 0.10 & \textit{(0.02)} & [38$|$53$|$71]\\
			~ & 0.90 & 1.11 & 9.95 & \textit{(1.07)} & -0.23 & \textit{(0.10)} & 0.52 & \textit{(0.17)} & 0.06 & \textit{(0.08)} & 0.09 & \textit{(0.03)} & [14$|$24$|$40]\\
			~ & 0.95 & 1.11 & 5.52 & \textit{(0.79)} & -0.27 & \textit{(0.13)} & 0.58 & \textit{(0.18)} & 0.06 & \textit{(0.09)} & 0.10 & \textit{(0.04)} & [6$|$15$|$27]\\
			\multicolumn{2}{c}{~} & \multicolumn{12}{c}{~}\\[-0.3cm]
			\hline\hline\end{tabular}
	}		
	\end{center}
	\caption{\footnotesize Simulation results for the quantile LASSO performance for the model in \eqref{sim1} for various quantiles levels and sample sizes based on 1000 Monte Carlo repetitions. Two models are always considered: the first one uses the prior knowledge that there are two change-points in the true model (the corresponding regularization parameter is denoted as $\lambda_{(2)}$) and the second model is based on the asymptotically appropriate value  $\lambda_{AS} = C (n^{-1} \log n)^{1/2}$. The estimation bias and the Mean Squared Error (MSE) quantity are provided with the corresponding standard errors. For the model with the regularization parameter $\lambda_{AS}$ we also provide an information about the estimated number of change-points (''[M$|$M$|$M]'' stands for the minimum, median, and maximum number of change-points estimated over 1000 Monte Carlo simulations). The model with $\lambda_{(2)}$ always contains two change-points and thus, three segments.}
	\label{tab1}
\end{table}

The proposed quantile LASSO method is also compared with the standard LASSO approach and the SMUCE estimator. The quantile LASSO is used to estimate a stepwise conditional median function while the standard LASSO approach and the SMUCE method are estimating the conditional mean value instead. However, the error distributions are all symmetric and, therefore, a mutual comparison of these three methods is quite straightforward. The results are summarized in Tables \ref{tab2} and \ref{tab3}.

\begin{figure}
	\centering
	\subfigure[Standard LASSSO $|$ Distribution $N(0,1)$] {\label{fig3:a}\includegraphics[width=0.48\textwidth, height = 0.25\textwidth]{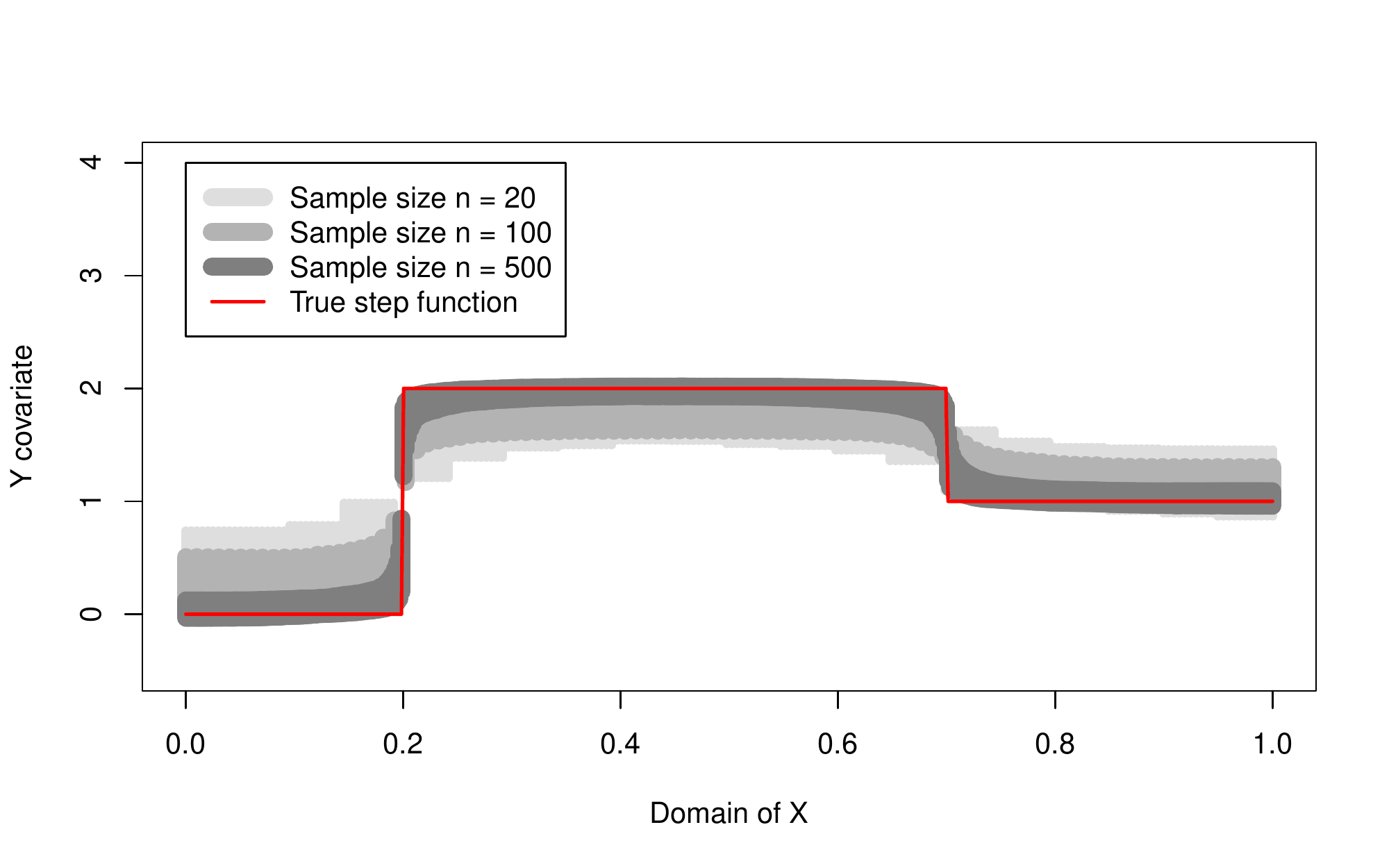}}
	\subfigure[Standard LASSSO $|$ Distribution $C(0,1)$] {\label{fig3:b}\includegraphics[width=0.48\textwidth, height = 0.25\textwidth]{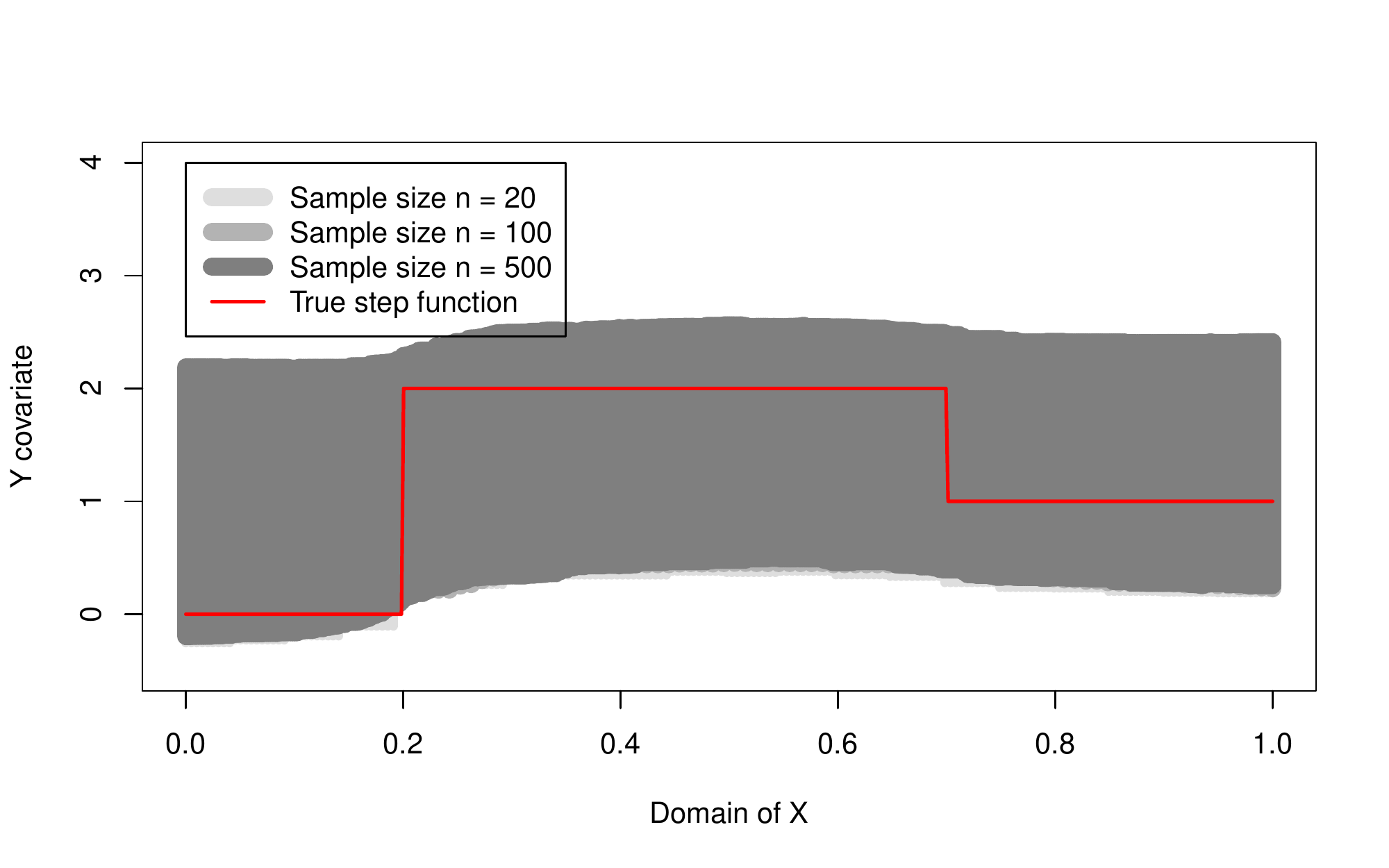}}
	
	\subfigure[SMUCE Method $|$ Distribution $N(0,1)$] {\label{fig3:a}\includegraphics[width=0.48\textwidth, height = 0.25\textwidth]{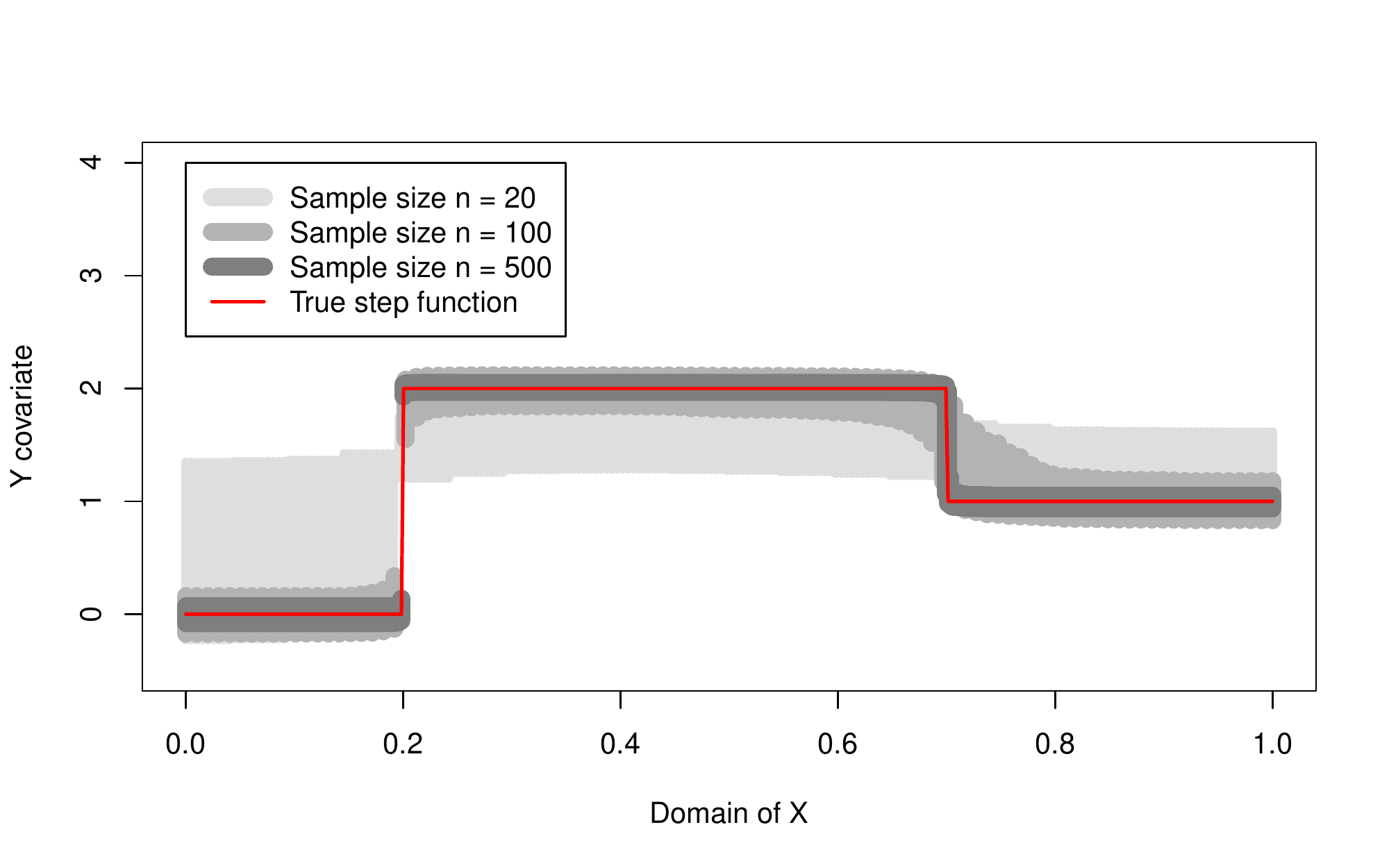}}
	\subfigure[SMUCE Method $|$ Distribution $C(0,1)$] {\label{fig3:b}\includegraphics[width=0.48\textwidth, height = 0.25\textwidth]{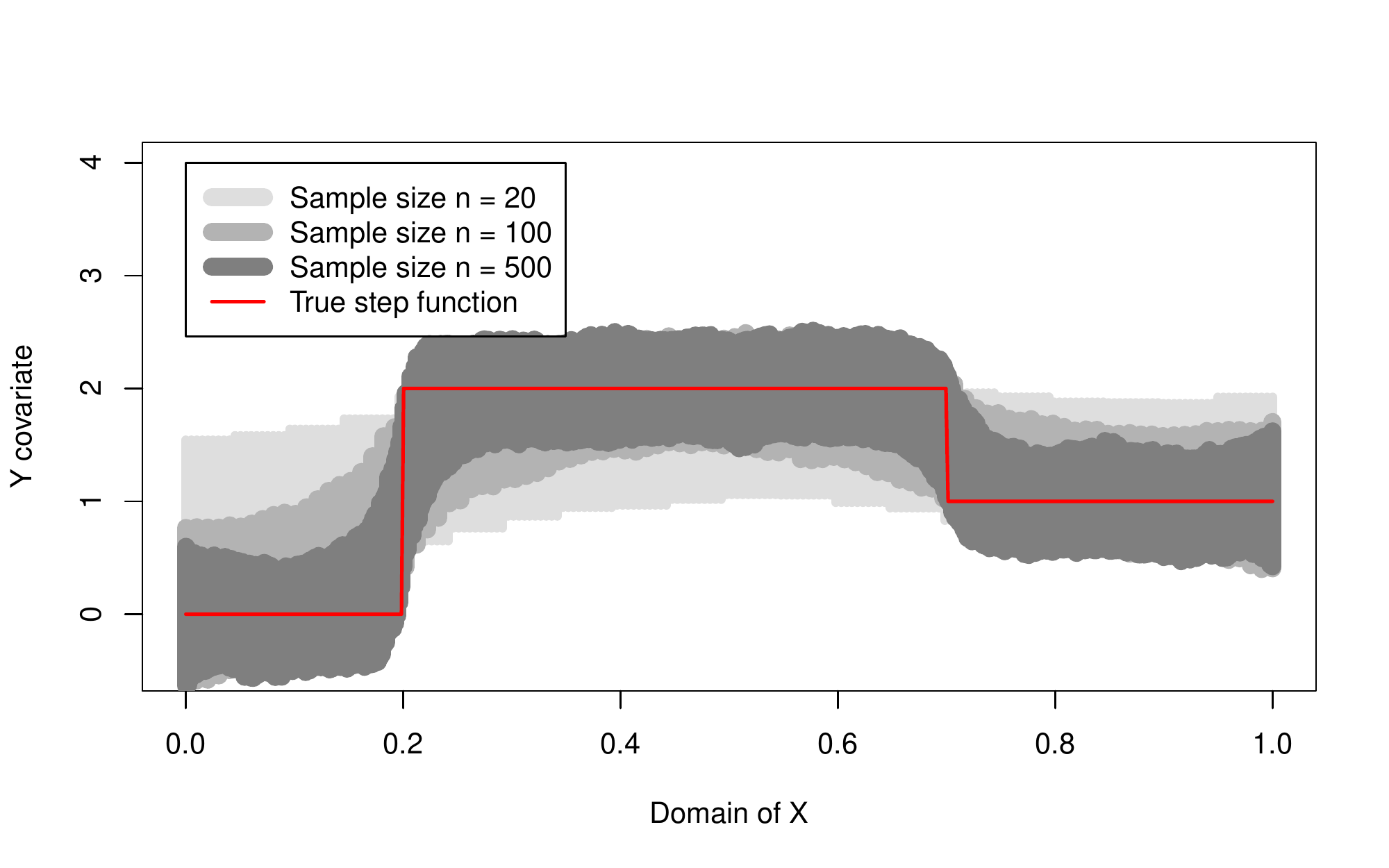}}
	
	\subfigure[Quantile LASSO $|$ Distribution $N(0,1)$] {\label{fig3:a}\includegraphics[width=0.48\textwidth, height = 0.25\textwidth]{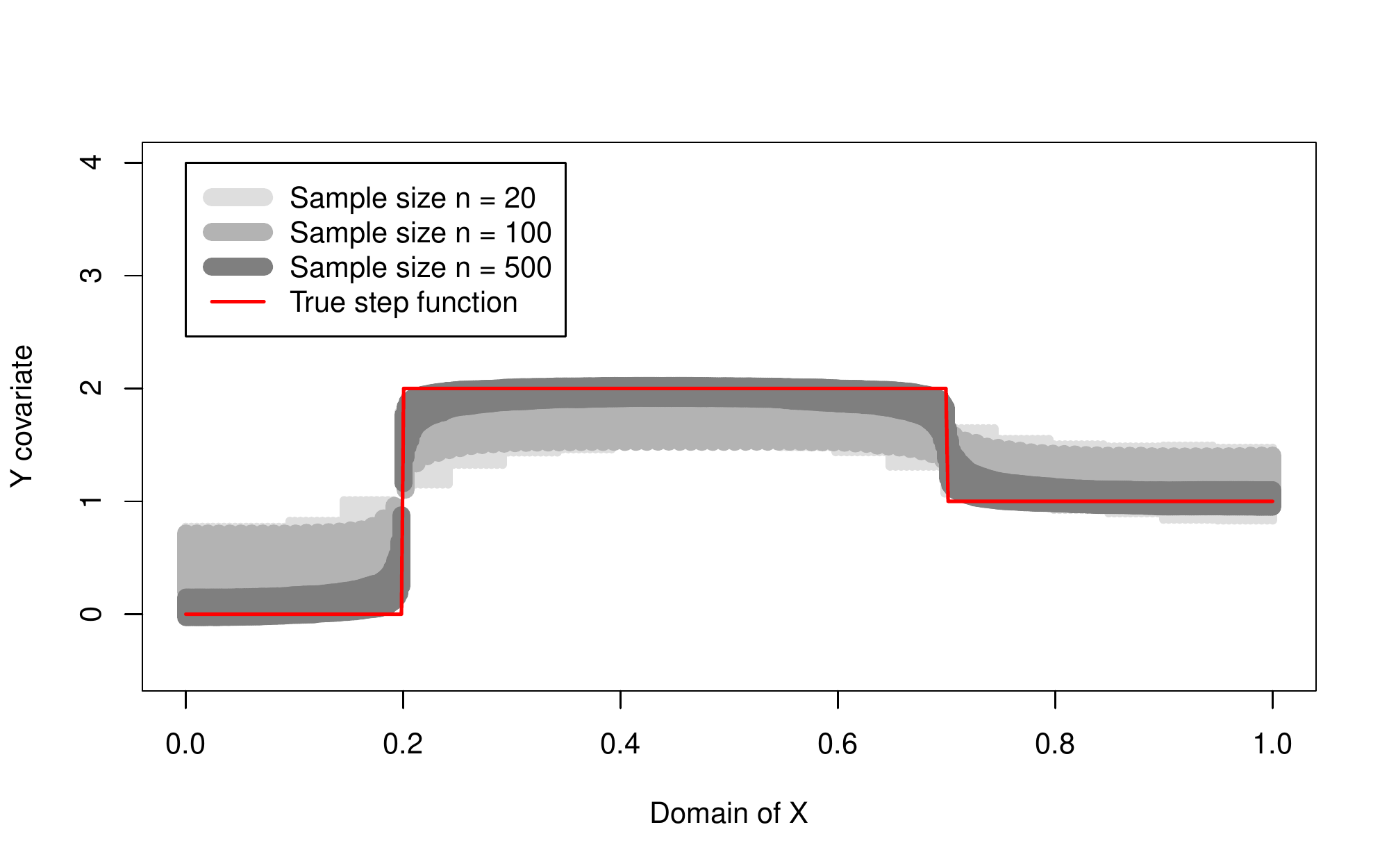}}
	\subfigure[Quantile LASSO $|$ Distribution $C(0,1)$] {\label{fig3:b}\includegraphics[width=0.48\textwidth, height = 0.25\textwidth]{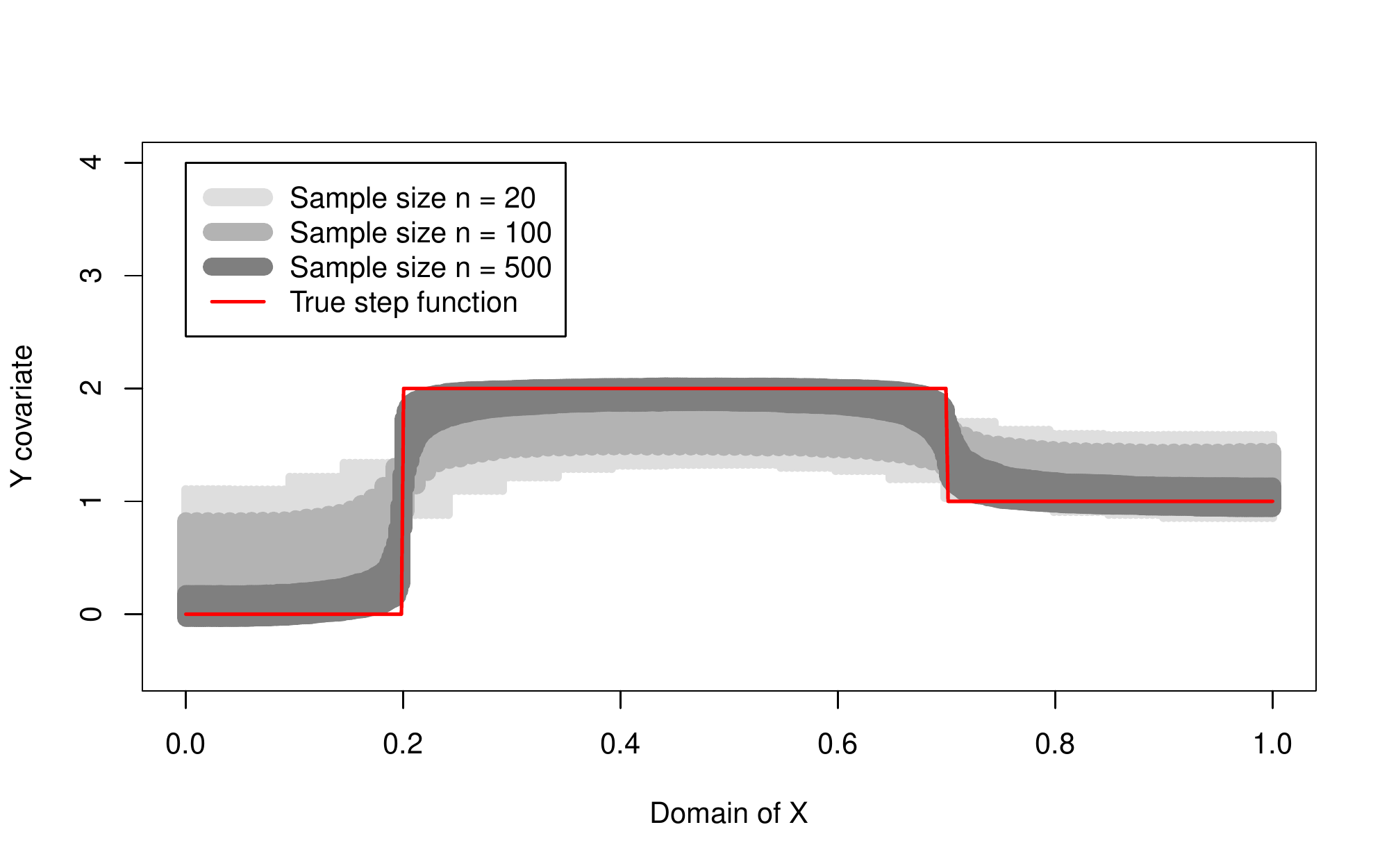}}
	
	\caption{\footnotesize Robustness and asymptotic performance of the standard LASSO method, SMUCE method and the proposed quantile LASSO estimator. The empirical point-wise interquartile ranges based on 1000 Monte Carlo repetitions are provided for three sample sizes $n \in \{20, 100, 500\}$ and two random error distributions (standard normal and Cauchy distribution). 
	The LASSO estimates are obtained for $\lambda_{MS}$	which minimizes the mean squared error term $\frac{1}{n} \sum_{i = 1}^{n} (\widehat{u}_i^* - u_{i}^*)^2$. The true underlying step function is plotted as a solid red line.   }
	\label{fig3}
\end{figure}

The performance of all three methods is very similar for normally distributed random errors, but the quantile LASSO (denoted as QLASSO) clearly outperforms the standard LASSO (denoted as SLASSO) and SMUCE in case of heavy tailed error distributions (Student's $t$ distribution and Cauchy distribution). The robust flavor of the quantile LASSO is evident at the first glance: while the quantile LASSO performs quite reasonably and a proper convergence is observed for all scenarios the standard LASSO fails for other than normal distributions---the estimation bias seems to increase with larger sample sizes and  the corresponding variance terms literally explode. Thus, no convergence can be observed for the standard LASSO estimates. The SMUCE method performs better than the standard LASSO but, it is still outperformed by the quantile LASSO for heavy tailed distributions (see Figure \ref{fig3}).

The reason why we observe such differences in the reported MSE values among the three models with the standard LASSO approach for heavy tailed distributions in Table \ref{tab2} can be understood when considering also Table \ref{tab3}. The standard LASSO models with $\lambda_{AS}$ and $\lambda_{MS}$ heavily overfit the data with respect to the number of detected change-points and therefore, the bias and variance terms are little suppressed by the huge number of change-points which are present in the model. The quantile LASSO, however, seems to perform more reasonably even for the heavy-tailed error distributions (median of the number of detected change-points is roughly 2 for the quantile LASSO, but the number of change-points for the standard LASSO is very unstable as it can range from zero up to the maximum number of parameters)---see Table \ref{tab3} for more details.  The vector of sparse parameters is more augmented for the models with $\lambda_{AS}$ and $\lambda_{MS}$ and therefore, the reported bias terms are slightly smaller than for the model with $\lambda_{(2)}$. The robust nature of the proposed quantile LASSO  approach can be also nicely visualized in Figure \ref{fig3}. The difference between the estimation performance with respect to the conditional median/mean of the quantile LASSO, standard LASSO, and the SMUCE method is obvious. While all three methods perform roughly at the same quality for the normally distributed error terms, the quantile LASSO only can handle heavy-tailed distributions---the Cauchy distribution in particular. Unlike the conditional median estimate produced by the quantile LASSO, the conditional mean estimates produced by the standard LASSO approach and SMUCE are totally unrealistic (with huge bias and variability and also too high and unstable number of estimated change-points).

\begin{table}[!ht]\footnotesize
	\begin{center}
		\scalebox{0.83}{
		    \renewcommand{\arraystretch}{0.5}
			\begin{tabular}{ccc|r<{\hspace{-\tabcolsep}}>{\hspace{-\tabcolsep}\,}lr<{\hspace{-\tabcolsep}}>{\hspace{-\tabcolsep}\,}lr<{\hspace{-\tabcolsep}}>{\hspace{-\tabcolsep}\,}l
					r<{\hspace{-\tabcolsep}}>{\hspace{-\tabcolsep}\,}lr<{\hspace{-\tabcolsep}}>{\hspace{-\tabcolsep}\,}lr<{\hspace{-\tabcolsep}}>{\hspace{-\tabcolsep}\,}l}
				\multirow{2}{*}{$\boldsymbol{\mathcal{D}}$} & & \multirow{2}{*}{$\boldsymbol{n}$} & 
				\multicolumn{4}{c}{\textbf{Model with $\boldsymbol{\lambda_{(2)}}$}} & \multicolumn{4}{c}{\textbf{Model with $\boldsymbol{\lambda_{AS}}$}} & \multicolumn{4}{c}{\textbf{Model w. $\boldsymbol{\lambda_{MS}} /$  SMUCE}} \\
				~ & ~ & ~ & 
				\multicolumn{2}{c}{\scalebox{1}{Est. Bias}} & \multicolumn{2}{c}{MSE} & \multicolumn{2}{c}{\scalebox{1}{Est. Bias}} & \multicolumn{2}{c}{MSE} & \multicolumn{2}{c}{\scalebox{1}{Est. Bias}} & \multicolumn{2}{c}{MSE}  \\\hline\hline
				\multicolumn{3}{c}{~} & \multicolumn{12}{c}{~}\\
				$\boldsymbol{N}$ & \multirow{3}{*}{\rotatebox{90}{\scalebox{0.9}{SLasso}}} & \textbf{20} & 0.00 & \textit{(0.23)} & 0.43 & \textit{(0.16)} & 0.00 & \textit{(0.23)} & 0.65 & \textit{(0.08)} & 0.00 & \textit{(0.23)} & 0.27 & \textit{(0.15)} \\
				& & \textbf{100} & 0.00 & \textit{(0.10)} & 0.35 & \textit{(0.12)} & 0.00 & \textit{(0.10)} & 0.09 & \textit{(0.04)} & 0.00 & \textit{(0.10)} & 0.13 & \textit{(0.08)} \\
				& & \textbf{500} & 0.00 & \textit{(0.04)} & 0.33 & \textit{(0.11)} & 0.00 & \textit{(0.04)} & 0.06 & \textit{(0.02)} & 0.00 & \textit{(0.04)} & 0.02 & \textit{(0.01)} \\
				\multicolumn{3}{c}{~} & \multicolumn{12}{c}{~}\\
				\rowcolor{Gray} &  & \textbf{20} & -0.01 & \textit{(0.31)} & 0.53 & \textit{(0.21)} & -0.04 & \textit{(0.29)} & 0.70 & \textit{(0.13)} & 0.00 & \textit{(0.24)} & 0.31 & \textit{(0.17)} \\
				\rowcolor{Gray} & & \textbf{100} & -0.01 & \textit{(0.14)} & 0.43 & \textit{(0.14)} & 0.00 & \textit{(0.13)} & 0.12 & \textit{(0.06)} & 0.00 & \textit{(0.13)} & 0.20 & \textit{(0.11)} \\
				\rowcolor{Gray} & \multirow{-4.5}{*}{\rotatebox{90}{\scalebox{0.9}{QLasso}}} & \textbf{500} & -0.01 & \textit{(0.06)} & 0.40 & \textit{(0.12)} & 0.00 & \textit{(0.05)} & 0.09 & \textit{(0.02)} & 0.00 & \textit{(0.05)} & 0.03 & \textit{(0.01)} \\
				\multicolumn{3}{c}{~} & \multicolumn{12}{c}{~}\\
				& \multirow{3}{*}{\rotatebox{90}{\scalebox{0.9}{~~SMUCE}}} & \textbf{20} &  ~  & \textit{~} & ~ & \textit{~} &  ~ & \textit{~} & ~ & \textit{~} & -0.01 & \textit{(0.23)} & 0.55 & \textit{(0.21)} \\
				& & \textbf{100} &  ~  & \textit{~} & ~ & \textit{~} &  ~ & \textit{~} & ~ & \textit{~} & 0.01 & \textit{(0.10)} & 0.12 & \textit{(0.10)} \\
				& & \textbf{500} &  ~  & \textit{~} & ~ & \textit{~} &  ~ & \textit{~} & ~ & \textit{~} & 0.00 & \textit{(0.04)} & 0.02 & \textit{(0.01)} \\
				\multicolumn{3}{c}{~} & \multicolumn{12}{c}{~}\\
				$\boldsymbol{t_3}$ & \multirow{3}{*}{\rotatebox{90}{\scalebox{0.9}{SLasso}}} & \textbf{20} & 0.00 & \textit{(0.38)} & 0.61 & \textit{(0.46)} & 0.00 & \textit{(0.38)} & 0.73 & \textit{(0.36)} & 0.00 & \textit{(0.38)} & 0.47 & \textit{(0.31)} \\
				& & \textbf{100} & 0.00 & \textit{(0.17)} & 0.41 & \textit{(0.14)} & 0.00 & \textit{(0.17)} & 0.33 & \textit{(0.81)} & 0.00 & \textit{(0.17)} & 0.20 & \textit{(0.11)} \\
				& & \textbf{500} & 0.00 & \textit{(0.08)} & 0.36 & \textit{(0.11)} & 0.00 & \textit{(0.08)} & 0.69 & \textit{(2.32)} & 0.00 & \textit{(0.08)} & 0.05 & \textit{(0.03)} \\
				\multicolumn{3}{c}{~} & \multicolumn{12}{c}{~}\\
				\rowcolor{Gray} &  & \textbf{20} & -0.01 & \textit{(0.35)} & 0.59 & \textit{(0.28)} & -0.04 & \textit{(0.34)} & 0.73 & \textit{(0.16)} & 0.01 & \textit{(0.28)} & 0.39 & \textit{(0.21)} \\
				\rowcolor{Gray} & & \textbf{100} & -0.01 & \textit{(0.15)} & 0.44 & \textit{(0.14)} & 0.00 & \textit{(0.14)} & 0.15 & \textit{(0.08)} & 0.00 & \textit{(0.14)} & 0.24 & \textit{(0.12)} \\
				\rowcolor{Gray} & \multirow{-4.5}{*}{\rotatebox{90}{\scalebox{0.9}{QLasso}}} & \textbf{500} & -0.02 & \textit{(0.07)} & 0.42 & \textit{(0.13)} & 0.00 & \textit{(0.06)} & 0.11 & \textit{(0.03)} & 0.00 & \textit{(0.06)} & 0.04 & \textit{(0.02)} \\
				\multicolumn{3}{c}{~} & \multicolumn{12}{c}{~}\\
				& \multirow{3}{*}{\rotatebox{90}{\scalebox{0.9}{SMUCE}}} & \textbf{20} &  ~  & \textit{~} & ~ & \textit{~} &  ~ & \textit{~} & ~ & \textit{~} & -0.01 & \textit{(0.39)} & 1.33 & \textit{(2.06)} \\
				& & \textbf{100} &  ~  & \textit{~} & ~ & \textit{~} &  ~ & \textit{~} & ~ & \textit{~} & 0.01 & \textit{(0.19)} & 1.08 & \textit{(1.60)} \\
				& & \textbf{500} &  ~  & \textit{~} & ~ & \textit{~} &  ~ & \textit{~} & ~ & \textit{~} & 0.00 & \textit{(0.09)} & 0.95 & \textit{(2.52)} \\
				\multicolumn{3}{c}{~} & \multicolumn{12}{c}{~}\\
				$\boldsymbol{C}$ & \multirow{3}{*}{\rotatebox{90}{\scalebox{0.9}{SLasso}}} & \textbf{20} & -1.57 & \textit{(74.19)} & 5867 & \textit{(154600)} & -1.57 & \textit{(74.19)} & 109736 & \textit{(3089261)} & -1.57 & \textit{(74.19)} & 5785 & \textit{(154589)} \\
				& & \textbf{100} & -1.15 & \textit{(22.21)} & 524 & \textit{(7333)} & -1.15 & \textit{(22.21)} & 50589 & \textit{(765657)} & -1.15 & \textit{(22.21)} & 512 & \textit{(7322)} \\
				& & \textbf{500} & -1.65 & \textit{(36.73)} & 1354 & \textit{(24483)} & -1.65 & \textit{(36.72)} & 671194 & \textit{(12245127)} & -1.65 & \textit{(36.72)} & 1353 & \textit{(24477)} \\
				\multicolumn{3}{c}{~} & \multicolumn{12}{c}{~}\\
				\rowcolor{Gray} &  & \textbf{20} & -0.02 & \textit{(0.46)} & 0.75 & \textit{(0.48)} & -0.03 & \textit{(0.44)} & 0.81 & \textit{(0.30)} & 0.01 & \textit{(0.36)} & 0.53 & \textit{(0.30)} \\
				\rowcolor{Gray} & & \textbf{100} & -0.02 & \textit{(0.20)} & 0.49 & \textit{(0.16)} & -0.02 & \textit{(0.18)} & 0.20 & \textit{(0.12)} & -0.01 & \textit{(0.18)} & 0.28 & \textit{(0.15)} \\
				\rowcolor{Gray} & \multirow{-4.5}{*}{\rotatebox{90}{\scalebox{0.9}{QLasso}}} & \textbf{500} & -0.02 & \textit{(0.09)} & 0.44 & \textit{(0.14)} & 0.00 & \textit{(0.08)} & 0.18 & \textit{(0.06)} & 0.00 & \textit{(0.07)} & 0.05 & \textit{(0.03)} \\
				\multicolumn{3}{c}{~} & \multicolumn{12}{c}{~}\\
				& \multirow{3}{*}{\rotatebox{90}{\scalebox{0.9}{SMUCE}}} & \textbf{20} &  ~  & \textit{~} & ~ & \textit{~} &  ~ & \textit{~} & ~ & \textit{~} & -1.58 & \textit{(74.19)} & 109953 & \textit{(3091453)} \\
				& & \textbf{100} &  ~  & \textit{~} & ~ & \textit{~} &  ~ & \textit{~} & ~ & \textit{~} & -1.16 & \textit{(22.21)} & 50683 & \textit{(766016)} \\
				& & \textbf{500} &  ~  & \textit{~} & ~ & \textit{~} &  ~ & \textit{~} & ~ & \textit{~} & -1.65 & \textit{(36.72)} & 671259 & \textit{(12245434)} \\
				\multicolumn{3}{c}{~} & \multicolumn{12}{c}{~}\\
				\hline\hline\end{tabular}
}		
	\end{center}
	\caption{\footnotesize Comparison of the quantile LASSO performance (QLasso) with the standard LASSO approach (SLasso) and the SMUCE method. The results are given for the model in \eqref{sim1} for three different (symmetric) error distributions with various signal-to-noise ratio ($N \equiv N(0,1)$,  $t_{3} \equiv$  Student's distribution with three degrees of freedom, and finally, $C \equiv Cauchy(0,1)$) and various sample sizes. Three models are considered: the model with the true number of change-points with the corresponding regularization parameter $\lambda_{(2)}$, the model with the asymptotically appropriate value $\lambda_{AS} = C \cdot (n^{-1} \log n)^{1/2}$, and the model with $\lambda_{MS}$ given by minimizing the mean squared error $\sum_{i = 1}^n (\widehat{u}_i^* - u_{i}^*)^2$. The reported values are given with the corresponding standard errors over 1000 Monte Carlo simulations.}
	\label{tab2}
\end{table}

\begin{table}[!ht]\footnotesize
	\begin{center}
		\scalebox{0.83}{
		     \renewcommand{\arraystretch}{0.5}
			\begin{tabular}{ccc|cccccr<{\hspace{-\tabcolsep}}>{\hspace{-\tabcolsep}\,}lr<{\hspace{-\tabcolsep}}>{\hspace{-\tabcolsep}\,}lc<{\hspace{-\tabcolsep}}>{\hspace{-\tabcolsep}\,}l}
				\multirow{2}{*}{$\boldsymbol{\mathcal{D}}$} & & \multirow{2}{*}{$\boldsymbol{n}$} & 
				$\boldsymbol{\lambda_{AS}}$ & \multicolumn{1}{c}{$\boldsymbol{\lambda_{(2)}}$} & \multicolumn{1}{c}{$\boldsymbol{\lambda_{MS}}$} & \multicolumn{2}{c}{Number of Jumps $\boldsymbol{|\widehat{\cal A}_n|}$} & 
				\multicolumn{6}{c}{\textbf{Change-point Detection Error}} \\
				~ & ~ & ~ & Value & Avg. & Avg. & $\lambda_{AS}$ & $\lambda_{MS}$/\scalebox{0.8}{SMUCE} & \multicolumn{2}{c}{(with $\lambda_{(2)}$)} & \multicolumn{2}{c}{(with $\lambda_{AS}$)} & \multicolumn{2}{c}{($\lambda_{CV}$/\scalebox{0.8}{SMUCE})} \\\hline\hline
				\multicolumn{3}{c}{~} & \multicolumn{11}{c}{~}\\
				$\boldsymbol{N}$ & \multirow{3}{*}{\rotatebox{90}{\scalebox{0.9}{SLasso}}} & \textbf{20} & 3.87 & 1.68 & 0.98 & [0$|$0$|$3] & [0$|$0$|$11] & 0.08 & \textit{(0.06)} & 0.17 & \textit{(0.03)} & 0.04 & \textit{(0.04)} \\
				&  & \textbf{100} & 2.15 & 7.77 & 3.29 & [2$|$2$|$14] & [2$|$2$|$21] & 0.02 & \textit{(0.02)} & 0.01 & \textit{(0.01)} & 0.01 & \textit{(0.01)} \\
				&  & \textbf{500} & 1.11 & 38.09 & 4.54 & [32$|$32$|$66] & [3$|$3$|$23] & 0.00 & \textit{(0.00)} & 0.00 & \textit{(0.00)} & 0.00 & \textit{(0.00)} \\
				\multicolumn{3}{c}{~} & \multicolumn{11}{c}{~}\\
				\rowcolor{Gray} & ~ & \textbf{20} & 3.87 & 1.58 & 1.17 & [0$|$0$|$0] & [0$|$0$|$13] & 0.10 & \textit{(0.07)} & NaN & \textit{(NA)} & 0.04 & \textit{(0.04)} \\
				\rowcolor{Gray}  &  & \textbf{100} & 2.15 & 6.54 & 4.00 & [2$|$2$|$21] & [1$|$1$|$28] & 0.03 & \textit{(0.04)} & 0.01 & \textit{(0.01)} & 0.01 & \textit{(0.03)} \\
				\rowcolor{Gray}  & \multirow{-4.5}{*}{\rotatebox{90}{\scalebox{0.9}{QLasso}}} & \textbf{500} & 1.11 & 31.31 & 4.68 & [35$|$35$|$70] & [4$|$4$|$65] & 0.00 & \textit{(0.00)} & 0.00 & \textit{(0.00)} & 0.00 & \textit{(0.00)} \\
				\multicolumn{3}{c}{~} & \multicolumn{11}{c}{~}\\
				& \multirow{3}{*}{\rotatebox{90}{\scalebox{0.9}{SMUCE}}} & \textbf{20} &  ~ & ~ & ~ &  ~ & [0$|$0$|$4] &  ~ & ~ &  ~ & ~ & 0.07 & \textit{(0.04)} \\
				&  & \textbf{100} &  ~ & ~ & ~ &  ~ & [0$|$1$|$4] &  ~ & ~ &  ~ & ~ & 0.02 & \textit{(0.03)} \\
				&  & \textbf{500} &  ~ & ~ & ~ &  ~ & [2$|$2$|$3] &  ~ & ~ &  ~ & ~ & 0.00 & \textit{(0.00)} \\
				\multicolumn{3}{c}{~} & \multicolumn{11}{c}{~}\\
				$\boldsymbol{t_3}$ & \multirow{3}{*}{\rotatebox{90}{\scalebox{0.9}{SLasso}}} & \textbf{20} & 3.87 & 2.14 & 2.25 & [0$|$0$|$5] & [0$|$0$|$9] & 0.10 & \textit{(0.06)} & 0.12 & \textit{(0.06)} & 0.06 & \textit{(0.05)} \\
				&  & \textbf{100} & 2.15 & 8.86 & 4.02 & [3$|$3$|$23] & [1$|$1$|$14] & 0.05 & \textit{(0.05)} & 0.02 & \textit{(0.02)} & 0.02 & \textit{(0.03)} \\
				&  & \textbf{500} & 1.11 & 41.10 & 7.61 & [67$|$67$|$124] & [3$|$3$|$23] & 0.01 & \textit{(0.02)} & 0.00 & \textit{(0.00)} & 0.00 & \textit{(0.00)} \\
				\multicolumn{3}{c}{~} & \multicolumn{11}{c}{~}\\
				\rowcolor{Gray} & ~ & \textbf{20} & 3.87 & 1.52 & 1.55 & [0$|$0$|$0] & [0$|$0$|$13] & 0.11 & \textit{(0.06)} & NaN & \textit{(NA)} & 0.05 & \textit{(0.05)} \\
				\rowcolor{Gray}  &  & \textbf{100} & 2.15 & 6.03 & 3.91 & [2$|$2$|$20] & [0$|$0$|$26] & 0.04 & \textit{(0.05)} & 0.01 & \textit{(0.02)} & 0.02 & \textit{(0.04)} \\
				\rowcolor{Gray}  & \multirow{-4.5}{*}{\rotatebox{90}{\scalebox{0.9}{QLasso}}} & \textbf{500} & 1.11 & 28.98 & 4.59 & [54$|$54$|$104] & [4$|$4$|$36] & 0.01 & \textit{(0.01)} & 0.00 & \textit{(0.00)} & 0.00 & \textit{(0.00)} \\
				\multicolumn{3}{c}{~} & \multicolumn{11}{c}{~}\\
				& \multirow{3}{*}{\rotatebox{90}{\scalebox{0.9}{SMUCE}}} & \textbf{20} &  ~ & ~ & ~ &  ~ & [0$|$0$|$6] &  ~ & ~ &  ~ & ~ & 0.09 & \textit{(0.05)} \\
				&  & \textbf{100} &  ~ & ~ & ~ &  ~ & [0$|$2$|$10] &  ~ & ~ &  ~ & ~ & 0.05 & \textit{(0.05)} \\
				&  & \textbf{500} &  ~ & ~ & ~ &  ~ & [2$|$2$|$23] &  ~ & ~ &  ~ & ~ & 0.02 & \textit{(0.02)} \\
				\multicolumn{3}{c}{~} & \multicolumn{11}{c}{~}\\
				$\boldsymbol{C}$ & \multirow{3}{*}{\rotatebox{90}{\scalebox{0.9}{SLasso}}} & \textbf{20} & 3.87 & 112.40 & 83.71 & [0$|$0$|$17] & [0$|$1$|$19] & 0.12 & \textit{(0.06)} & 0.08 & \textit{(0.06)} & 0.10 & \textit{(0.06)} \\
				&  & \textbf{100} & 2.15 & 386.53 & 237.47 & [7$|$7$|$96] & [0$|$0$|$96] & 0.12 & \textit{(0.06)} & 0.01 & \textit{(0.01)} & 0.07 & \textit{(0.06)} \\
				&  & \textbf{500} & 1.11 & 1414.77 & 1353.09 & [164$|$164$|$493] & [0$|$0$|$499] & 0.10 & \textit{(0.06)} & 0.00 & \textit{(0.00)} & 0.08 & \textit{(0.06)} \\
				\multicolumn{3}{c}{~} & \multicolumn{11}{c}{~}\\
				\rowcolor{Gray} & ~ & \textbf{20} & 3.87 & 1.45 & 2.66 & [0$|$0$|$1] & [0$|$0$|$12] & 0.12 & \textit{(0.06)} & 0.07 & \textit{(NA)} & 0.06 & \textit{(0.05)} \\
				\rowcolor{Gray}  &  & \textbf{100} & 2.15 & 5.27 & 3.65 & [2$|$2$|$22] & [0$|$0$|$39] & 0.05 & \textit{(0.05)} & 0.02 & \textit{(0.02)} & 0.03 & \textit{(0.04)} \\
				\rowcolor{Gray}  & \multirow{-4.5}{*}{\rotatebox{90}{\scalebox{0.9}{QLasso}}} & \textbf{500} & 1.11 & 24.58 & 4.86 & [54$|$54$|$192] & [2$|$2$|$42] & 0.01 & \textit{(0.01)} & 0.00 & \textit{(0.00)} & 0.00 & \textit{(0.00)} \\
				\multicolumn{3}{c}{~} & \multicolumn{11}{c}{~}\\
				& \multirow{3}{*}{\rotatebox{90}{\scalebox{0.9}{SMUCE}}} & \textbf{20} &  ~ & ~ & ~ &  ~ & [0$|$0$|$7] &  ~ & ~ &  ~ & ~ & 0.10 & \textit{(0.05)} \\
				&  & \textbf{100} &  ~ & ~ & ~ &  ~ & [1$|$1$|$22] &  ~ & ~ &  ~ & ~ & 0.05 & \textit{(0.04)} \\
				&  & \textbf{500} &  ~ & ~ & ~ &  ~ & [25$|$25$|$72] &  ~ & ~ &  ~ & ~ & 0.01 & \textit{(0.01)} \\
				\multicolumn{3}{c}{~} & \multicolumn{11}{c}{~}\\
				\hline\hline\end{tabular}
}		
	\end{center}
	\caption{\footnotesize Comparison of the quantile LASSO performance (QLasso) with the standard LASSO approach (SLasso) and the SMUCE method.  Three different (symmetric) error distributions are considered ($N \equiv N(0,1)$, $t_{3} \equiv$  Student's distribution with three degrees of freedom, and finally, $C \equiv Cauchy(0,1)$) and the number of estimated change-points (where ''[m$|$m$|$m]'' stands for the minimum, median, and maximum number of change-points estimated over 1000 Monte Carlo simulations) and the change-point detection rate together with the corresponding standard errors are provided. The models with three different values of $\lambda_n > 0$ are considered: the model with $\lambda_{(2)}$, the model with $\lambda_{AS}$, and the model with $\lambda_{MS}$. The change-point detection rate is calculated only for models where at lest two change-points were discovered, otherwise NA values are reported.}
	\label{tab3}
\end{table}

Moreover, the same can be also told about the change-point detection performance. If we use the prior knowledge that two change-points (three segments) are supposed to be estimated then all three methods  perform quite well if the error terms are normally distributed but, for the Cauchy distribution, the detection of the standard LASSO and SMUCE approach is way aside from the true change-points locations. The quantile LASSO, however, can still provide a valid detection.

The behavior of the quantile LASSO estimator which can be observed in the simulation results is, indeed, in a concordance with the theoretical results proved in Section \ref{results} and the common knowledge of the LASSO performance. The LASSO penalty, in general, is well-known for recovering usually more non-zero coefficients than really needed---this is also confirmed by the simulation study. Secondly, the estimated parameters are always shrunk towards zero and thus, the estimates tend to underestimate the underlying structure, introducing a systematic bias, which is also observed in the simulation study.

\section{Conclusion and Final Remarks}
\label{Conclusion}
In this paper we proposed the quantile LASSO estimator and we investigated its main theoretical and empirical properties. The quantile LASSO is robust with respect to outlying observations and heavy-tailed random error distributions: it clearly outperforms the standard LASSO method in both---the estimation of the unknown underlying structure and, also, in detection of the unknown change-point locations (both under the heavy-tailed error distributions). 

From the theoretical point of view, the main advantage of the proposed method lies in much weaker  distributional assumptions: the quantile LASSO performance does not rely on any normal or sub-Gaussian distributions which are typically required for the standard LASSO approach and, moreover, much complex insight into the data can be obtained by estimating an arbitrary quantile rather than the mean value only. Another convenient property of the proposed method is that instead of proving its oracle properties or sign consistency results and thus, requiring strong assumptions for the design matrix, we rather show the performance with respect to the change-point detection and therefore, only some mild assumptions are required and the method, in general, is much widely applicable.

The proposed simulations study confirms the theoretical results and it markedly emphasizes the robust nature of the quantile LASSO estimator.

\section*{Acknowledgement}
	The work was partially supported by a bilateral grant between France and the Czech Republic provided by the PHC Barrande 2017 grant of Campus France (CG, grant number 38105NM) and the Ministry of Educations, Youth, and Sports in the Czech Republic (MM, Mobility grant 7AMB17FR030).

%\section*{References}

%\begin{appendices}

\renewcommand\thesection{\Alph{section}}
\setcounter{section}{0}

 \section{APPENDIX: Proofs}
\label{Append}
  \subsection{Auxiliary lemmas and their proofs}
 In this section we state three important lemmas which are crucial for proving the results from Section \ref{results}. The first lemma  is a direct consequence of the  Karush-Kuhn-Tucker (KKT) optimality conditions. It is useful not only for proving the asymptotic behavior of the change-point number estimator, but also for deriving the properties of the change-points location estimators given by
the sequence $\widehat{t}_{1} < \widehat{t}_{2} < \dots < \widehat{t}_{|\widehat{\cal A}_n|}$. 
\begin{lem}
\label{Lemma 3}
		For the model described in \eqref{eq1} and any solution $\widehat{\eb^n} \in \mathbb{R}^n$ of the minimization problem in \eqref{eq5}, it holds, with probability one, for any $n \in \mathbb{N}$ and $\lambda_{n} > 0$, that
	\begin{equation}
	\label{eq8}
	\tau (n -\widehat t_l) - \sum^n_{i=\widehat t_l} \e1_{\{Y_i < \widehat u_i\}}= n \lambda_n \widehat \alpha_l, \qquad \forall l \in \{1, \cdots , |\widehat{\cal A}_n | \},
	\end{equation}
	and 
	\begin{equation}
	\label{eq9}
	\left|\tau (n -j) - \sum^n_{i=j} \e1_{\{Y_i < \widehat u_i\}} \right|\leq  n \lambda_n  , \qquad \forall j \in \{1, \cdots , n \},
	\end{equation}
	with
	$
	\widehat \alpha_l \equiv  
	\e1_{\{\widehat u_{\widehat t_l} > \widehat u_{\widehat t_l- 1}\}}-\e1_{\{\widehat u_{\widehat t_l} \leq \widehat u_{\widehat t_l- 1} \}} 
	$.
\end{lem}

\noindent {\bf Proof of Lemma \ref{Lemma 3}.}\\
By the Karush-Kuhn-Tucker (KKT) optimality conditions, we have, for    all $j \in \widehat{\cal A}_n$, that
\[
\tau\sum^n_{i=1} X_{ij} - \sum^n_{i=1} X_{ij} \e1_{\big\{Y_i < (\eeX_n \widehat{\eb^n})_i \big\}} =n \lambda_n sign(\widehat \beta_j).
\]
Taking into account the form of $\eeX_n$, we obtain the relation in \eqref{eq8}. Similarly we also obtain the relation in \eqref{eq9}.
\hspace*{\fill}$\blacksquare$ 

\begin{lem}
\label{Lemma}
Let $A$ and $B$ be two random variables and $x>0$ is some positive real value such that $\PP[|A+B| \leq x]=1$, then, for any constant  $v>1$ we have that
\[
1 \leq \PP\bigg[x \geq \frac{|A|}{v} \bigg] +\PP\bigg[|B| \geq \frac{v-1}{v}|A| \bigg].
\] 
\end{lem}

\noindent {\bf Proof of Lemma \ref{Lemma}.}\\
Obviously, it holds that $1=\PP\big[x \geq \frac{|A|}{v} \big]+\PP\big[x < \frac{|A|}{v} \big]$. The inequality $|A+B| \geq |A|-|B|$ implies that:  $\PP[|A+B| \leq x] \leq \PP[|A|-|B| \leq x] $. Then, by using the fact that$\PP[|A+B| \leq x]=1$, we can write: $\PP\big[x < \frac{|A|}{v} \big] = \PP\big[\big\{x < \frac{|A|}{v} \big\} \cap \big\{|A+B| \leq x \big\} \big] \leq \PP\big[\big\{x < \frac{|A|}{v} \big\} \cap \big\{|A|-|B| \leq x \big\} \big]=\PP \big[|B| \geq \frac{v-1}{v}|A| \big]$ and the lemma follows. 
\hspace*{\fill}$\blacksquare$ 

\begin{lem}
Let $\{v_n\}$ and $\{x_n\}$ be two positives sequences such that $v_n x^2_n (\log n)^{-1} {\underset{n \rightarrow \infty}{\longrightarrow} } \infty$. Then, under Assumption (A2) imposed for error terms $\{\varepsilon_i\}_{1 \leqslant i \leqslant n}$, we have
\[
\PP\bigg[ \max_{
\substack{
1 \leq r_n< s_n \leq n\\
s_n-r_n \geq v_n
}} 
\;
\sup_{t \in \R}
\bigg| \frac{1}{s_n-r_n} \sum ^{s_n-1}_{i=r_n} \e1_{\{\varepsilon_i \leq t\}}-F(t) \bigg| \geq x_n
 \bigg]  {\underset{n \rightarrow \infty}{\longrightarrow} } 0,
\]
for $F$ being the distribution function of $\varepsilon_i$.
\label{Lemma 4}
\end{lem}

\noindent {\bf Proof of Lemma \ref{Lemma 4}.}\\
Firstly, we have that 
\begin{align}
\label{tx1}
\PP\bigg[ \max_{\substack{1 \leq r_n< s_n \leq n\\s_n-r_n \geq v_n}} &\; \sup_{t \in \R} \bigg| \frac{1}{s_n-r_n} \sum^{s_n-1}_{i=r_n} \e1_{\{\varepsilon_i \leq t\}}-  F(t) \bigg| \geq x_n  \bigg]\nonumber\\
& \leq \sum_{\substack{1 \leq r_n< s_n \leq n\\s_n-r_n \geq v_n}} \PP\bigg[\sup_{t \in \R} \bigg| \frac{1}{s_n-r_n} \sum^{s_n-1}_{i=r_n} \e1_{\varepsilon_i \leq t}-F(t) \bigg| \geq x_n  \bigg].
\end{align}
By  Dvoretzky-Kiefer-Wolfowitz's  inequality (see \cite{Dvoretzky.Kiefer.Wolfowitz.56}) for the independent  Bernoulli random variables $ \e1_{\{\varepsilon_i \leq t\}}$, we obtain for all $\epsilon >0$, that
\[
\PP\bigg[\sup_{t \in \R} \bigg| \frac{1}{s_n-r_n} \sum^{s_n-1}_{i=r_n} \e1_{\varepsilon_i \leq t}-F(t) \bigg| \geq x_n  \bigg] \leq 2 \exp\big(-2(s_n-r_n)x^2_n \big).
\]
Then, taking into account (\ref{tx1}) and the fact that $v_n x^2_n (\log n)^{-1} {\underset{n \rightarrow \infty}{\longrightarrow} } \infty$, we also obtain that 
\[
\PP\bigg[ \max_{\substack{1 \leq r_n< s_n \leq n\\s_n-r_n \geq v_n}} \; \sup_{t \in \R} \bigg| \frac{1}{s_n-r_n} \sum^{s_n-1}_{i=r_n} \e1_{\varepsilon_i \leq t}-F(t) \bigg| \geq x_n  \bigg]  \leq 2 n^2 \exp\big(-2(s_n-r_n)x^2_n \big) {\underset{n \rightarrow \infty}{\longrightarrow} } 0,
\]
which proofs the assertion of Lemma \ref{Lemma 4}.
\hspace*{\fill}$\blacksquare$

\subsection{Proofs of Theorems}
In this Section we proof the main results formulated in the three theorems in Section \ref{results}.\\

\noindent {\bf Proof of Theorem  \ref{Proposition 5}.}\\
Let us start by defining two random events, for any $k=1, \cdots , K^*$:
\[
V_{n,k} \equiv \big\{|\widehat t_k -t^*_k | \geq n \delta_n \big\} 
\qquad \textrm{and} \qquad
W_n \equiv \Big\{ \max_{1 \leqslant j \leqslant K^*} |\widehat t_j -t^*_j | <\frac{I^*_{min}}{2} \Big\}.
\]
By Assumption (A5), since  $K^*< \infty$, the theorem is proved if we show that for any $k =1, \cdots , K^*$, it holds that
\begin{equation}
\label{PVnk}
\lim_{n \rightarrow \infty} \PP[V_{n,k} ]=0.
\end{equation}
In order to prove the relation in (\ref{PVnk}), we suppose that random event $V_{n,k}$ occurs.  The event $V_{n,k}$, for any $k=1, \cdots , K^*$, can be also expressed as $V_{n,k}= \big(V_{n,k} \cap W_n\big) \cup \big(V_{n,k} \cap \overline W_n\big)$, with $\overline W_n$ being the complementary event of $W_n$. \\
If event $V_{n,k}$ occurs, without any loss of generality, we can assume that $t^*_k-\widehat t_k \geq [n \delta_n]$. The opposite case for $\widehat t_k -t^*_k \geq [n \delta_n]$ follows similarly. We now consider two steps for proving (\ref{PVnk}): firstly,  we study $\PP[V_{n,k} \cap W_n]$ and, later, we focus on $\PP[V_{n,k} \cap \overline W_n]$.

\noindent\textbf{Step 1.} We will show that for any $k=1, \cdots , K^*$, it holds that
\begin{equation}
\label{eq13}
\lim_{n \rightarrow \infty} \PP[V_{n,k} \cap W_n] =0.
\end{equation}
Let us start by considering the relation in \eqref{eq9}, for $j=t^*_k$,
\begin{equation}
\label{eq11}
\bigg| \tau(n- t^*_k) -\sum^n_{i=t^*_k} \e1_{\{Y_i \leq \widehat u_i\}} \bigg| \leq n \lambda_n,
\end{equation}
and the relation in \eqref{eq8}, for $l=k$:
\[
\tau(n- t^*_k +t^*_k -\widehat t_k) -\bigg(\sum^{t^*_k-1}_{i=\widehat t_k}+\sum^n_{i=t^*_k} \bigg)\e1_{\{Y_i \leq \widehat u_i\}} = n \lambda_n \widehat \alpha_k,
\]
where we assume, without any loss of generality, that $t^*_{k-1} \leq \widehat t_k \leq t^*_k$. Thus, we have
\[
\tau(n - t^*_k)- \sum^n_{i=t^*_k} \e1_{\{Y_i \leq \widehat u_i\}}+\tau(t^*_k - \widehat t_k) - \sum^{t^*_k-1}_{i=\widehat t_k} \e1_{\{Y_i \leq \widehat u_i\}}= n \lambda_n \widehat \alpha_k.
\]
Next, we apply the following general result: for any $a, b, c \in \mathbb{R}$, such that 
\begin{equation}
\label{abc}
a+b=\pm c \quad \textrm{ and } \quad |b| \leq c, \quad \textrm{it holds that } |a| \leq 2c.
\end{equation}
Using \eqref{abc}  for $a=\Big[\tau(t^*_k - \widehat t_k) - \sum^{t^*_k-1}_{i=\widehat t_k} \e1_{\{Y_i \leq \widehat u_i\}}\Big]$, $b=\Big[\tau(n - t^*_k)- \sum^n_{i=t^*_k} \e1_{\{Y_i \leq \widehat u_i\}}\Big]$, and  $c= n \lambda_n$, we have that  event $U_{n,k}$ occurs with probability 1, where
\[
%\label{eq12}
U_{n,k} \equiv \bigg\{  \big| \tau(t^*_k - \widehat t_k) - \sum^{t^*_k-1}_{i=\widehat t_k} \e1_{\{Y_i \leq \widehat u_i\}}\big| \leq 2 n \lambda_n  \bigg\}.
\]

Next, we use  Lemma \ref{Lemma}, for $x=2n \lambda_n$, some constant $v$ such that $\displaystyle{v>\frac{\max(\tau, F(\mu^*_{k+1}- \mu^*_k))}{|\tau - F(\mu^*_{k+1}- \mu^*_k) |}}$ and random variables $A$ and $B$ defined as follows:
\begin{itemize}
\item if $\tau < F(\mu^*_{k+1}- \mu^*_k)$, then $A =\sum^{t^*_k-1}_{i=\widehat t_k} \e1_{\{Y_i \leq \widehat u_i\}} $, $B= \tau(t^*_k-\widehat t_k)$;
\item  if $\tau > F(\mu^*_{k+1}- \mu^*_k)$, then $A= \tau(t^*_k-\widehat t_k)$,
 $B=\sum^{t^*_k-1}_{i=\widehat t_k} \e1_{\{Y_i \leq \widehat u_i\}} $.
\end{itemize}
Then, for the probability $\PP [V_{n,k} \cap W_n]$ we obtain
\begin{align}
\hskip-0.5cm\PP [V_{n,k} \cap W_n] &= \PP [V_{n,k} \cap W_n \cap U_{n,k} ]  \nonumber\\
& =\PP \bigg[\bigg\{  \big|\tau(t^*_k - \widehat t_k) - \sum^{t^*_k-1}_{i=\widehat t_k} \e1_{\{Y_i \leq \widehat u_i\}} \big| \leq 2 n \lambda_n  \bigg\} \cap  V_{n,k} \cap W_n \bigg] \nonumber\\
%& \leq \PP \bigg[\bigg\{\tau|t^*_k-\widehat t_k| \leq n \lambda_n  \bigg\} \cap  V_{n,k} \cap W_n \bigg] + \PP %\bigg[\bigg\{\sum^{t^*_k-1}_{i=\widehat t_k} \e1_{\{Y_i \leq \widehat u_i\}} \leq n \lambda_n  \bigg\} \cap  V_{n,k} \cap W_n %\bigg]\nonumber\\
& \textcolor{red}{\leq } {\cal P}_1 + {\cal P}_2, \label{eq_mm1}
\end{align}
with ${\cal P}_1 \equiv \PP \big[ \big\{ \frac{|A|}{v} \leq x \big\} \cap V_{n,k} \cap W_n \big] $ and ${\cal P}_2 \equiv \PP \big[ \big(|B| \geq  \frac{v-1}{v} |A| \big) \cap V_{n,k} \cap W_n \big] $ and we distinguish for two individual cases where we either have $\tau < F(\mu^*_{k+1}- \mu^*_k)$, or $\tau > F(\mu^*_{k+1}- \mu^*_k)$.

We start with the situation for which  $\tau < F(\mu^*_{k+1}- \mu^*_k)$. We consider the first term in \eqref{eq_mm1} where we have ${\cal P}_1= \PP \bigg[ \bigg\{ \sum^{t^*_k-1}_{i=\widehat t_k} \e1_{\{\varepsilon_i +\mu^*_k < \widehat \mu_{k+1}\}} \leq 2 n \lambda_n v \bigg\} \cap  V_{n,k} \cap W_n \bigg]$, with the constant $v$, such that  $v>\frac{F(\mu^*_{k+1}- \mu^*_k)}{F(\mu^*_{k+1}- \mu^*_k) - \tau}$. 

Under Assumptions (A2), (A5), and (A6), by applying Theorem 2 of \cite{Fan.Fan.Barut.14},  we obtain that the relation in (\ref{Bel3}) holds. Then, we have that 
\begin{equation}
\label{T21}
 \widehat \mu_{k+1} -\mu^*_{k+1} =O_{\PP} \left( \sqrt{\frac{\log n}{n}}\right),
\end{equation}
which implies that  there exists a constant $c_1 >0$ not depending on $n$, that 
\begin{equation}
\label{num}
|\widehat \mu_{k+1} - \mu^*_{k+1}| \leq c_1\sqrt{\frac{\log n}{n}}, \qquad \textrm{with probability converging to 1.}
\end{equation}

Next, we recall two general results, which are needed to complete the proof:
\begin{description}
	\item \textit{(i)} Let $X$ and $Z$ be two  real random variables and $x \in \R$. Then the following holds:
	\begin{description}
	\item \textit{(i1)} If  $\PP[Z \geq x]=1$, then $\e1_{\{X \leq Z\}} \geq \e1_{\{X \leq x\}}$ with probability 1.
	\item \textit{(i2)}  If  $\PP[Z \leq x]=1$, then $\e1_{\{X \leq Z\}} \leq \e1_{\{X \leq x\}}$  with probability 1.
	\end{description}
		\item \textit{(ii)}  Let  $S_1$ and $S_2$ be two real random variables such that $S_1 \leq S_2$ with probability one. Then for any $x \in \mathbb{R}$ we have that $\PP[S_1 \leq x] \geq \PP[S_2 \leq x]$.
\end{description}
Using now relation  \textit{(i1)}  together with (\ref{num}), we have, with probability converging to 1,  \[\e1_{\big\{\varepsilon+\mu^*_k-\mu^*_{k+1} \leq \widehat{\mu}_{k+1} - \mu^*_{k+1}\big\}} \geq \e1_{\big\{\varepsilon+\mu^*_k-\mu^*_{k+1} \leq - c_1\sqrt{\frac{  \log n}{n}}\big\}} . \]
Using this last inequality  together with \textit{(ii)}, we obtain for ${\cal P}_1$, that
\begin{eqnarray}
{\cal P}_1 & = &  \PP \bigg[ \bigg\{ \sum^{t^*_k-1}_{i=\widehat t_k} \e1_{\{\varepsilon_i +\mu^*_k < \widehat \mu_{k+1}\}} \leq 2 n \lambda_n v \bigg\} \cap  V_{n,k} \cap W_n \bigg]    \nonumber \\
& \leq &   \PP \bigg[ \bigg\{ \sum^{t^*_k-1}_{i=\widehat t_k}\e1_{\big\{\varepsilon_i+\mu^*_k-\mu^*_{k+1} \leq - c_1\sqrt{\frac{\log n}{n}}\big\}} \leq 2 n \lambda_n v \bigg\} \cap  V_{n,k} \cap W_n \bigg]+o(1)  \nonumber \\
& \leq &   \PP \bigg[ \bigg\{ \sum^{t^*_k-1}_{i=t^*_k - [n \delta_n]}\e1_{\big\{\varepsilon_i+\mu^*_k-\mu^*_{k+1} \leq - c_1\sqrt{\frac{ \log n}{n}}\big\}} \leq 2 n  \lambda_n v \bigg\} \cap  V_{n,k} \cap W_n \bigg]+o(1), \nonumber
\end{eqnarray}
where for the  last inequality we used the fact that $t^*_k - \widehat t_k$ must the smallest possible value, that is $(t^*_k - [n \delta_n])$, with $[n \delta_n]$ being the integer part of $n \delta_n$.

By the random events $V_{n, k}$ and $W_n$ we have that $n \delta_n < t^*_k - \widehat t_k < {I^*_{min}}/{2}$. By Assumption (A5) and the  Strong Law of Large Numbers for independent $\varepsilon_i$, we obtain
\[
\frac{1}{[n \delta_n]} \sum^{t^*_k-1}_{i=t^*_k-[n \delta_n]} \left( \e1_{\big\{\varepsilon_i \leq \mu^*_{k+1} - \mu^*_k -c_1 \sqrt{\frac{K^* \log n}{n}}\big\}}- F\Big(\mu^*_{k+1} - \mu^*_k - c_1 \sqrt{\frac{ \log n}{n}}\Big)  \right) \overset{a.s.} {\underset{n \rightarrow \infty}{\longrightarrow}} 0.
\]
Since by Assumption (A2) we have $F(x)>0$ for all $x \in \R$,  there exists a constant $ C>0$, such that $F\big(\mu^*_{k+1} - \mu^*_k - c_1 \sqrt{\frac{ \log n}{n}} \big)   > C$. Thus, there also exists a positive constant $\tilde C > 0$, such that
\[
\sum^{t^*_k-1}_{i=t^*_k-[n \delta_n]} \e1_{\big\{\varepsilon_i \leq \mu^*_{k+1} - \mu^*_k - c_1 \sqrt{\frac{\log n}{n}}\big\}} \geq \widetilde{C} [n \delta_n],
\] 
with probability converging to one as $n$ tends to infinity. Taking into account Assumption (A4), we finally get
\begin{equation}
\label{e16P2}
{\cal P}_1 \leq \PP \cro{{\widetilde{C}} [n \delta_n] \leq \sum^{t^*_k-1}_{i=t^*_k-[n \delta_n]} \e1_{\big\{\varepsilon_i \leq \mu^*_{k+1} - \mu^*_k - c_1 \sqrt{\frac{ \log n}{n}}\big\}} \leq n \lambda_n  } + o(1) \textcolor{violet}{\underset{n \rightarrow \infty}{\longrightarrow} 0.}
\end{equation}
Analogously, for ${\cal P}_2= \PP \Big[ \big(\tau (t^*_k - \widehat t_k) \geq  \frac{v-1}{v} \sum^{t^*_k-1}_{i=\widehat t_k} \e1_{\{\varepsilon_i \leq \widehat \mu_{k+1}- \mu^*_k\}} \big) \cap V_{n,k} \cap W_n \Big]$, we have
\[
{\cal P}_2 \leq  \PP \left[ \Big\{\frac{1}{t^*_k - \widehat t_k}   \sum^{t^*_k-1}_{i=\widehat t_k} \e1_{\big\{\varepsilon_i \leq  \mu^*_{k+1}- \mu^*_k   - c_1 \sqrt{\frac{ \log n}{n}}\big\}} \leq \frac{\tau v}{v-1} \Big\} \cap V_{n,k} \cap W_n \right].
\]
Using Lemma \ref{Lemma 4} for $v_n=[n \delta_n]$ and $x_n=\big| \frac{\tau v}{v-1} - F(\mu^*_{k+1}- \mu^*_k) \big|$,  and due to Assumption (A4) where we have  $\delta_n/\lambda_n \rightarrow \infty$, we get that 
\[
\frac{n \delta_n}{\log n}= \frac{n \delta_n}{\sqrt{\log n}}  \cdot \frac{\sqrt{n}}{\sqrt{n \log n}}= \frac{n \delta_n}{n \lambda_n}\cdot \sqrt{\frac{n}{\log n}} \rightarrow \infty.
\]  
Thus, we  obtain
\begin{align*}
\PP \bigg[ \max_{t^*_k-\widehat t_k \geq [n \delta_n]} \bigg| \frac{1}{t^*_k - \widehat t_k}   \sum^{t^*_k-1}_{i=\widehat t_k} \e1_{\Big\{\varepsilon_i \leq  \mu^*_{k+1}- \mu^*_k   - c_1 \sqrt{\frac{ \log n}{n}}\Big\}}  & -  F\Big(\mu^*_{k+1}- \mu^*_k-  c_1 \sqrt{\frac{ \log n}{n}}\Big) \bigg| \\
& \geq \bigg| \frac{\tau v}{v-1} - F(\mu^*_{k+1}- \mu^*_k) \bigg| \bigg] \underset{n \rightarrow \infty}{\longrightarrow}  0.
\end{align*}
We used the fact that $ F\Big(\mu^*_{k+1}- \mu^*_k-  c_1 \sqrt{\frac{ \log n}{n}}\Big) \rightarrow F(\mu^*_{k+1}- \mu^*_k)$ as $n$ converges to infinity, and   ${\tau v}(v-1)^{-1} - F(\mu^*_{k+1}- \mu^*_k) < 0$. Finally, we have
\begin{align}
\label{e16P1}
{\cal P}_2  \leq \PP \bigg[\max_{t^*_k-\widehat t_k \geq [n \delta_n]} \bigg( \frac{1}{t^*_k - \widehat t_k}   \sum^{t^*_k-1}_{i=\widehat t_k} &\e1_{\big\{\varepsilon_i \leq  \mu^*_{k+1}- \mu^*_k   - c_1 \sqrt{\frac{ \log n}{n}}\big\}}  -  F(\mu^*_{k+1}- \mu^*_k) \bigg) \\
& < \frac{\tau v}{v-1} - F(\mu^*_{k+1}- \mu^*_k)  \bigg] \underset{n \rightarrow \infty}{\longrightarrow}  0,\nonumber
\end{align} 
and combining relations \eqref{eq_mm1},  \eqref{e16P2}, and \eqref{e16P1}, we get that \eqref{eq13} holds true.

Let us now focus on the second case, where $\tau > F(\mu^*_{k+1}- \mu^*_k)$. For  ${\cal P}_1$, we can write
\begin{align*}
\label{e16P2b}
{\cal P}_1 & = \PP \bigg[ \bigg\{  \tau(t^*_k- \widehat t_k) \leq 2 n \lambda_n v \bigg\} \cap  V_{n,k} \cap W_n \bigg] \\
& \leq \PP[n \delta_n \tau \leq \tau(t^*_k- \widehat t_k) \leq 2 n \lambda_n]  \underset{n \rightarrow \infty}{\longrightarrow}  0.
\end{align*}
So, we only need to  deal with ${\cal P}_2= \PP \big[ \big\{ \frac{1}{t^*_k- \widehat t_k} \sum^{t^*_k-1}_{i=\widehat t_k} \e1_{\{\varepsilon_i \leq \widehat \mu_{k+1}- \mu^*_k\}} \geq \frac{v-1}{v} \tau \big\} \cap V_{n,k} \cap W_n \big]$. Using  \textit{(i2)}  together with \eqref{num}, we have that \[\e1_{\big\{\varepsilon+\mu^*_k-\mu^*_{k+1} \leq \widehat{\mu}_{k+1} - \mu^*_{k+1}\big\}} \leq \e1_{\big\{\varepsilon+\mu^*_k-\mu^*_{k+1} \leq c_1\sqrt{\frac{  \log n}{n}}\big\}},  \] 
with probability converging to 1. Thus, by \textit{(ii)} we obtain
\begin{equation}
\label{PP2}
{\cal P}_2 \leq \PP \left[ \Big\{ \frac{1}{t^*_k- \widehat t_k} \sum^{t^*_k-1}_{i=\widehat t_k} \e1_{\big\{\varepsilon_i+\mu^*_k-\mu^*_{k+1} \leq c_1\sqrt{\frac{  \log n}{n}}\big\}} \geq \frac{v-1}{v} \tau \Big\} \cap V_{n,k} \cap W_n \right].
\end{equation}
By Lemma \ref{Lemma 4} we obtain
\begin{align*}
\PP \bigg[  \max_{t^*_k-\widehat t_k \geq [n \delta_n]} \bigg| \frac{1}{t^*_k-\widehat t_k} \sum^{t^*_k - 1}_{i=\widehat t_k}& \e1_{\big\{\varepsilon_i+\mu^*_k-\mu^*_{k+1} \leq c_1\sqrt{\frac{  \log n}{n}}\big\}}  - F\Big(\mu^*_{k+1} -\mu^*_k+c_1\sqrt{\frac{  \log n}{n}} \Big)  \bigg| \\
& \geq \bigg|\frac{v-1}{v}\tau -F(\mu^*_{k+1} -\mu^*_k) \bigg| \bigg]\underset{n \rightarrow \infty}{\longrightarrow}  0,
\end{align*}
and, since $\frac{v-1}{v}\tau -F(\mu^*_{k+1} -\mu^*_k) >0$, we also have
\begin{align*}
\PP \bigg[  \max_{t^*_k-\widehat t_k \geq [n \delta_n]} \bigg( \frac{1}{t^*_k-\widehat t_k} \sum^{t^*_k - 1}_{i=\widehat t_k}& \e1_{\big\{\varepsilon_i+\mu^*_k-\mu^*_{k+1} \leq c_1\sqrt{\frac{  \log n}{n}}\big\}}  - F\big(\mu^*_{k+1} -\mu^*_k \big)  \bigg)\\
& \geq  \frac{v-1}{v}\tau -F(\mu^*_{k+1} -\mu^*_k)   \bigg] \underset{n \rightarrow \infty}{\longrightarrow}  0.
\end{align*}
Combining now the last expression with \eqref{PP2} and \eqref{e16P2b}, we obtain that \eqref{eq13} holds true also for the case $\tau > F(\mu^*_{k+1}- \mu^*_k)$.

\noindent\textbf{Step 2.}
Now, we study the probability $\PP[V_{n,k} \cap\overline{W}_n]$, with $\overline W_n \equiv \left\{ \max_{1 \leqslant j \leqslant K^*} |\widehat t_j - t^*_j| > \frac{I^*_{min}}{2}\right\}$. We consider the following three random events (using the same notations as in \cite{Harchaoui.Levy.10}):
\begin{align*}
D_{n}^{(l)} & \equiv \{ \exists k \in \{ 1, \cdots ,   K^*\}; \,\, \widehat t_k \leq t^*_{k-1} \} \cap \overline{W}_n;\\
D_{n}^{(m)} & \equiv \{ \forall k \in \{ 1, \cdots ,   K^*\}; \,\, t^*_{k-1} < \widehat t_k < t^*_{k+1} \} \cap \overline{W}_n;\\
D_{n}^{(r)} &\equiv \{ \exists k \in \{ 1, \cdots ,   K^*\}; \,\, \widehat t_k \geq t^*_{k+1} \} \cap \overline{W}_n.
\end{align*} 

Then, $\PP[V_{n,k} \cap\overline{W}_n]= \PP[V_{n,k} \cap D_{n}^{(l)}] + \PP[V_{n,k} \cap D_{n}^{(m)}] + \PP[V_{n,k} \cap D_{n}^{(r)}]$  and we deal with each probability term on the right side separately. For $\PP[V_{n,k} \cap D_{n}^{(m)}]$ we have 
\begin{align*}
\PP[V_{n,k} \cap D_{n}^{(m)}] \leq & \PP[V_{n,k}\cap \{\widehat t_{k+1} - t^*_k \geq I^*_{min}/2 \} \cap D_{n}^{(m)}]\\ 
& +\sum^{K^*}_{i=k+1} \PP [\{t^*_i - \widehat t_i \geq I^*_{min}/2 \} \cap \{\widehat t_{i+1} - t^*_i \geq I^*_{min}/2 \} \cap D_n^{(m)}].
\end{align*}
We have that $t^*_{k-1}< \widehat t_k < t^*_k < \widehat t_{k+1} < t^*_{k+2}$ and applying \eqref{eq9} for $j=t^*_k$ and \eqref{eq8} for $l=k$, we obtain that 
\begin{displaymath}
\PP \left[\tau (t^*_k -\widehat t_k) \leq 2 n \lambda_n + \sum^{t^*_k}_{i=\widehat t_k} \e1_{\{\varepsilon_i \leq \widehat \mu_{k+1} -\mu^*_k\}}\right] = 1.
\end{displaymath}
On the other hand, using \eqref{eq9} for $j=t^*_k$ and \eqref{eq8} for $l=k+1$, we also get that 
\begin{displaymath}
\PP \left[\tau (\widehat t_{k+1} - t^*_k) \leq 2 n \lambda_n + \sum_{i=t^*_k}^{\widehat t_{k+1}} \e1_{\{\varepsilon_i \leq \widehat \mu_{k+1} - \mu^*_{k+1}\}}\right] = 1.
\end{displaymath}
Therefore, 
\[
 \PP[V_{n,k}\cap \{\widehat t_{k+1} - t^*_k \geq I^*_{min}/2 \} \cap D_{n}^{(m)}] \leq  \PP[\{t^*_k-\widehat t_k \geq n \delta_n  \}\cap \{\widehat t_{k+1} - t^*_k \geq I^*_{min}/2 \} \cap D_{n}^{(m)}].
\]
Since $\mu^*_k$, $\mu^*_{k+1}$ don't depend on $n$ and $\mu^*_k \neq \mu^*_{k+1}$, there is at least one of the differences $\widehat \mu_{k+1} -\mu^*_k$ or $\widehat \mu_{k+1} - \mu^*_{k+1}$ which does not converge to 0 as $n \rightarrow \infty$. Suppose it's $\widehat \mu_{k+1} -\mu^*_k$. Then  
\begin{align*}
&\PP[V_{n,k} \cap  \{\widehat t_{k+1} - t^*_k \geq I^*_{min}/2 \} \cap D_{n}^{(m)}]\\
& \leq  \PP\Big[\{t^*_k-\widehat t_k \geq n \delta_n  \}\cap \{\widehat t_{k+1} - t^*_k \geq I^*_{min}/2 \} \cap \{ \tau (t^*_k -\widehat t_k) \leq 2 n \lambda_n + \sum^{t^*_k}_{i=\widehat t_k} \e1_{\{\varepsilon_i \leq \widehat \mu_{k+1} -\mu^*_k\}}\}\Big].
\end{align*}

Similarly as in Step 1 we obtain that the last probability converges to 0 as $n\rightarrow \infty$. Analogously we can show that $\lim_{n \rightarrow \infty}\PP [\{t^*_i - \widehat t_i \geq I^*_{min}/2 \} \cap \{\widehat t_{i+1} - t^*_i \geq I^*_{min}/2 \} \cap D_n^{(m)}] =0$, for any $i=k+1, \cdots , K^*$, and since $K^*$ is bounded we obtain that  $\lim_{n \rightarrow \infty} \PP[V_{n,k} \cap D_{n}^{(m)}] =0.$\\
 For $\PP[D_n^{(l)}]$ we have, similarly as in \cite{Harchaoui.Levy.10}, that
\begin{align*}
 \PP[D_n^{(l)}] \leq 2^{K^*-1} \sum^{K^*-1}_{k=1} \sum^{K^*-1}_{m \geq k} & \PP \bigg[ \{ t^*_m -\widehat t_m >I^*_{min}/2 \} \cap  \{ \widehat t_{m+1} - t^*_m>I^*_{min}/2 \}  \bigg]\\
 & + 2^{K^*-1} \PP[t^*_{K^*} -\widehat t_{K^*} >I^*_{min}/2].
\end{align*}

Repeating the same arguments as above we can also show  that $\lim_{n \rightarrow \infty} \PP[  D_{n}^{(l)}]=0  $ and $\lim_{n \rightarrow \infty} \PP[  D_{n}^{(r)}]=0$, therefore, also $\lim_{n \rightarrow \infty} \PP[V_{n,k} \cap D_{n}^{(l)}]=0$ and $\lim_{n \rightarrow \infty} \PP[V_{n,k} \cap D_{n}^{(r)}]=0$. Putting everything together we have that $\PP \cro{V_{n,k} \cap \overline W_n }   {\underset{n \rightarrow \infty}{\longrightarrow} } 0$ which competes the proof. \hspace*{\fill}$\blacksquare$ \\

\noindent\textbf{Proof of Theorem  \ref{Proposition 6}}.\\
 In order to prove the theorem we take into account the relation in (\ref{Bel2bis}) and we study the  probability
\begin{align} 
&\PP  \bigg[\bigg\{ {\cal E} \big( \widehat{{\cal T}}_{|\widehat{\cal A}_n|} || {\cal T}^* \big) \geq n \delta_n \bigg\} \cap\left\{ K^* \leq | \widehat{\cal A}_n| \leq C_1 K^*\right\} \bigg] \leq\nonumber\\
& \leq \PP \bigg[ {\cal E} \big( \widehat{{\cal T}}_{K^*} || {\cal T}^* \big) \geq n \delta_n \bigg| |\widehat{\cal A}_n|=K^* \bigg]+\sum^{C_1K^*}_{K>K^*} \PP \bigg[ {\cal E} \big( \widehat{{\cal T}}_{K} || {\cal T}^* \big) \geq n \delta_n  \bigg| |\widehat{\cal A}_n|=K \bigg],\label{eq16}
\end{align}
where used conditional probabilities, conditioned on the number of the estimated jumps $|\widehat{\cal A}_n|$. 
For the first term in \eqref{eq16} we have 
\begin{align*}
\PP \bigg[ {\cal E} \big( \widehat{{\cal T}}_{K^*} || {\cal T}^* \big) \geq n \delta_n \bigg| |\widehat{\cal A}_n|=K^* \bigg] & \leq \PP \bigg[ \sup_{1 \leq k \leq K^*} \inf_{1 \leq j \leq K^*} | \widehat t_k - t^*_j| \geq n \delta_n \bigg] \\
& = \PP \bigg[ \sup_{1 \leq k \leq K^*} | \widehat t_k - t^*_k| \geq n \delta_n \bigg],
\end{align*}
and taking into account the assertion of Theorem \ref{Proposition 5}, we have that the last probability converges to 0 as $n \rightarrow \infty$. Therefore
\begin{equation}
\label{PPK*}
\lim_{n \rightarrow \infty}  \PP \bigg[ {\cal E} \big( \widehat{{\cal T}}_{K^*} || {\cal T}^* \big) \geq n \delta_n  \bigg| |\widehat{\cal A}_n|=K^*\bigg] =0.
\end{equation}

For the second term in \eqref{eq16} we have
\begin{equation}
\label{eq23bis}
\begin{split}
\sum^{C_1K^*}_{K>K^*} \PP \bigg[ {\cal E} \big( \widehat{{\cal T}}_{K} || {\cal T}^* \big) \geq n \delta_n \bigg| |\widehat{\cal A}_n|=K \bigg] & \leq \sum^{C_1K^*}_{K>K^*} \sum^{K^*}_{k=1}\PP \cro{ \forall 1 \leq l \leq K, |\widehat t_l - t^*_k| \geq n \delta_n} \\
& = \sum^{C_1K^*}_{K>K^*} \sum^{K^*}_{k=1} \left( \PP[E_{n,k,K,1}]+\PP[E_{n,k,K,2}]+\PP[E_{n,k,K,3}]] \right),
\end{split}
\end{equation}
where (using the same notations as in \cite{Harchaoui.Levy.10}), the random events $E_{n,k,K,1}$, $E_{n,k,K,2}$, and $E_{n,k,K,3}$ are defined as follows:
\begin{align*}
E_{n,k,K,1} & \equiv \{ \forall 1 \leq l \leq K; \,\, |\widehat t_l - t^*_k| \geq n \delta_n , \, \widehat t_l < t^*_k  \};\\
E_{n,k,K,2} & \equiv \{ \forall 1 \leq l \leq K; \,\, |\widehat t_l - t^*_k| \geq n \delta_n , \, \widehat t_l > t^*_k  \};\\
E_{n,k,K,3} & \equiv \{ \exists 1 \leq l \leq K; \,\, |\widehat t_l - t^*_k| \geq n \delta_n , \, |\widehat t_{l+1} - t^*_k| \geq n \delta_n , \, \widehat t_l < t^*_k < \widehat t_{l+1} \}.
\end{align*} 
Let us start with $\PP[E_{n,k,K,1}]$: since $\mu^*_{K^*} \neq \mu^*_{K^*+1}$, we can deduct by using the relation in \eqref{Bel3} and Assumption (A1) that for fixed $K$, the only option for $E_{n,k,K,1}$ to occur with  probability not converging to zero as $n$ goes to infinity, is for $k=K^*$ and $t^*_{K^*-1} < \widehat t_K < t^*_{K^*}$. Therefore, we study the random event $E_{n,K^*,K,1}=\bigg\{  \big\{ t^*_{K^*}- \widehat t_K >n \delta_n \big\} \cap \big\{t^*_{K^*-1} < \widehat t_K < t^*_{K^*} \big\} \cap \big\{ \widehat t_K -t^*_{K^*-1}  \geq n \delta_n\big\} \bigg\} $.

Applying now the relation in \eqref{eq9} from Lemma \ref{Lemma 3}, for $j=t^*_{K^*}$, we have
\begin{equation}
\label{eq17}
\bigg| \tau(n -t^*_{K^*})- \sum^n_{i=t^*_{K^*}} \e1_{\{Y_i < \widehat u_i\}}   \bigg| \leq n \lambda_n,
\end{equation}
and, analogously, using the relation in \eqref{eq8}, for $l=K$, we obtain
\begin{equation}
\label{eq18}
\tau(n -\widehat t_K)- \sum^n_{i=\widehat t_K} \e1_{\{Y_i < \widehat u_i\}}  = n\lambda_n \widehat \alpha_K.
\end{equation}
Next, the expression in (\ref{eq18}) can be also rewritten as
\begin{equation}
\label{eq19}
\bigg[ \tau(n-t^*_{K^*}) - \sum^n_{i=t^*_{K^*}}\e1_{\{Y_i < \widehat u_i\}} \bigg]+\bigg[ \tau(t^*_{K^*}-\widehat t_K) - \sum^{t^*_{K^*}-1}_{i=\widehat t_K}\e1_{\{Y_i < \widehat u_i\}}  \bigg] = n\lambda_n \widehat \alpha_K,
\end{equation}
and we can use  the property already given in \eqref{abc} for  $a=\tau(t^*_{K^*}-\widehat t_K) - \sum^{t^*_{K^*}-1}_{i=\widehat t_K}\e1_{\{Y_i < \widehat u_i\}}$, $b= \tau(n-t^*_{K^*}) -  \sum^n_{i=t^*_{K^*}}\e1_{\{Y_i < \widehat u_i\}}$, and $c=n \lambda_n$. Then, taking into account \eqref{eq17}, \eqref{eq19}, and \eqref{abc}, we have,  with probability one, that  
\begin{equation}
\label{eq20}
\bigg| \tau(t^*_{K^*}-\widehat t_K) - \sum^{t^*_{K^*}-1}_{i=\widehat t_K}\e1_{\{\varepsilon_i < \widehat \mu_{K+1} - \mu^*_{K^*}\}} \bigg| \leq 2 n \lambda_n. 
\end{equation}
By Lemma \ref{Lemma}, for $x=2n \lambda_n$, some constant $v$ such that $\displaystyle{v>\frac{\max(\tau, F(\mu^*_{K^*+1}- \mu^*_{K^*}))}{|\tau - F(\mu^*_{K^*+1}- \mu^*_{K^*}) |}}$, and random variables $A$ and $B$ defined as 
\begin{itemize}
\item $A =\sum^{t^*_{K^*}-1}_{i=\widehat t_K} \e1_{\{Y_i \leq \widehat u_i\}} $ and $B= \tau(t^*_{K^*}-\widehat t_K)$, if $\tau < F(\mu^*_{K^*+1}- \mu^*_{K^*})$,
\item $A= \tau(t^*_{K^*}-\widehat t_K)$ and $B=\sum^{t^*_{K^*}-1}_{i=\widehat t_K} \e1_{\{Y_i \leq \widehat u_i\}} $, if $\tau > F(\mu^*_{K^*+1}- \mu^*_{K^*})$,  
\end{itemize}
we have,
\begin{align*}
\PP[E_{n,K^*,K,1}] & \textcolor{red}{\leq } \PP\Big[E_{n,K^*,K,1} \cap \Big\{\frac{|A|}{v} \leq x \Big\}  \Big] + \PP\Big[E_{n,K^*,K,1} \cap  \Big\{|B| \geq  \frac{v-1}{v} |A| \Big\}\Big]\\
& \equiv {\cal P}_{1,K^*}+{\cal P}_{2,K^*}.
\end{align*}
To show that ${\cal P}_{1,K^*} \underset{n \rightarrow \infty}{\longrightarrow}  0$ we can use the same idea as 
for the probability in \eqref{e16P2} and, similarly, to show that ${\cal P}_{2,K^*} \underset{n \rightarrow \infty}{\longrightarrow}  0$, we use the same principle as in  \eqref{e16P1}. Finally, by Assumption (A5), we have
\begin{equation}
\label{eq25}
\lim_{n \rightarrow \infty} \sum^{C_1 K^*}_{K>K^*}  \PP \cro{ E_{n,K^*,K,1}} =0.
\end{equation}

Next, we consider $E_{n,k,K,2}$: again, the only option for $E_{n,k,K,2}$ to occur with some probability not converging to zero, is for $k=1$ and $t^*_{1} < \widehat t_1 < t^*_{2}$. Therefore, we only need to focus on $E_{n,1,K,2}=\bigg\{  \big\{ t^*_2- \widehat t_1 >n \delta_n \big\} \cap \big\{t^*_1 < \widehat t_1 < t^*_2 \big\} \cap \big\{ \widehat t_1 -t^*_{1}  \geq n \delta_n\big\} \bigg\} $.
Applying  Lemma \ref{Lemma 3} for $j=t^*_{2}$ and $l=1$, we obtain, same as before, that
\begin{equation}
\label{PE1}
\lim_{n \rightarrow \infty} \sum^{C_1 K^*}_{K>K^*}   \PP[E_{n,1,K,2}]=0.
\end{equation}
Finally, we deal with $E_{n,k,K,3}$. We can apply Lemma \ref{Lemma 3} for the same indexes $j$ and $l$ as in the proof of Proposition 4 in \cite{Harchaoui.Levy.10}, and by following the same idea as above we get that 
\begin{equation}
\label{PE2}
\lim_{n \rightarrow \infty} \sum^{C_1 K^*}_{K>K^*} \sum^{K^*}_{k=1}  \PP[E_{n,k,3}]=0.
\end{equation}
Using now the relations in \eqref{eq25}, \eqref{PE1}, and (\ref{PE2}), taking also into account the expression in  \eqref{eq23bis}, we obtain that
\[
\lim_{n \rightarrow \infty}\sum^{C_1K^*}_{K>K^*} \PP \bigg[ {\cal E} \big( \widehat{{\cal T}}_{K} || {\cal T}^* \big) \geq n \delta_n \bigg| |\widehat{\cal A}_n|=K^* \bigg]=0,
\]
which, together with \eqref{PPK*} and  \eqref{eq16},  gives
\[
\lim_{n \rightarrow \infty}  \PP \bigg[\bigg\{ {\cal E} \big( \widehat{{\cal T}}_{|\widehat{\cal A}_n|} || {\cal T}^* \big) \geq n \delta_n \bigg\} \cap\left\{ K^* \leq | \widehat{\cal A}_n| \leq C K^*\right\} \bigg] =0,
\]
which also implies the relation in \eqref{eq15}.
\hspace*{\fill}$\blacksquare$ \\

\noindent {\bf Proof of Theorem  \ref{theorem 3.3}}.\\ 
Let $\widehat{{\cal T}}_{\widehat K} \equiv \{ \widehat t_1, \cdots , \widehat t_{\widehat K}  \}$ be the set of the change-point location the estimates by the quantile LASSO method, such that $|\widehat{\cal A}_n | =\widehat{K}$. Let us consider two quantile processes
\begin{align*}
Q(\emu(K),K) & \equiv \sum^n_{i=1} \rho_\tau (Y_i-u_i)=\sum^{K+1}_{k=1} \sum^{t_k-1}_{i=t_{k-1}} \rho_\tau (Y_i - \mu_k),\\
S(\emu(K),K) & \equiv Q(\emu(K),K)  + n \lambda_n \sum^{n-1}_{i=1} |u_{i+1}- u_i  | \\
& = \sum^n_{i=1} \rho_\tau (Y_i-(\eeX_n \eb)_i) +n \lambda_n \sum^{n}_{i=2} |\beta_i  |,
\end{align*}
where $\boldsymbol{\mu}(K) = (\mu_1, \dots, \mu_{K + 1})^\top$ for some $K \in \mathbb{N}$ fixed and $(\beta_i)_{1 \leqslant i \leqslant n}$, and $(u_i)_{1 \leqslant i \leqslant n}$ being defined in Section 2.  Let us  define the quantile LASSO estimator of the $(K+1)$-dimensional vector $\emu(K)$ and of the cahnge-point number $K$, as
\[
\big(\widehat \emu(\widehat{K}), \widehat{K}\big)  \equiv \argmin_{\emu \in \mathbb{R}^{K + 1}, K \in \mathbb{N}} \quad S(\emu(K),K),
\]
with $\widehat \emu(\widehat{K})=\big(\widehat{\mu}_{1}, \cdots , \widehat{\mu}_{\widehat{K}+1} \big)^\top$ obtained by estimating   $K$ and $\emu(K)$ simultaneously. Let us also define another estimator for the same vector $\emu(K)=(\mu_1 , \cdots , \mu_{K+1})^\top$, however, for some $K \in \mathbb{N}$ fixed, defined as
\[
\overset{\vee}{\emu}(K) \equiv \argmin_{\emu \in \R^{K+1}} \quad S(\emu(K),K) , 
\]
where $\overset{\vee}{\emu}(K)  = \big(\overset{\vee}{ \mu}_1, \cdots, \overset{\vee}{ \mu}_{K+1}  \big)^\top$. The assertion of the theorem will be proved if we show that under the supposition that $\widehat{K} < K^*$ we have
\begin{equation}
\label{eq28}
\PP \cro{S(\widehat \emu(\widehat K),\widehat K) > S(\overset{\vee} \emu (K^*), K^*) } {\underset{n \rightarrow \infty}{\longrightarrow}} 1.
\end{equation}
For  $\widehat K < K^*$, let us consider then the difference
\begin{align}
D & \equiv S(\widehat \emu(\widehat K),\widehat K) - S( \overset{\vee} \emu (K^*), K^*)\nonumber\\
& = \big[ S(\widehat \emu(\widehat K),\widehat K) - S(\emu^*, K^*) \big] - \big[ S( \overset{\vee} \emu (K^*), K^*) - S(\emu^*, K^*) \big]\nonumber\\
& \equiv  D_1 - D_2, \label{D}
\end{align}
where the $(K^*+1)$-vector of the true values $(\mu^*_k)_{1 \leqslant k \leqslant K^*+1}$ is denoted as $\emu^*$. 
In order to study the difference $D_1$, we can rewrite it as a sum of two terms
\begin{align*}
D_1 & = \underbrace{\sum^{\widehat K+1}_{k=1} \sum^{\widehat t_k-1}_{i=\widehat t_{k-1}} \rho_\tau(Y_i - \widehat \mu_k) - \sum^{K^*+1}_{k=1} \sum^{t^*_k-1}_{i=t^*_{k-1}} \rho_\tau(Y_i - \mu^*_k)}_{D_{1,1}} \\
& \hskip1cm + \underbrace{ n \lambda_n \cro{\sum^{\widehat K+1}_{k=1}  | \widehat \mu_{k+1}- \widehat \mu_{k}|  - \sum^{K^*+1}_{k=1} |  \mu^*_{k+1}- \mu^*_k| }}_{D_{1,2}}\\ & \equiv D_{1,1}+D_{1,2},
\end{align*}
and, similarly, the difference $D_2$ can be further rewritten as 
\begin{align*} 
D_2 & =\sum^{K^*+1}_{k=1} \sum^{t^*_k-1}_{i=t^*_{k-1}}  \rho_\tau (Y_i - \overset{\vee}{\mu}_{k+1} ) - \sum^n_{i=1} \rho_\tau(\varepsilon_i)+n \lambda_n \cro{\sum^{K^*+1}_{k=1} \pth{ |\overset{\vee}{\mu}_{k+1}- \overset{\vee}{\mu}_k| -  | \mu^*_{k+1}-\mu^*_k| }} \\
& \equiv D_{2,1}+D_{2,2}.
\end{align*}
We start by studying the difference $D_2$: firstly, we focus on $D_{2,2}$ and  afterwards on $D_{2,1}$.  Using the inequality $\big| |a|-|b|   \big| \leq |a-b|$, Assumption (A5), and the relation in (\ref{T21}), we have
\begin{equation}
\label{D2.2}
D_{2,2}  = O_{\PP} \left(  n \lambda_n \sqrt{\frac{  \log n}{n}} \right).
\end{equation}
On the other hand, from \cite{Knight.98}, we have for any $x, y \in \R$ that 
\[
% \label{Knight}
\rho_\tau(x-y)- \rho_\tau(x)=y(\e1_{\{x \leq 0\}} - \tau)+\int^y_0 (\e1_{\{x \leq t\}} -\e1_{\{x \leq 0\}})dt.
\]
Using this  relation for $x=\varepsilon_i$ and $y=C \sqrt{\frac{  \log n}{n}}$, we obtain that $D_{2,1}$ can be expressed as
\begin{align*}
D_{2,1}&=C \sqrt{\frac{  \log n}{n}}\sum^n_{i=1}(\e1_{\{\varepsilon_i < 0\}}-\tau)+\sum^{K^*+1}_{k=1} \sum^{t^*_k-1}_{i=t^*_{k-1}} \int^{ C \sqrt{\frac{ \log n}{n}}}_{0} \big[\e1_{\{\varepsilon_i < t\}} -\e1_{\{\varepsilon_i <0\}} \big] dt \\
&  \equiv D_{2,11}+D_{2,12}.
\end{align*}
By the Limit Central Theorem for i.i.d. Bernoulli random variables we obtain that $\sum^n_{i=1}(\e1_{\varepsilon_i < 0}-\tau)=O_{\PP}(n^{1/2})$ and, thus $D_{2,11} =O_{\PP} \big(\sqrt{\log n}\big)$.

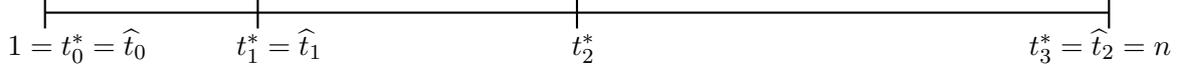
\begin{figure}
		{\centering 
		\setlength{\unitlength}{.7cm}
		\begin{center}
			\begin{picture}(25,3)(-5.5,-1.5)
			\thicklines
			\put(-5,0){\line(1,0){20}}
			\put(-5,-0.3){\line(0,1){0.6}}
			\put(-5.7,-0.8){\scalebox{1}{$1 = t_{0}^* = \widehat{t}_{0}$}}
			
			\put(-1,-0.3){\line(0,1){0.6}}
			\put(-1.4,-0.8){\scalebox{1}{$t_{1}^* = \widehat{t}_{1}$}}
			
			\put(5,-0.3){\line(0,1){0.6}}
			\put(4.9,-0.8){\scalebox{1}{$t_{2}^*$}}	
			
			\put(15,-0.3){\line(0,1){0.6}}
			\put(13.5,-0.8){\scalebox{1}{$t_{3}^* = \widehat{t}_2 = n$}}	
			\end{picture}
		\end{center}
	}	
	\caption{\footnotesize An illustration of the change-point scenario with two estimated jumps ($K^* = 2$ and $K = 1$).}
	\label{fig22}
\end{figure}

For $ D_{2,12}$, we use the following identity: for all $a >0$ (the situation for $a <0$ is quite analogous) it holds that  $\eE[\int^a_0 (\e1_{\{\varepsilon < t \}} - \e1_{\{\varepsilon < 0\}})dt ]= \int^a_0 \eE[\e1_{\{0 < \varepsilon < t\}}] dt = \int^a_0 [F(t) - F(0)]dt$.  Now,  for some $t$ in a neighborhood of zero, we can write $F(t)- F(0)= t f(\tilde t) $, for some $\tilde t \in (0, t)$, and using the fact that $f(t)>0$ for all $t \in \R$, which follows from Assumption (A2), we have that:
\[
\eE [D_{2,1} ] = \eE [D_{2,12} ]
%\sum^{K^*+1}_{k=1} \sum^{t^*_k-1}_{i=t^*_{k-1}} \bigg[ \rho_\tau \bigg( \varepsilon_i - C \sqrt{\frac{K^* %\log n}{n}} \bigg) - \rho_\tau(\varepsilon_i) \bigg] 
= C_+ \sum^{K^*+1}_{k=1}  \sum^{t^*_k}_{i=t^*_{k-1}} \frac{K^* \log n}{n} =C_+K^* \log n .
\] 
Similarly, we obtain that the variance of $ D_{2,12}$ is $C_+K^* \log n$. Thus, by the Bienaym\'e-Tchebychev   inequality for ${n}^{-1}D_{2,12}$, we have with probability converging to 1 that  $D_{2,12}  =C_+K^* \log n  $,   and  $D_{2,1}= D_{2,12} \big( 1+o_{\PP}(1)\big)=C_+K^* \log n$. 
Then, taking also into account relation (\ref{D2.2}),  Assumption (A6), we obtain
\begin{equation}
\label{D2}
D_2  = C_+^{(3)}K^* \log n +O_{\PP} \left(   (n \lambda_n) \sqrt{\frac{  \log n}{n}} \right) = C_+^{(3)}K^* \log n +O_{\PP} \left(\log n \right),
\end{equation}
with $C_+^{(3)}>0$ bounded and not depending on $n$.\\
\bigskip
Finally, we study $D_1$ of (\ref{D}). We recall that $\widehat{\cal T}_{\widehat{K}} \equiv \{  \widehat t_1,  \cdots , \widehat t_{\widehat K}\}$. Then, $D_{1,2}$ can be rewritten as $D_{1,2} = n\lambda_n \big[{\cal R}(\widehat{\cal T}_{\widehat{K}})- {\cal R}({\cal T}^*)  \big] $, with
\[
{\cal R}(\widehat{\cal T}_{\widehat{K}}) \equiv \sum^{\widehat K+1}_{k=1} | \widehat \mu_{k+1} - \widehat \mu_k|,
\qquad \textrm{and} \qquad
{\cal R}({\cal T}^*) \equiv \sum^{K^*+1}_{k=1}  | \mu^*_{k+1}- \mu^*_k |.
\]
Thus, since $\widehat \mu_k$ is bounded for all $k$, and since  $\widehat{K} < K^*$,  we have $D_{1,2}=  O_{\PP} ( K^* ( n \lambda_n)) =O_{\PP}(n \lambda_n)$.\\
We study now $D_{1,1}$. For this, let us also consider the sets   ${\cal T}^* \equiv \{ t^*_1, \cdots, t^*_{K^*} \}$, $\overset{\sim}{\cal T} \equiv   \widehat{\cal T}_{\widehat{K}} \cup  {\cal T}^*= \{ \overset{\sim}{t}_1, \cdots , \overset{\sim}{t}_{|\overset{\sim}{\cal T} |}\}$, and $\overset{\sim}{\emu}(|\overset{\sim}{\cal T} |) \equiv \argmin_{ \emu(|\overset{\sim}{\cal T} |)\in \R^{|\overset{\sim}{\cal T} |+1}} S \big( \emu(|\overset{\sim}{\cal T} |), |\overset{\sim}{\cal T} | \big)= \big( \overset{\sim}{\mu}_1, \cdots , \overset{\sim}{\mu}_{|\overset{\sim}{\cal T} |+1}  \big)$.  Then, we can write,
\begin{align}
D_{1,1} & = \bigg[ \sum^{\widehat K+1}_{k=1} \sum^{\widehat t_k}_{i=\widehat t_{k-1}+1} \rho_\tau (Y_i -\widehat \mu_k) && - \sum_{k \in \{ 1, \cdots, |\overset{\sim}{\cal T}| \}}  \sum^{\overset{\sim}{t}_k-1}_{i=\overset{\sim}{t}_{k-1}} \rho_\tau(Y_i - \widetilde \mu_k)   \bigg] \nonumber \\ 
& ~ && + \bigg[ \sum_{k \in \{ 1, \cdots, |\overset{\sim}{\cal T}|\}}   \sum^{\overset{\sim}{t}_k-1}_{i=\overset{\sim}{t}_{k-1}}\rho_\tau(Y_i - \widetilde \mu_k)  - \sum^{K^*}_{k=1} \sum^{t^*_k}_{i=t^*_{k-1}+1} \rho_\tau (\varepsilon_i)  \bigg] \nonumber\\
& \equiv D_{1,11}+D_{1,12}. && \nonumber
\end{align}
For a better illustration when dealing with $D_{1,11}$ we take a particular case for $K^*=2$ and  $K=1$ (see the illustration in Figure \ref{fig22}). The other cases are the same, but more painful to do. 
We start by expressing $D_{1, 11}$ as a sum of three terms where
\begin{align*}
D_{1,11} & = \sum^{t^*_1}_{i=t^*_0} \bigg(\rho_\tau (\varepsilon_i+\mu^*_1-\widehat \mu_1)  -   \rho_\tau (\varepsilon_i + \mu^*_1-\widetilde \mu_1) \bigg)\\
& \quad \quad + \sum^{t^*_2}_{i=t^*_1+1}\bigg( \rho_\tau (\varepsilon_i+\mu^*_2-\widehat \mu_2) 
-   \rho_\tau (\varepsilon_i+\mu^*_2-\widetilde \mu_2) \bigg)\\
& \quad \quad  +\sum^{t^*_3}_{i=t^*_2+1}\bigg( \rho_\tau (\varepsilon_i+\mu^*_3-\widehat \mu_2) -   \rho_\tau (\varepsilon_i+\mu^*_3-\widetilde \mu_2) \bigg)\\
& \equiv D_{1,111}+ D_{1,112}+ D_{1,113}.
\end{align*}
For $D_{1,111}$, since $|\widehat \mu_1 - \mu^*_1 |=  O_{\PP}(|\widetilde \mu_1 - \mu^*_1 |)= O_{\PP}\bigg( \sqrt{\frac{  \log n}{n}} \bigg) $, we have, 
\[
\eE\bigg[ \rho_\tau \big(\varepsilon_i +C_1  \sqrt{\frac{ \log n}{n}}  \big)- \rho_\tau\big(\varepsilon_i +C_2  \sqrt{\frac{  \log n}{n} } \big) \bigg]=O \bigg( \frac{ \log n}{n}\bigg).
\]
Then, same as for $D_{2,12}$, we have that $D_{1,111}=O_{\PP} \pth{I^*_{max} n^{-1} \log n } $. For $D_{1,112}$, the estimator $\widehat \mu_2$ is different from at least one of the true values $\mu^*_2$ or $\mu^*_3$. Suppose that it is different to $\mu^*_2$. Then, since $\mu^*_k $ does not depend on $n$, we have that for all $\epsilon >0$, there exists some constant $C > 0$, such that $\PP\big[ |\widehat \mu_2 - \mu^*_2 | > C \big] >1- \epsilon$.

Let us now define ${\cal D} \equiv \rho_\tau (\varepsilon_i+\mu^*_2-  \mu^*_2+C_{(+)})  -   \rho_\tau (\varepsilon_i+\mu^*_2-  \mu^*_2+ C \sqrt{\frac{\log n}{n}})$, with $C_{(+)}$ being a positive constant not depending on $n$. By Assumption (A2) we have $f(x)>0$ and the density is bounded for all $x \in \R$,  therefore
\[
\eE[{\cal D}]=\int^{C_{(+)}}_{ C \sqrt{\frac{ \log n}{n}}}\left[ F(t)- F \bigg( C \sqrt{\frac{  \log n}{n}} \bigg)  \right] dt \geq C_{(++)},
\]
with  $C_{(++)} > 0$ being some positive constant. Then, with probability converging to 1, as $n \rightarrow \infty$, we also have  $D_{1,112} \geq C_+^{(1)} I^*_{min}$, for some  $C_+^{(1)} > 0$  not depending on $n$.
   
If $\widehat \mu_3 - \overset{\sim}{\mu}_3=  O_{\PP}\pth{\sqrt{\frac{  \log n}{n}}}$, then $D_{1,113}$ is same as  $D_{1,111}$, otherwise, since $\mu^*_k $ does not depend on $n$, we have that, there exists $C>0$ such that for all $\epsilon >0$, $\PP \big[|\widehat \mu_2 - \widetilde \mu_3 | >C \big] >1- \epsilon$. Then, $D_{1,113} >0$ and it is same as $D_{1,112}$. To conclude, we have that the following holds  
\[
D_{1,11} \geq C_+^{(1)} I^*_{min} - I^*_{max} \frac{K^* \log n}{n} (K^*-1),
\]  
with probability converging to 1, as $n \rightarrow \infty$.

For  $D_{1,12}$, as for $D_{2,1}$, we have with probability converging to 1, as $n \rightarrow \infty$, that
\[
D_{1,12}= C_+^{(2)} K^* \log n .
\]
On the other hand, the relation in \eqref{A7} implies that $I^*_{min} / ( I^*_{max} \frac{ \log n}{n}) \rightarrow \infty$ for $n \rightarrow \infty$, and thus, we have with probability converging to one, that 
\[
D_{1,1} \geq C_+^{(1)} I^*_{min} +C_+^{(2)} K^* \log n.
\]
Finally, since $D_{1,2} = O_{\PP}(   ( n \lambda_n))$, and taking into account that by  Assumptions (A3) and (A4) we have $I^*_{min} / (n \lambda_n) \rightarrow \infty$ for $n \rightarrow \infty$, we obtain that
\[
D_1 \geq C_+^{(1)} I^*_{min}+C_+^{(2)} K^* \log n,
\]
holds with probability converging to one as $n$ tends to infinity. Taking now into account the expression in \eqref{D2}, we have
\[
D=D_1-D_2 \geq C_+^{(1)} I^*_{min}+ \bigg( C_+^{(2)} - C_+^{(4)}\bigg)   \log n,
\]
where the inequality holds with probability converging to 1, for $n \rightarrow \infty$,  with some $C^{(4)}_+ \geq 0$. By relation (\ref{A7}), the right side of the last relation  is dominated by $C_+^{(1)} I^*_{min}$, which is greater then zero.
Thus, we obtain  (\ref{eq28}), when $\widehat K < K^*$, which completes the proof.
\hspace*{\fill}$\blacksquare$ \\

%\end{appendices}


\begin{thebibliography}{99}
	\bibliographystyle{abbrv}
	
	\bibitem[Antoch et al.(2006)]{Antoch06}
	Antoch, J., Gregoire, G., and Hu\v{s}kov\'a M. (2006).
	  Test for Continuity of Regression Function.
	  {\it Journal for Statistical Planning and Inference}, {\bf 137}(1), 753 -- 777.
	
		 
	\bibitem[Boysen(2009)]{boysen2009}
	Boysen, L., Kempe, A., Munk, A., Liebscher, V., and Wittich, O. (2009).  
	  Consistencies and rates of conference of jump penalized least squares estimators.
	 {\textit{Annals of Statistics}}, {\bf 37}(1), 157--183.
	
	\bibitem[Chen et al.(2001)]{chen2001}
	Chen, S., Donoho, D., and Saunders, M.A. (2001).
	  Atomic decomposition by basis pursuit. 
	 {\textit{SIAM Reviews}}, {\bf 43}(1), 129--159.
	
	\bibitem[Ciuperca(2014)]{Ciuperca.14}
	Ciuperca, G. (2014).
	  Model selection by LASSO methods in a change-point model.
	 {\textit{Statistical Papers}}, {\bf 55}(1), 349--374.
	
	\bibitem[Ciuperca(2016)]{Ciuperca-16}
	Ciuperca, G. (2016). 
	  Adaptive LASSO model selection in a multiphase quantile regression. 
	 {\it Statistics}, {\bf 50}(5), 1100--1131. 
 
    \bibitem[Cs\" org\H o and Horv\'ath(1988)]{Csorgo1988}
    Cs\" org\H o, M. and Horv\'ath, L. (1988).
      20 Nonparametric methods for changepoint problems.
      {\it Handbook of Statistics}, {\bf 7}, 403 -- 425.

		\bibitem[Cs\" org\H o and Horv\'ath(1997)]{Csorgo1997}
	Cs\" org\H o, M. and Horv\'ath, L. (1997).
	  Limit Theorems in Change-Point Analysis.
	  {\it Wiley Series in Probability \& Statistics}, Chichester, England.
 
		\bibitem[Desmet and Gijbels(2011)]{desmet2011}%1
		Desmet, L. and Gijbels, I. (2011). 
		  Curve Fitting Under Jump and Peak Irregularities Using Local Linear Regression.
		   {\itshape Communications in Statistics - Theory and Methods}, {\bf 40}, 4001 -- 4020.
	
	\bibitem[Dvoretzky et al.(1956)]{Dvoretzky.Kiefer.Wolfowitz.56}
	Dvoretzky, A., Kiefer, J., and Wolfowitz, J. (1956). 
	  Asymptotic minimax character of the sample distribution function and of the classical multinomial estimator.
	 {\textit{Annals of Mathematical Statistics}}, {\bf 27}, 642--669.
	
	\bibitem[Fan et al.(2014)]{Fan.Fan.Barut.14}
	Fan, J., Fan, Y., and Barut, E. (2014). 
	  Adaptive robust variable selection.
	 {\it Annals of Statistics}, {\bf 42}(1), 324 --351.
	
	\bibitem[Frick, Munk and Sieling(2014)]{Frick.14}
	Frick, K., Munk, A., and Sieling, H. (2014).
	  Multiscale change point inference.
	  {\it Journal of the Royal Statistical Society. Series B: Statistical Methodology}, {\bf 76}(3), 495--580.
 
		\bibitem[Gao et al.(2008)]{gao08}%1
		Gao, J., Gijbels, I., and Van Bellegem, A. (1995).
		  Nonparametric Simultaneous Testing for Structural Breaks.
		  {\it Journal of Econometrics}, {\bf 143}, 123 -- 142.
	
	\bibitem[Harchaoui and Levy(2010)]{Harchaoui.Levy.10}
	Harchaoui, Z. and L\'evy-Leduc, C. (2010).
	  Multiple change-point estimation with a total variation penalty.
	  {\it Journal of the American Statistical Association}, {\bf 105}(492), 1480--1493.
 
	\bibitem[He and Shao(1996)]{he1996}
	He, X. and Shao, Q. (1996).
	  A General Bahadur Representation of M-estimators and
	Its Applications to Linear Regression with Nonstochastic Design.
	  {\it Annals of Statistics}, {\bf 24}(6), 2608 – 2630.
	
	\bibitem[Horv\'ath and Kokoszka(2002)]{Horvath2002}
	Horv\'ath, L. and Kokoszka, P. (2002).
	  Change-Point Detection With Non-Parametric Regression.
	  {\it Statistics 36}, {\bf 23}(1), 9--31.
	
	\bibitem[Hyun, G'Sell and Tibshirani(2016)]{Tibshirani.16}
	Hyun, S., G'Sell, M., and Tibshirani, R.J. (2016).
	  Exact Post-Selection Inference for Changepoint Detection and
	Other Generalized Lasso Problems.
	  {\it arxiv.org/abs/1606.03552}.
 
	\bibitem[Kim et al.(2009)]{kim2009} 
	Kim, H.J., Yu, B., and Feuer, E.J. (2009). 
	  Selecting the Number of Change-poins in Segmented Line Regression. 
	  {\it Statistica Sinica}, {\bf 19}, 597--609.
	
	\bibitem[Knight(1998)]{Knight.98}
	Knight,  K. (1998). 
	  Limiting distributions for L1 regression estimators under general conditions.
	  {\it Annals of Statistics},  {\bf 26}(2),  755--770.
	
	\bibitem[Koenker(2005)]{Koenker.05}
	Koenker, R. (2005).
	  Quantile Regression.
	  {\it Cambridge University Press}, Cambridge, United Kingdom.
	
	\bibitem[Lee et al.(2016)]{Lee.Seo.Shin.16}
	Lee, S., Seo, M.H., and Shin, Y. (2016).
	  The lasso for high dimensional regression with a possible change point.
	  \textit{Journal of the Royal Statistical Society: Series B}, \textbf{78}(1), 193--210.
 
    \bibitem[Lin et al.(2016)]{lin2016}
	Lin K., Sharpnack J., Rinaldo A., and  Tibshirani R. J.(2016).
	  Approximate Recovery in Changepoint Problems, from $l_2$ Estimation Error Rates.
	  {\it arxiv.org/abs/1606.06746}
	
	\bibitem[Maciak and Hu\v{s}kov\'a(2017)]{Maciak.Huskova.17}
	Maciak, M. and Hu\v{s}kov\'a, M. (2017).
	  Discontinuities in Robust Nonparametric Regression
	with $\boldsymbol{\alpha}$-mixing Dependence.
	  {\it Journal of Nonparametric Regression}, {\bf 29}(2), 447--475.
	
	\bibitem[Maciak and Mizera(2016)]{Maciak.Mizera.16}
	Maciak, M. and Mizera, I. (2016).
	  Regularization Techniques in Joinpoint Regression.
	  {\it Statistical Papers}, {\bf 57}(4), 939--955.
 
	\bibitem[Mammen and Van De Geer(1997)]{mammen1997}
	Mammen, E. and Van De Geer, S. (1997).
	  Locally Adaptive Regression Splines.
	 {\itshape Annals of Statistics}, {\bf 25}(1), 387--413.
 
	\bibitem[Meinshausen and B\"uhlmann(2006)]{Meinshausen2006}
	Meinshausen, N. and  B\"uhlmann, P. (2006).
	  High-dimensional graphs and variable selection with the Lasso.
	 {\itshape Annals of Statistics}, {\bf 34}(3), 1436--1462.
	
	\bibitem[Peštová and Pešta(2016)]{pestova.pesta.16}
	Pe\v{s}tov\'a, B. and Pe\v{s}ta, M. (2016).
	  Testing structural changes in panel data with small fixed panel size and bootstrap.
	  {\it Metrika}, {\bf 78}(6), 665--689.
	
	\bibitem[Qian and Su(2016)]{Qian.Su.16}
	Qian, J. and Su, L. (2016).
	  Shrinkage estimation of regression models with multiple structural changes.
	  {\it Econometric Theory}, {\bf 32}(6), 376--1433.
 
		\bibitem[Qiu and Yandell(1998)]{qiu_yandell}
		Qiu, P. and Yandell, B. (1998).
		  A Local Polynomial Jump Detection Algorithm in Nonparametric Regression.
		  {\itshape Technometrics}, {\bf 40}, 141 -- 152.
		

		\bibitem[Tibshirani(1996)]{Tibshirani:96}
Tibshirani, R. (1996). 
  Regression shrinkage and selection via the LASSO.
 {\textit{Journal of the Royal Statistical Society: Series B}}. \textbf{58}, 267--288.

		\bibitem[Tibshirani et al.(2005)]{tibs2005}
		Tibshirani, R., Saunders, M., Rosset, S., Zhu, J., and Knight, K. (2005). 
		  Sparsity and smoothness via the fused lasso.
		 {\textit{Journal of the Royal Statistical Society: Series B}}. \textbf{67}, 91--108.

\bibitem[Tibshirani(2011)]{Tibshirani:11}
Tibshirani, R. (2011). 
  Regression shrinkage and selection via the lasso: a retrospective.
 {\textit{Journal of the Royal Statistical Society: Series B}}. \textbf{73}, 273--282.

 
	\bibitem[Tibshirani (2014)]{tibs2014}
	Tibshirani, R. J. (2014).
	  Adaptive Piecewise Polynomial Estimation via Trend Filtering. 
	  {\it Annals of Statistics}, {\bf 42}(1), 285 -- 323.
	
 
	\bibitem[Zhao and Yu(2006)]{zhao2006}
	Zhao, P. and Yu, B. (2006)
	  On Model Selection Consistency of Lasso. 
	  {\it Journal of Machine Learning Research}, {\bf 6}, 2541--2567.
	
	\bibitem[Zou and Hastie (2005)]{zou2005}
	Zou, H. and Hastie, T. (2005)
	  Regularization and variable selection via the elastic net. 
	  {\it Journal of the Royal Statistical Society: Series B}, {\bf 67}(2), 301--320.
	
	\bibitem[Zou (2006)]{zou2006}
	Zou, H. (2006)
	  The Adaptive Lasso and Its Oracle Properties. 
	  {\it Journal of the American Statistical Association}, {\bf 101}, 1418--1429.
\end{thebibliography}
\end{document}